\documentclass[%
 aip,
 amsmath,amssymb,
 reprint,%
]{revtex4-1}

\usepackage{graphicx}% Include figure files
\usepackage{dcolumn}% Align table columns on decimal point
\usepackage{bm}% bold math
\usepackage{algorithm,algcompatible,amsmath}
\usepackage{svg}
\usepackage{bbm}
\usepackage[caption=false]{subfig}
\usepackage[caption=false]{subfig}
\usepackage[english]{babel}
\DeclareGraphicsExtensions{.pdf,.png}

\definecolor{edit_col}{rgb}{0, 0, 0}

\draft % marks overfull lines with a black rule on the right

\begin{document}

%\preprint{}

\title{Selecting embedding delays: An overview of embedding techniques and a new method using persistent homology.} %Title of paper

\author{Eugene Tan}
\altaffiliation[Corresponding Author, ]{eugene.tan@uwa.edu.au}
\affiliation{Complex Systems Group, Department of Mathematics and Statistics, The University of Western Australia, Crawley, Western Australia 6009, Australia}
    
\author{Shannon Algar}
\affiliation{Complex Systems Group, Department of Mathematics and Statistics, The University of Western Australia, Crawley, Western Australia 6009, Australia}
\affiliation{Forrest Research Foundation, The University of Western Australia, Crawley, Western Australia 6009, Australia}

\author{D\'{e}bora Corr\^{e}a}
\affiliation{Complex Systems Group, Department of Mathematics and Statistics, The University of Western Australia, Crawley, Western Australia 6009, Australia}
\affiliation{ARC Centre for Transforming Maintenance Through Data Science, The University of Western Australia, Crawley, Western Australia 6009, Australia}

\author{Michael Small}
\affiliation{Complex Systems Group, Department of Mathematics and Statistics, The University of Western Australia, Crawley, Western Australia 6009, Australia}
\affiliation{ARC Centre for Transforming Maintenance Through Data Science, The University of Western Australia, Crawley, Western Australia 6009, Australia}
\affiliation{Mineral Resources, CSIRO, Kensington, Western Australia 6151, Australia}

\author{Thomas Stemler}
\affiliation{Complex Systems Group, Department of Mathematics and Statistics, The University of Western Australia, Crawley, Western Australia 6009, Australia}

\author{David Walker}
\affiliation{Complex Systems Group, Department of Mathematics and Statistics, The University of Western Australia, Crawley, Western Australia 6009, Australia}

\date{\today}


\begin{abstract}
Delay embedding methods are a staple tool in the field of time series analysis and prediction. However, the selection of embedding parameters can have a big impact on the resulting analysis. This has led to the creation of a large number of methods to optimise the selection of parameters such as embedding lag. This paper aims to provide a comprehensive overview of the fundamentals of embedding theory for readers who are new to the subject. We outline a collection of existing methods for selecting embedding lag in both uniform and non-uniform delay embedding cases. Highlighting the poor dynamical explainability of existing methods of selecting non-uniform lags, we provide an alternative method of selecting embedding lags that includes a mixture of both dynamical and topological arguments. The proposed method, {\em Significant Times on Persistent Strands} (SToPS), uses persistent homology to construct a characteristic time spectrum that quantifies the relative dynamical significance of each time lag. We test our method on periodic, chaotic and fast-slow time series and find that our method performs similar to existing automated non-uniform embedding methods. Additionally, $n$-step predictors trained on embeddings constructed with SToPS was found to outperform other embedding methods when predicting fast-slow time series. 
\end{abstract}

\maketitle %\maketitle must follow title, authors, abstract and \pacs

\textbf{Embedding methods are commonly used to analyse time series whose full system state cannot be fully or directly observed. However, most embedding methods require the careful selection of parameters to achieve a faithful reconstruction of the system dynamics. One common class of embedding methods --- time delay embedding --- requires the careful selection of embedding lags. In this paper, we provide an outline of embedding theory and a collection of existing methods and principles for guiding the selection of embedding lags. Finally, we present an argument for the usage of non-uniform embedding and propose a new persistent homology based method, SToPS, to inform the selection for multiple embedding lags.}

\section{The Case for Embedding}
\label{embedding_case}
Since the significant results proposed by Whitney \cite{whitney1936differentiable} and Takens \cite{takens1981detecting}, the ideas of mathematical embedding have pervaded through almost all aspects of the nonlinear dynamics literature. The related theorems were then subsequently formalised, codified and discussed in the seminal paper ``Embedology", by Sauer, Yorke and Casdagli \cite{sauer1991embedology}. This paved the way for the development of numerous embedding techniques such as the method of derivatives \cite{packard1980geometry}, time delay embedding \cite{abarbanel1994predicting, takens1981detecting} and PCA embedding \cite{broomhead1986extracting} among others, that have been subsequently applied to a wide variety of study areas. Today, the embedding approach remains an invaluable tool in the study of nonlinear time-series analysis. 

\textcolor{edit_col}{A time series $x(t)\in \mathbb{R}^m$ may be generically viewed as the product of a data generating process consisting of successive, though not always regular, observations of some dynamical process with state $\vec{s}(t) \in \mathbb{R}^n$ via a measurement function $h: \mathbb{R}^n \to \mathbb{R}^m$. Typically, $m<n$ as the full state of the underlying dynamical process cannot be observed. For the purposes of illustration, we will consider in this paper the simplest case where $x(t)$ is scalar (i.e. $m=1$).} 

The main goals in time series analysis often fall into two main categories: system identification or classification, and prediction. The aims of the former are focused on characterising and understanding the dynamics and operating mechanisms of the underlying dynamical system. This can range from the study of system invariants such as Lyapunov spectra \cite{lyapunov1992general, wolf1985determining} and correlation dimension \cite{grassberger1983characterization, grassberger2004measuring,kantz2004nonlinear} to bifurcation analysis \cite{guckenheimer2013nonlinear, kantz2004nonlinear}. The latter task of time-series prediction has a more practical aim that is clearly stated in its name: ``given some history $x(t_0)...x(t)$, find the best predicted estimates of $x(t+\tau)$''. It is worth clarifying that these two areas are not exclusive and often work synergistically. However, we will focus on the latter problem of time series prediction in this paper. 

Though simple in its aim, practically fulfilling the task of time series prediction presents many challenges. Outside the study of simple toy models, many interesting systems exhibit high-dimensional and also chaotic behaviour \cite{gershenfeld1999nature, govindan1998evidence, weigend1993results}. For the cases of time delay systems, the dimension of the system dynamics may not even be finite. \textcolor{edit_col}{The potential inaccessibility of the full system state also adds to the challenge. As such, it is common to reframe the time series prediction problem in terms of the measurement function $h: \mathbb{R}^n \to \mathbb{R}^m$ where $m<n$.} Thus, time series prediction can be rewritten: ``given some observed history $h(x(t_0)),...,h(x(t))$, predict the future state $x(t+\tau)$''. This can be framed in terms of probability theory with aim to calculate the conditional probability $P(x(t+\tau)|h(x(t_0)),...,h(x(t)))$. 

Given the observational restrictions, it is desirable for any time-series predictor to extract and utilise as much information that is contained within the observed time series $x(t)$. Therein lies the value of embedding, vis-\`a-vis dimension augmentation. Fortunately, under sufficient precision and noise-free assumptions, Takens' embedding theorem guarantees that a time delay embedding with dynamics $\Phi$ defined in a space of sufficiently high dimension $\mathbb{R}^d$ constructed from scalar observed time series is generically diffeomorphic to the full state space dynamics of the underlying dynamical system \cite{takens1981detecting}. In essence, one can reconstruct a proxy image of the full state dynamics using only part of the observed system variable. The usual approach to achieve this is to embed the observed time series into a sufficiently large dimension containing the `full' state dynamics, learn the state space dynamics or vector field using one's favoured modelling tool (e.g. neural networks \cite{Sangiorgio2020, Zhang2000}, reservoir computing \cite{Haluszczynski2019, jaeger2001echo, Pathak2017} support vector machines \cite{muller1997predicting, sapankevych2009time} etc.) and predict forward, performing all the required calculations in the new embedded space. %This still remains a common approach in the field of time-series prediction.

Despite the elegance and utility afforded by embedding theorems, their reliance on infinite data precision, length and noise-free signals pose practical problems. The effectiveness of embedding is highly dependent on the choice of embedding parameters used to augment the data \cite{cellucci2003comparative, theiler1992testing}. For example, in the case of time delay embedding, the selection of time lag \cite{fraser1986independent} and embedding dimensions \cite{kennel2002false, sauer1991embedology} can have a profound effect on the quality of the resulting reconstruction.

Thus, we may summarise the main challenge of embedding as the following question: ``How do we select good embedding parameters?". This will be the focus of our discussion. Numerous methods have been proposed to tackle  this problem ranging from purely statistical and dynamical arguments such as mutual information \cite{fraser1986independent, kantz2004nonlinear} and continuity statistic \cite{pecora2007unified}, to more purely topological arguments like distortion \cite{casdagli1991state}, noise amplification \cite{casdagli1991state, uzal2011optimal} and fill-factor \cite{buzug1992comparison}. In many cases, each of these methods only perform well for specific types of systems and not others. 

This paper has three main objectives. Firstly, to provide a simple overview of the challenges of selecting good embedding parameters. Secondly, to collate and compare the various popular methods across the dynamics-topology spectrum that have been proposed to tackle the problem of embedding parameter selection. We will focus on the particular case of optimising time delay embedding. Finally, to present a different approach based on the growing field of persistent homology --- the significance score --- which attempts to incorporate both dynamical and topological arguments into the selection of embedding parameters. 

This paper is structured as follows. We begin in  Section \ref{section:Embedding Methods} by providing an overview of embedding theory and various common embedding methods in. This is followed by a short discussion on guiding principles on the selection of embedding parameters in Section \ref{embedding_considerations}. Sections \ref{uniform_delay_embeddings}-\ref{simultaneous} introduces several embedding parameter selection methods for both uniform and non-uniform delay embedding. Finally, we present in Section \ref{Persistent_Strands} SToPS, a persistent homology approach to embedding parameter selection, which is our contribution to the embedding parameter selection problem.

\section{Embedding Methods}
\label{section:Embedding Methods}
Multiple embedding methods exist to perform state space reconstruction. However, they all aim to perform a similar task, augment an observed time series into a high enough dimension that is useful for describing the underlying dynamics. For completeness, we provide a brief overview on the embedding and three common embedding methods, time delay embedding, derivatives embedding and global principal value embedding. However, the ideas present in this paper will focus on the selection of time delay embedding parameters. A deeper discussion on the other embedding methods can be found elsewhere \cite{casdagli1991state}.

\subsection{Embedding Theory}
For a given dynamical system with state $\vec{s}(t)$ with dynamics on state space $S \subseteq \mathbb{R}^d$ and evolution operator $f$ such that,
\begin{equation}
    \vec{s}(t+T) = f_T(\vec{s}(t)),
\end{equation}
we can define a measurement function $h: \mathbb{R}^d \to \mathbb{R}$ that simulates the process of observing the system and extracting information evaluated at given time steps to produce a time-series,
\begin{equation}
    x(t) = h(\vec{s}(t)).
\end{equation}

An embedding can be defined as a transformation $\Psi : \mathbb{R} \to \mathbb{R}^m$ that augments the dimension of the time series using observed values across some window of time $[t_1, t_2]$. For example, uniform delay embedding is given by the following:
\begin{equation}
    \Psi(x(t)) = \vec{x}(t) = (x(t),x(t-\tau),...),
\end{equation}
where $\vec{x}(t)$ is the embedding vector with dynamics defined in a reconstructed state space $X \subseteq \mathbb{R}^m$ and a transformed evolution operator $F$. If $\Psi$ is a valid embedding, Takens' embedding theorem guarantees that generically, there exists a diffeomorphism $\Phi: S \to X$ that preserves the dynamics of the system such that,
\begin{equation}
    F = \Phi \circ f \circ \Phi ^{-1}.
\end{equation}

\begin{figure*}
    \centering
    \includegraphics[width = 0.95\textwidth]{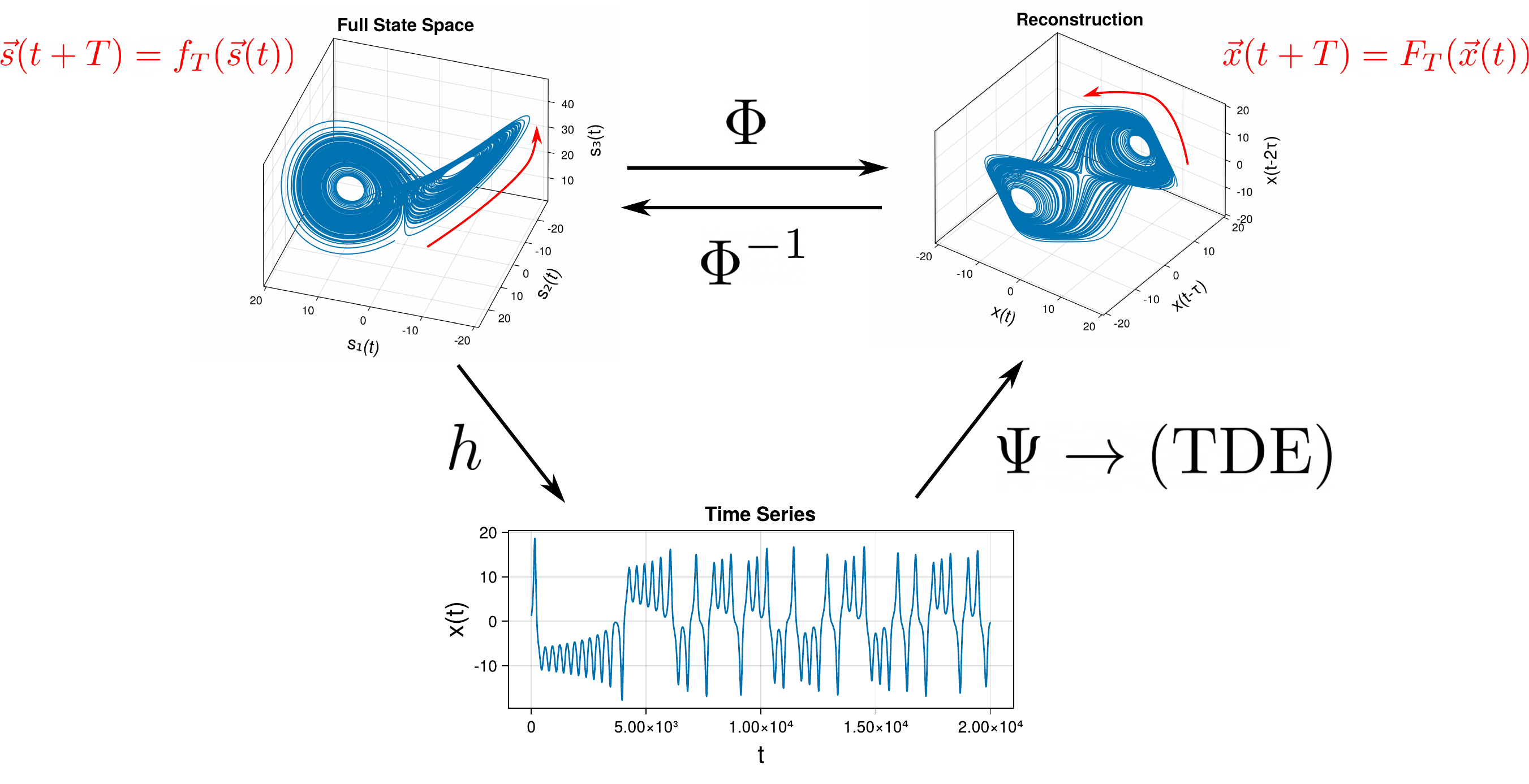}
    \caption{Schematic of the embedding process and the relationship between its components.}
    \label{fig:Embedding}
\end{figure*}

Learning the dynamics along the reconstructed state space $F$ is equivalent to learning the true system dynamics $f$ (see Figure \ref{fig:Embedding}). Therefore, the task of time series prediction using embedding simplifies to learning the evolution operator $F$ where
\begin{equation}
    x(t+T)=\Psi ^{-1}  \circ F_T \circ \Psi(x(t)).
\end{equation}

\subsection{Time Delay Embedding}
First described by Packard et al. \cite{packard1980geometry}, time delay embedding involves the augmentation of a scalar time series $x(t)$ into a higher dimension through the construction of delay vector $\vec{x}(t)$ given as
\begin{equation}
    \vec{x}(t) = (x(t), x(t-\tau),...,x(t-(m-1)\tau)),
\end{equation}
where the embedding parameters to be selected are the delay lag $\tau$ and embedding dimension $m$. According to the guarantees of Takens' theorem, any value of $\tau$ will yield a valid embedding given sufficiently large $m$ and measurement values of infinite precision. However, this is not achievable in practice and different selections of delay lag and embedding dimension can yield varying results. A further discussion of this is given in Section \ref{embedding_quality}. 

\textcolor{edit_col}{We also note that the task of selecting ideal delay lag $\tau$ and embedding dimension $m$ is not unique to time delay embedding. Selecting values of $\tau$ and $m$ are also key decisions in other time series analysis methods such as permutation entropy \cite{bandt2002permutation} and ordinal partition networks \cite{mccullough2015time, zhang2017constructing}. In both of these instances, a delay vector is constructed and represented by an encoding based on the size order each component. The time series may then be viewed as a transitions between different encoding states and used for further analysis. }

Many automated time-series prediction methods such as recurrent neural networks \cite{Sangiorgio2020,Zhang2000} and reservoir computing \cite{jaeger2001echo, Pathak2017, Haluszczynski2019, schrauwen2007overview} may also be related to delay embedding. In both cases, input time series is fed into a dynamical network that contains some notion of memory. The forward propagation of this memory of past states on future states effectively acts as a time delay embedding with small delay lag $\tau$ and large embedding dimension $m$. 

\subsection{Derivatives Embedding}

The embedding method of derivatives reconstructs an embedding vector using successively increasing order of time derivatives from the observed time series \cite{packard1980geometry}. This is given by:
\begin{equation}
    \vec{x}(t) = \left( x(t),\frac{dx(t)}{dt},...,\frac{d^mx(t)}{dt^m} \right)
    \label{eq:derivatives_embedding}
\end{equation}
Derivatives are taken via numerical approximations. The derivatives embedding method is a valid embedding for sufficiently large $m$ if one is able to accurately calculate the required derivatives. 

\subsection{Integral-Differential Embedding}
One weakness of the derivatives embedding approach is the need to evaluate numerical derivatives from data. Whilst this may be acceptable for the first derivative, approximations of successive higher order derivatives are generally inaccurate as the signal to noise ratio tends to be negatively impacted. This is true even for the cases of very clean data sets.

An alternative to derivatives embedding, is integral-differential embedding \cite{gilmore1998topological}. This approach avoids the calculation of successive higher order derivatives by replacing the second order the derivative with an integral instead. This yields the following embedding construction:
\begin{equation}
    \vec{x}(t) = \left( \int_{-\infty}^t x(t)-\langle x(t)\rangle_t \, dt,x(t),\frac{dx(t)}{dt} \right),
    \label{eq:int_diff_embedding}
\end{equation}
where the first component is first set to zero mean before integration. The usage of a first order integral and numerical derivative results in a degradation of the signal to noise ratio by only one order each for the first and third embedded components. This is in contrast with the derivatives embedding approach where each successive numerical derivative has a signal to noise ratio that is degraded with increasing orders of magnitude.  However, the integral-differential embedding approach suffers from the same noise effects as the pure derivatives method for higher dimensional embedding. This limits its applicability to systems where system dynamics are presumed to be high dimensional.

\subsection{Global Principal Value Embedding}
The method of principal value embedding was proposed by Broomhead and King \cite{broomhead1986extracting} as a modified alternative to time delay embedding using the theorem's by Takens. This method draws upon the ideas of principal component analysis \cite{lever2017points} to find an ideal rotation of the time delay embedding with a sufficiently high dimension. Given a time series $x(t)$ of length $N_T$ and a sliding window of length $M$, we can construct a collection of $N=N_T-(M-1)$ delay vectors,
\begin{equation}
    \mathbf{X}=N^{-1/2} 
    \begin{bmatrix}
    \vec{x}_1\\
    \vec{x}_2\\
    \vdots\\
    \vec{x}_N\\
    \end{bmatrix}
\end{equation}
where $\vec{x}_i$ is the delay vector constructing using the $i^{th}$ value in the time series as the first component,
\begin{equation}
    \vec{x}_i = (x(t_i), x(t_{i-1}),...,x(t_{i-(M-1)})).
\end{equation}

An $m\times m$ covariance matrix $C$ can be calculated from $\mathbf{X}$. The elements $C_{ij}$ of this matrix can be simply given as,
\begin{equation}
    C_{ij} = \langle x(t) \,x(t+(i-j))\rangle_t
\end{equation}
where $\langle ... \rangle_t$ denotes a time average. The principal components of $C$ are then found by calculating its respective eigenvalues and eigenvectors. Taking the first $d$ principal components corresponding to the desired number of embedding dimension, the eigenvector matrix can be used to calculate a projection of $\mathbf{X}$ corresponding to the final embedded coordinates. Readers are advised to refer elsewhere \cite{broomhead1986extracting,casdagli1991state} for more details.

Principal component value embedding essentially aims to distill and simplify a high dimensional delay embedding (usually obtained by taking a large number of lagged components) into a lower dimensional subspace. The remaining subspaces are argued correspond to component directions with little dynamical variation and importance. One application of this method was as an attempt to simplify the selection of the optimal embedding dimension \cite{broomhead1986extracting}, where the ideal embedding dimension $m$ corresponds to the number of singular values that are distinctly greater than some `noise floor'. However, this approach has received several criticisms \cite{mees1987singular,paluvs1992singular}. The main of which arguing that the onset of a plateau noise floor can be attributed to the precision and noise strength in the data, rather than the importance of the corresponding eigenvector direction.

Within the general context of embedding, Palu{\v{s}} and Dv{\v{o}}r{\'{a}}k test the quality of a delay embedding with reduced dimension obtained using the first $k$ principal eigenvector directions \cite{paluvs1992singular}. The authors show that the reduced dimension embeddings' estimates of dynamical invariants such as the correlation dimension vary with different time delay and number of components. They argue that the usage of principal components of the covariance matrix is restricted to linear correlations. Therefore, whilst components in the embedding may be independent in the linear sense, they may still be nonlinearly dependent. Instead, the truncation of embedding dimensions can result in the exclusion of important nonlinear components.

The inclusion of large time lags within a given principal component direction also may not make much sense for chaotic systems where temporally close observations decorrelate exponentially in time. From the perspective of selecting time lags, each principal component will almost invariably contain contributions from all $M$ possible lagged components. Apart from the dimensional reduction argument \cite{paluvs1992singular} (of which care must be taken in its interpretation), global principal value embedding does not present any significant difference to the general time delay embedding with large embedding dimension $m$.

\section{Embedding Considerations}
\label{embedding_considerations}
\subsection{Embedding Quality}
\label{embedding_quality}
As previously discussed, the theoretical guarantees of Takens' fail in the presence of finite precision and noise \cite{cellucci2003comparative, fraser1986independent} leading to the concept of `optimal' embedding parameters. The existence of such an `optimal' set implies that not all embeddings are of equal quality. However, this requires a measure by which embedding quality can be compared against. An attempt to quantify embedding quality was studied extensively by Casdagli \cite{casdagli1991state} and Potapov \cite{potapov1997distortions}.

There are large variations in the definition of embedding quality such as those based on information theoretic arguments \cite{fraser1986independent}, prediction tasks \cite{casdagli1991state} and attractor topology \cite{buzug1992comparison, nichkawde2013optimal}. It is worthwhile to note that the `optimality' of a set of parameters is dependent on the task that the embedding is being used for. As such, whilst one may find similar results between methods based on different notions of embedding quality, disagreement between results will likely always be present \cite{kantz2004nonlinear, casdagli1991state, cellucci2003comparative}. Therefore, it is better to avoid the claim that a particular set of parameters are more favourable unless there are dynamical and topological reasons within the data itself that support it. However, we will highlight in this section the general considerations that are often used in defining the quality of an embedding. 

The different notions of embedding quality can be summarised in two broad categories or arguments, prediction-based and topological arguments. Prediction-based notions of embedding quality can be seen to be inspired by application of embeddings in the context of time-series prediction. Fundamentally, good embeddings should enable better predictions \cite{casdagli1991state, potapov1997distortions}. 

In time series prediction tasks, the presence of an unknown measurement function $h$ and noisy data introduces some degree of uncertainty to the inference of the real system state $s(t)$. Casdagli argued that an ideal embedding should minimise the uncertainty of inferring the true state $\vec{s}(t)$ given a position $\vec{x}(t)$ in the reconstructed state space. In essence, the inverse transformation $\Phi^{-1}$ applied to a constructed delay embedding in the presence of noise should have little ambiguity on the true state of the system, if $\vec{s}(t)$ could be fully known \cite{casdagli1991state}. This robustness to noise and low ambiguity should in theory be beneficial for time series prediction and also forms the basis of the ideas of noise amplification and distortion used to quantify embedding quality. In poor embeddings, such as those whose attractor manifolds are laminar with little separation, the effect of noisy perturbations across manifold layers result in significant uncertainty of the true state $s(t)$, making time series prediction difficult.

The information theoretic arguments for choosing embedding parameters (e.g. autocorrelation, minimum mutual information, continuity statistics) are also closely related with the ideas of prediction. These methods generally try to maximise the amount of new information incorporated in each delay dimension with the aim that it will provide more information of the true system state $\vec{s}(t)$ and aid in time series prediction.

The other broad category of defining embedding quality are those based on topological and geometrical arguments. Many of these methods focus on the study of the attractor structure and distribution of the manifold in its ambient state space. In essence, a good embedding with respect to topology and geometry should aim to be well spaced out and unfolded in its ambient space \cite{buzug1992comparison, kantz2004nonlinear, nichkawde2013optimal}. This notion of quality has parallels with the noise amplification arguments of Casdagli. Some methods based on geometrical arguments include statistics such as the fill factor \cite{buzug1992comparison} and displacement from diagonal \cite{nichkawde2013optimal}. Ultimately, many of the considerations outlined above for determining the ideal lag and embedding dimension for time delay embedding can be summarised with the concepts of irrelevance and redundancy \cite{casdagli1991state, gibson1992analytic}.

\subsection{Irrelevance and Redundancy}
\label{Irrelevance_Redundancy}
The selection of time lag $\tau$ and embedding dimension $m$ are the main challenges when constructing a time delay embedding. There is uncertainty on the relative importance between embedding lag and dimension. Furthermore, it has also been proposed that these embedding parameters may not be independent. Instead the quantity $\tau_w=m\tau$ termed the embedding window has been proposed as a more important parameter to optimise \cite{casdagli1991state, gibson1992analytic}. However, for sufficiently large embedding dimension, it could be argued that a selection of $\tau_w$ may be simplified to an appropriate selection of $\tau$.

\textcolor{edit_col}{The selection of the embedding window $\tau_w$ (and by extension embedding lag $\tau$) may be summarised by a notion that it must be neither too short (redundance) nor too long (irrelevance). This explanation applies for chaotic or aperiodic signals. However, periodic signals may be successfully embedded with large lags $\tau$ where the effective lag $\tau^*$ is related to the period $T$,
\begin{equation}
    \tau = nT + \tau^*, \quad n \in \mathbb{N}.
\end{equation}
}

Embeddings with high redundance result in trajectories that lie in layers roughly parallel to each other (e.g. close to the diagonal). In the presence of sufficient noise, the clear separation between layers is affected. This results in a greater degree of uncertainty of the true system state $\vec{s}(t)$, given some noisy observation $\vec{x}(t)$ in reconstructed state space.

Similarly, embeddings with high irrelevance contain components that are highly decorrelated with the true state \cite{casdagli1991state}. This is also unfavourable as it may introduce unwanted crossings between trajectories in the reconstructed manifold. Therein lies the Goldilocks problem of selecting an embedding window $\tau_w$ that is neither too large or too small.

\section{Uniform Delay Embeddings}
\label{uniform_delay_embeddings}
The simplest form of delay embedding is the case of uniform delays where single constant values for $\tau$ and $m$ are selected. In this case, embedding vectors are selected with uniformly increasing time delays as given in Equation \ref{eq:uniform_embedding},

\begin{equation}
    \vec{x}(t) = (x(t),x(t-\tau),...,x(t-(m-1)\tau)).
    \label{eq:uniform_embedding}
\end{equation}

Due to the debate between the selection priority of $\tau$ and $m$, multiple methods have been proposed to simultaneously estimate both parameters. Some methods include those of Gao and Zheng \cite{gao1993local, gao1994}, characteristic lengths \cite{cellucci2003comparative} and Schuster (wavering product) \cite{liebert1991}. An overview of these methods are provided in Section \ref{simultaneous}.

Other common methods attempt to simplify the problem by assuming the independence of $\tau$ and $m$ and choose to estimate both values separately. Generally, an embedding $\tau$ is first determined using a choice of various measures. Once selected, uniform delay vectors of increasing dimensions are constructed and tested with algorithms such as the Grassberger-Procaccia \cite{grassberger2004measuring} algorithm or False Nearest Neighbours \cite{hegger1999,kennel1992,kennel2002false} until convergence is achieved. The length of delay vector when stability is reached is used to decide the embedding dimension $m$.

The methods that are used for determining embedding lag $\tau$ vary from simple heuristics to more complex statistical arguments. One common heuristic is the selection of embedding lag as one quarter of the signal period (or quasi-period for chaotic signals). Delving into more statistically-grounded arguments, autocorrelation \cite{kantz2004nonlinear} and its nonlinear generalisation, mutual information \cite{garcia2005nearest}, continuity statistics \cite{pecora2007unified} and L-statistics \cite{uzal2011optimal} are occasionally used to determine good values for $\tau$. A comprehensive overview and further discussion on these methods is provided in Sections \ref{optimising}-\ref{non-uniform_embedding}.

\section{Simultaneous Optimisation of Uniform Embedding Parameters}
\label{simultaneous}
In contrast with many of the current methods that involve the selection of $\tau$ and $m$ independently, several embedding methods have been proposed that simultaneously estimate both values using a single measure. This measure is often calculated across multiple lags and repeated for increasing $m$. A value for $m$ is first selected according to some criterion followed by the selection of $\tau$. Detailed steps on the implementation of these methods are outlined by Celluci  \cite{cellucci2003comparative}.

\subsection{Method of Gao and Zheng and Characteristic Lengths}
The embedding method proposed by Gao and Zheng is based on the incidence of false nearest neighbors \cite{gao1993local, gao1994}. False nearest neighbours can be attributed to either redundancy (insufficiently unfolded) and irrelevancy (spurious intersections in the attractor). The method proposed by Gao and Zheng operates on the notion that the separation distance and proportion of false nearest neighbours, should be minimised in an ideal embedding.

Consider a pair of points in embedded space $\vec{x}_i, \vec{x}_j$ and their evolution $k$ steps into the future $\vec{x}_{i+k}, \vec{x}_{j+k}$. Points that are false nearest neighbours will tend to separate faster than real neighbours as the attractor unfolds in a time delay embedding. As a result, the ratio between their distances $\lvert \vec{x}_{i+k}- \vec{x}_{j+k}\rvert / \lvert \vec{x}_i- \vec{x}_j\rvert$ will be larger for pairs of false nearest neighbours and approximately equal to 1 for real neighbours. Gao and Zheng then propose the following measure $\Lambda$ to optimise the embedding parameters,
\begin{equation}
    \Lambda(k,m,\tau) = \frac{1}{N_{ref}} \sum_{i,j} \ln \frac{\lvert \vec{x}_{i+k}- \vec{x}_{j+k} \rvert}{\lvert \vec{x}_i- \vec{x}_j \rvert},
\end{equation}
where $N_{ref}$ is the number of randomly sampled point pairs over which the distance ratio is averaged. There are several additional restrictions on the selection of point pairs $\vec{x}_{i}, \vec{x}_{j}$. Firstly, the initial separation of these points should satisfy $\lvert \vec{x}_i - \vec{x}_j \rvert \leq r$ where $r$ is a small selected threshold, i.e. the initial separation of points should be small enough such that the calculation of growing separation is sensible. Secondly, the selection of pairs of points should not have an intersecting Theiler window $\lvert i-j \rvert>l_{\rm Theiler}$, where $l_{\rm Theiler}\in \mathbb{N}^+$. This is done to prevented unwanted correlations between points on the same local trajectory \cite{theiler1986spurious, theiler1990estimating}. Finally, the constant $k$ should not be too large and selected with respect to the natural time scale of the system dynamics.

To identify good embedding parameters, profiles of $\Lambda(\tau)$ are calculated for increasing values of embedding dimensions $m$. The value of $m$ that corresponds to the largest decrease across the profile $\Lambda(\tau)$ is selected as the embedding dimension. The embedding lag $\tau$ is then selected as the first minimum of $\Lambda(\tau)$.

\subsection{Characteristic Lengths}

The method of characteristic lengths is an extension of Gao and Zheng \cite{cellucci2003comparative} that attempts to solve the problem of selecting an evolution time $k$. Instead of arbitrarily selecting $k$, a characteristic length $J(m, \tau)$ describing the natural spatial scale of the system attractor is calculated,
\begin{equation}
    J(m,\tau) = \langle \lvert \vec{x}_i - \vec{x}_j \rvert \rangle,
\end{equation}
where $\langle ... \rangle$ denotes an average over sampled pairs of points of the attractor. The characteristic length is then used to calculate the separation time $T_J(\vec{x_i},\vec{x_j})$ defined as the time taken for pairs of nearest neighbours to diverge by some proportion of the characteristic length $J(m, \tau)$. For real neighbours, $T_J$ will converge to a value related to the Lyapunov exponent of the system with increasing embedding dimension $m$, whilst false nearest neighbours will result in a smaller value $T_J$ as trajectories quickly separate. The new measure that is used to determine the embedding parameters is given by,
\begin{equation}
    C(m, \tau) = \frac{1}{N_{ref}} \sum_{i,j}  T_J(\vec{x_i},\vec{x_j}),
\end{equation}
where $N_{ref}$ is the number of sampled pairs of nearby neighbours. The values for $m$ and $\tau$ that maximise $C(m,\tau)$ are selected as the embedding dimension and lag.

\subsection{Wavering Product}
The wavering product \cite{liebert1991} is similar to that of Gao and Zheng and characteristic lengths in that all are based around the concepts of nearest neighbours. The authors propose that good embeddings should preserve the correspondence between the order of nearby neighbours of a given reference point (i.e. the order of neighbours sorted according to distance from some reference point $\vec{x}_i$ should be preserved). This is done by comparing the order of $N$ nearby neighbours of a point $\vec{x}_i$ between a given embedding $\Phi_k$, whose ordered sequence neighbours are given by,
\begin{equation}
    X_{\Phi_k}=\{ \vec{x}_{i,1},...,\vec{x}_{i,N} \}_{\Phi_{k}},
\end{equation}
and its projection onto its next order embedding $\Phi_{k+1}$ (by increasing $m$ or $\tau$) with the sequence given by,
\begin{equation}
    Z_{\Phi_{k+1}}=\{ \vec{z}_{i,1},...,\vec{z}_{i,N} \}_{\Phi_{k+1}}.
\end{equation}
Here, $\vec{x}_{i,n}$ corresponds to the $n^{th}$ nearest neighbour of the $i^{th}$ reference point $\vec{x}_i$. The projection $\vec{z}_{i,n}$ corresponds to the same neighbour data point $\vec{x}_{i,n}$ whose position is recalculated from the next order embedding $\Phi_{k+1}$.

Similarly, comparisons can also be made into a projection into an embedding of lower order (by decreasing $m$ or $\tau$) giving a new set of ordered points,
\begin{equation}
    V_{\Phi_{k-1}}=\{ \vec{v}_{i,1},...,\vec{v}_{i,N} \}_{\Phi_{k-1}}.
\end{equation}

Ideally, a good embedding should preserve a one to one correspondence in these ordered sequences. This will yield a value equal to 1 for the following ratios,
\begin{equation}
    \frac{\lvert \vec{x}_i-\vec{z}_{i,n} \rvert}{\lvert \vec{x}_i-\vec{x}_{i,n} \rvert}
,\;{\rm and}\;
\frac{\lvert \vec{x}_i-\vec{x}_{i,n} \rvert}{\lvert \vec{x}_i-\vec{v}_{i,n} \rvert}.
\end{equation}

The method presented by Schuster and Liebert propose the following measure as the product of the above two ratios,

\begin{equation}
    W_i(m\tau)=\prod_{n=1}^{N}
\left\{
\left(
\frac{\lvert \vec{x}_i-\vec{z}_{i,n} \rvert}{\lvert \vec{x}_i-\vec{x}_{i,n} \rvert}
\right)
\cdot
\left(
\frac{\lvert \vec{x}_i-\vec{x}_{i,n} \rvert}{\lvert \vec{x}_i-\vec{v}_{i,n} \rvert}
\right)
\right\}.
\end{equation}

The measure to be optimised is given by the average over $N_{ref}$ randomly sampled reference reference points,
\begin{equation}
    W(m,\tau)=\ln \left( 
\frac{1}{N_{ref}}\sum_{i=1}^{N_{ref}} W_i(m,\tau)
\right)
\end{equation}
with $m$ being selected as the dimension which achieves the limiting behaviour of $W$ and $\tau$ corresponds to the first minimum of the resulting profile. 

\section{Selecting Uniform Embedding Lag}
\label{optimising}
In contrast to the embedding methods in Section \ref{simultaneous} that utilise a single measure to optimise the selection of embedding lag and dimension, the most practiced approach still selects the embedding lag $\tau$ and dimension $m$ independently according to separate metrics. Here, we will focus on the methods used to select the lag for time delay embedding. This lag may then be used to construct delay embeddings of arbitrary dimension.

The simpler case of uniform delay embedding requires the selection of a singular value of $\tau$ which is increased in multiples to construct the required delay vector. The methods used to inform the selection of delay lag $\tau$ mirror the broad dichotomy in notions of embedding quality outlined in Section \ref{embedding_quality}. In this section, we divide the various methods for optimising embedding lag into two categories: the first based on dynamical and information theoretic arguments, and the second based on topological arguments. We also discuss the topic of non-uniform embedding in Section \ref{non-uniform_embedding}.

Methods rooted in a dynamical approach can be interpreted as focusing on the mechanism behind the data generating process and statistical relationships between measurements. Information theoretic and statistical approaches are also included in this category. The methods that we will review in this category include the autocorrelation, minimum mutual information and quarter period. This is in contrast with those in the second category whose methods are more topologically. Some examples include the fill-factor \cite{buzug1992comparison} and noise amplification \cite{casdagli1991state, uzal2011optimal}.

A brief note on embedding dimension: There are several methods that are used to identify the required embedding dimension such as the false nearest neighbours (FNN) \cite{hegger1999, kennel2002false} and the Grassberger-Procaccia algorithm \cite{grassberger2004measuring} used to estimate the correlation sum and dimension. Other invariants similar to correlation dimension such as the Kaplan-Yorke dimension \cite{kaplan1979chaotic} and box-counting dimension \cite{farmer1983dimension} are also often used to determine the embedding dimension. These results are usually used in conjunction with Whitney's theorem stating that any $d$-dimensional manifold (such as an attractor) can be embedded in at least $m>2d+1$ dimensions \cite{whitney1936differentiable}.  However, it has been noted that this direct application has its flaws as Whitney's theorem is only proved for integer dimensions $d$ , which is rarely the case for a majority of systems of interest, such as those exhibiting fractal and chaotic behaviour. However, we note that an extension of Whitney's theorem to generalise the inequality to the box counting dimensions ($m>2D_F$) was given by Sauer et al \cite{kantz2004nonlinear, ott1994coping}.

\subsection{Dynamical Approaches}
\subsubsection{Autocorrelation and Minimum Mutual Information}
Many commonly analysed dynamical systems tend to exhibit chaotic behaviour where nearby trajectories rapidly diverge and quickly become uncorrelated. The method of autocorrelation is based on the idea that each component in a reconstruction should include as much new information regarding the true state as possible. It has been suggested that components in the delay vector should aim to minimise the correlation \cite{kantz2004nonlinear}. This has similar effects to minimising the redundancy of the reconstruction.

For delay embedding, the lag corresponding to the first minimum of the autocorrelation function is taken as the embedding time lag. Alternatively, the decay time to $1/e$ of the autocorrelation signal has also been proposed \cite{kantz2004nonlinear}. A variation of this approach based on the first root of the mean local autocovariance has also been proposed as a robust alternative to the minimum autocorrelation approach \cite{moore2020mean}.

\textcolor{edit_col}{One weakness of autocorrelation is its inability to account for non-stationarities in the time series (e.g. drifts in phase, frequency and magnitude). Additionally, its application is limited to linear signals \cite{abarbanel1993analysis}. In all but the simplest cases, dynamical systems exhibit some level of non-stationary and nonlinear behaviour. Fraser and Swinney \cite{fraser1986independent} proposed that the mutual information between the system be used in place of autocorrelation. In their original paper, Fraser and Swinney first provide a geometrical interpretation to complement the theoretic arguments for mutual information. Namely, consider a set of points whose values in one component $x$ lie within some fixed window. From this set, track their positions $\tau$ steps into the future and calculate the distribution of values $p_\tau(x)$ in the same component for the same set of points. A value of $\tau$ that results in a wider distribution $p_\tau(x)$ should correspond to a good lag, which also corresponds to small values in the mutual information.}

The mutual information $I(\tau)$ can be interpreted as the nonlinear analogue of the autocorrelation function,
\begin{equation*}
    I(\tau) = \int P(x(t),x(t+\tau)) \log_2 \left(
\frac{P(x(t),x(t+\tau))}{P(x(t)) P(x(t)+\tau)} \right) dt,
\end{equation*}
where $P(x(t))$ is the probability of observing a state $x(t)$ at any given time and $P(x(t),x(t+\tau))$ is the joint probability defined similarly for both time $t$ and a future time $t+\tau$. Drawing from information theory, mutual information $I(\tau)$ aims to quantify the amount of information about a future state at time $t+\tau$ that is contained in an observation at time $t$. High levels of mutual information for a given lag $\tau$ imply a high degree of correlation between states and will result in higher redundancy for the delay reconstruction. 

The strengths of the minimum mutual information and autocorrelation lies in its ability to provide reasonable estimates for lag with relatively simple and quick computation. However, there are no guarantees for the existence of a clear minimum for a given mutual information profile $I(\tau)$ \cite{wallot2018calculation}. Additionally, calculating mutual information requires the numerical estimation of probability density functions $P(x(t))$ and $P(x(t), x(t+\tau))$, and thus requires consideration regarding optimal histogram bin size and data length requirements \cite{kraskov2004estimating, papana2009evaluation}. \textcolor{edit_col}{Numerous alternative methods to more effectively estimate mutual information have been proposed including the usage of adaptive binning \cite{darbellay1999estimation,fraser1986independent}, kernel density estimators \cite{moon1995estimation} and $k$-nearest neighbours \cite{kraskov2004estimating}.}

\subsubsection{Quarter of Period}
A commonly used heuristic for selecting an embedding lag is to set $\tau$ to be quarter of the most dominant period in the signal \cite{judd1998embedding}. This approach allows the natural time scale of the system dynamics to be encoded within the embedding procedure. This heuristic is inspired from the problem of embedding a sine wave in 2D $x(t)=\sin (\omega t)$. In this case, $\tau= \frac{2\pi}{4\omega}$ yields a 2D delay embedding that is the most circular with other values resulting in elliptical trajectories instead. However, this heuristic cannot be directly applied to chaotic systems where signals are aperiodic. Instead, an estimation of some form of pseudo-period is required, which will be the focus of our proposed method in this paper.

\subsection{Topological Approaches}
\subsubsection{Fill-Factor}
The fill-factor approach first proposed by Buzug et al. \cite{buzug1992comparison} is an entirely geometrical approach to calculating the quality of a given embedding. This method assumes that an ideal embedding should be able to unfold an attractor and maximise the separation between the trajectories. The authors argue that such an embedding optimally utilises the ambient space and reduces the ambiguity of the true state of the system for different points in the reconstructed state space.

The fill-factor is calculated by first sampling $m+1$ random points from an $m$ dimensional delay embedding of the data. A reference point $\vec{x}_r$ is then selected from this collection and the corresponding relative distance vectors can be calculated,

\begin{equation}
    \vec{d}_i (\tau)=
\begin{bmatrix}
x_i(t)-x_r(t)\\
x_i(t-\tau)-x_r(t-\tau)\\
\vdots\\
x_i(t-(m-1)\tau)-x_r(t-(m-1)\tau)\\
\end{bmatrix}.
\end{equation}

The corresponding $m \times m$ matrix can then be expressed as

\begin{equation}
    M(\tau) = (\vec{d}_1, \vec{d}_2,...,\vec{d}_m),
\end{equation}
and the volume of the resulting parallelepiped is given by calculating the determinant of $M$,

\begin{equation}
    V(\tau) = \det(M(\tau)).
\end{equation}

The final expression for the fill-factor is given by calculating the average volume over a collection of randomly sampled parallelepipeds $V_i(\tau)$, normalised by the range of the sampled data points,

\begin{equation}
    f = \log \left( \frac{\frac{1}{N}\sum_{i=1}^N V_i(\tau)}{(\max_{k}x(t_k)-\min_kx(t_k))^m}\right).
\end{equation}

The authors recommend the selection of $\tau$ that maximises the fill-factor $f$ over the interval $\tau \in (0,T_c/2)$, where $T_c$ is the characteristic recurrence time. The value of $T_c$ is given by,

\begin{equation}
    T_c=\frac{1}{\omega_c},
\end{equation}
where $\omega_c$ is the most dominant frequency from the power spectrum of the time series. 

\subsubsection{Noise Amplification}
Noise amplification was a measure proposed by Casdagli in an attempt to quantify the quality of an embedding \cite{casdagli1991state, potapov1997distortions}. This is supported by the notion that a good embedding should be useful in performing predictions. Additionally, good embeddings should be able to still perform relatively well even in the presence of noise. Noise amplification for a given embedding $\vec{x}(t)$ is defined with respect to predictability of the system $T$ steps into future under the presence of noise. Generally, this is given by:

\begin{equation}
    \sigma(T,\vec{x}) = \lim_{\epsilon \to 0} \sigma_\epsilon (T,\vec{x}),
\end{equation}
where
\begin{equation}
    \sigma_\epsilon (T,\vec{x}) = \frac{1}{\epsilon}\sqrt{\text{Var}[x(T)|B_\epsilon(\vec{x})]}.
\end{equation}

Here, $\text{Var}[x(T)|B_\epsilon(\vec{x})]$ corresponds to the conditional variance of $T$ step predictions into the future in $\mathbb{R}$ from an initial condition $\vec{x}$ in embedding space $\mathbb{R}^m$ contaminated with added small observation noise $\epsilon$. In this case, it is assumed that predictions have no model errors. This condition may be fulfilled by choosing nearby neighbours in the embedding $\mathbb{R}^m$ as a proxy for noisy initial conditions \cite{uzal2011optimal}. 

Finally, the noise amplification quantity $\sigma(T,\vec{x})$ is averaged over a collection of reference points sampled across the time series in order to calculate the noise amplification value $\sigma$. Embeddings with high noise amplification imply that nearby neighbours in embedded space $\mathbb{R}^m$ tend to have future trajectories that rapidly diverge because they do not correspond to real neighbors in the true manifold $\mathcal{M}$ state space. Therefore, the impact of noise is greatly amplified as small perturbations in the reconstructed space $\mathbb{R}^m$ result in large uncertainties in the true state of the system.

\subsubsection{L-Statistic}
One weakness of the noise amplification measure is its requirement to define $T$, the prediction horizon over which to calculate the noise amplification. This was addressed by Uzal \cite{uzal2011optimal} by modifying the definition of noise amplification to the following:

\begin{equation}
    \sigma_\epsilon^2 (\vec{x}) = \frac{1}{T_M}\int_0^{T_M} \sigma_\epsilon^2 (T,\vec{x}) \,dT.
\end{equation}

This definition calculates the noise amplification with respect to a range of prediction horizons up to a maximum value of $T_M$ and is found to be relatively robust for sufficiently large $T_M$. 

The algorithm used to calculate $\sigma$ relies on using $k$ nearest neighbours from a reference point $\vec{x}_i$ as a proxy. Based on the distribution of points, this can result in effective noise levels $\epsilon$ of different sizes for each point. Therefore, Uzal proposed a normalisation constant $\alpha_k$ accommodate for this variation given by:

\begin{equation}
    \alpha_k^2 = \left[ \sum_i \epsilon_k^{-2}(\vec{x}_i)\right]^{-1}
\end{equation}

Combining these two ideas, the authors propose that noise amplification $\sigma$ measures some notion of redundancy, and $\alpha_k$ measures some notion of irrelevance. The L-statistic is then described as a cost function to minimise both of these values simultaneously,

\begin{equation}
    L = \log(\sigma \alpha_k)=\log(\sigma)+\log(\alpha_k)
\end{equation}

\section{Non-Uniform Embedding}
\label{non-uniform_embedding}

The popularity of uniform delay embedding can be attributed to its ease of implementation and optimisation. In a direct application, uniform delays only require the selection of two parameters, $\tau$ and $m$. However, the convenience of such an approach comes at the cost of reduced versatility and limitations, particularly when analysing systems with dynamics occurring on multiple disparate timescales \cite{judd1998embedding,hirata2006reconstructing}.

Firstly, the choice to use a single delay limits the ability for the reconstruction to highlight features across multiple disparate time-scales \cite{pecora2007unified}. For example, a fast-slow system with characteristic time scales $\tau_1$ and $\tau_2$ where $\tau_1/ \tau_2 \gg1$ , the choice of selecting $\tau_1$ (i.e. slow dynamics) as the embedding lag can limit the reconstruction's ability to fully unfold attractor topologies corresponding to the fast dynamics. The dynamics the time scale of $\tau_2$ (i.e. fast dynamics) will appear as noisy fluctuations within the reconstructed state space.

Secondly, reconstruction from a uniform delay embedding that is sufficient is not necessarily optimal. Here, we must clarify that the definition of optimal presumes some criterion or notion of quality. Casdagli noted that the quality of an embedding, defined as the reconstruction's robustness to noisy data for prediction, can vary locally throughout different regions of the attractor \cite{casdagli1991state}. This behaviour was also highlighted by Uzal in his extension of Casdagli's noise amplification and distortion methods \cite{uzal2011optimal}. Additionally, we should also consider that invariant measures such as the Lyapunov exponent also vary locally \cite{abarbanel1991variation, algar2015noise}. Hence, the selection of a single embedding lag implies that all these variations may be averaged. 

Non-uniform delay embedding has been proposed as an natural extension of uniform embedding that aims to address some of the latter's limitations. Non-uniform delay embedding requires the selection of multiple delay lags $\{\tau_1, \tau_2,...\tau_{m-1}\}$ in order to construct a delay vector,

\begin{equation}
    \vec{x}(t) = (x(t), x(t-\tau_1), ..., x(t-\tau_{m-1})).
\end{equation}

The selection of delay lags represents a combinatorially hard problem that grows with increasing embedding dimension. The methods proposed for constructing non-uniform delay embedding often involve the iterative selection of time lags to gradually construct a delay vector until the required embedding dimension is reached.  In this section, we give an overview of the various methods that have been proposed to solve and automate this problem. These methods include the continuity statistic \cite{pecora2007unified}, PECUZAL \cite{kramer2021unified}, maximising derivatives on projections (MDOP) \cite{nichkawde2013optimal}, reduced autoregressive models \cite{judd1998embedding,hirata2006reconstructing}, \textcolor{edit_col}{and search optimisation algorithms such as ant colony optimisation \cite{shen2013optimal} and Monte Carlo decision tree search (MTCDS) \cite{kraemer2022optimal}.}

\subsection{Garcia and Almeida}
One of the earliest proposed methods of choosing non-uniform delays was proposed by Garcia and Alemeida \cite{garcia2005nearest}. They proposed a variation of the nearest neighbours methods of Kennel and Hegger applied to the problem of selecting time delays. Their method also recursively selected lags using a proposed $N$-statistic over multiple embedding cycles. At the end of each cycle, the false nearest neighbours algorithm is used to assess the quality of the newly constructed embedding. This process is repeated until the false nearest neighbour statistic $F$ decreases below a critical threshold. 

For the selection of the first time lag $\tau_1$, a 2D delay embedding $\vec{x}(t)$ is first done with respect to some prospective time lag $\tau^*$ to be tested,

\begin{equation}
    \vec{x}(t) = (x(t), x(t-\tau^*)).
\end{equation}

The closest neighbour $\vec{x}(t_j)$ for each point $\vec{x}(t_i)$ in the embedding reconstruction is identified. Neighbours should be chosen such that they are not temporally close (i.e. with respect to some Theiler window) \cite{theiler1986spurious, theiler1990estimating}. This is to ensure that their spatial proximity is not purely due to their temporal proximity. The two Euclidean distances $d_{1,\tau^*}(t_i), d_{2,\tau^*}(t_i)$ between any given two points are then calculated as follows,

\begin{align}
    d_{1,\tau^*}(t_i) &= || \vec{x}(t_i) - \vec{x}(t_j) || \\
    d_{2,\tau^*}(t_i) &= || \vec{x}(t_i+\delta t) - \vec{x}(t_j+\delta t) ||,
\end{align}
where $\delta t$ is the sampling time of the data. Simply put $d_{1,\tau^*}$ is the spatial separation between pairs of nearest neighbours in the reconstructed state space, and $d_{2,\tau^*}$ is the resulting separation one step forward in time. The $N$-statistic is taken as the proportion of points whose distances ratio $d_{2,\tau^*}/d_{2,\tau^*}>10$,

\begin{equation}
    N(\tau^*)= \frac{1}{N}\sum_{i=1}^N \mathbbm{1} \left( \frac{d_{2,\tau^*}(t_i)}{d_{1,\tau^*}(t_i)}>10 \right),
\end{equation}
where $N$ is the length of the time series and $\mathbbm{1}$ is the indicator function. The threshold of 10 was heuristically selected by the authors based on the numerical calculations of Kennel at el. \cite{kennel1992}. The time lag corresponding to the first minimum in $N(\tau^*)$ is taken to be the embedding lag. 

For non-uniform delay embedding, the selection of additional lags for each subsequent embedding cycle is done using a similar procedure. However, the reconstructed space used to calculate nearest neighbours and pairwise distances are calculated conditional on previously selected lags. Therefore, the selection of the $m^{th}$ embedding lag in a non-uniform embedding procedure will require neighbours and distances $d_{2,\tau^*},d_{2,\tau^*}$ to be calculated using the embedding with $m-1$ lags that have already been chosen and the new candidate lag $\tau^*$,
\begin{equation}
    \vec{x}(t) = (x(t),x(t-\tau_1),...,x(t-\tau_{m-1}),x(t-\tau^*)).
\end{equation}

\subsection{Continuity Statistic}
The continuity statistic was first proposed by Pecora et al. as a way to procedurally construct non-uniform delay vectors based on the idea of functional independence between vector coordinates \cite{pecora2007unified}. Takens' and Sauer both discussed the requirement that an embedding reconstruction requires vectors whose coordinates are independent \cite{rand2006dynamical, sauer1991embedology}. Pecora et al. proposed using a test for calculating the functional dependence between the components of a delay vector's components in order to assess the quality of an embedding. A functional dependence between vector coordinates implies,

\begin{equation}
    x(t-\tau_{m}) = F(x(t), x(t-\tau_1),...,x(t-\tau_{m-1})),
    \label{eq:functional_dependence}
\end{equation}
where $F$ is some arbitrary function. Constructing a non-uniform delay embedding requires iteratively building of a collection of time lags $\vec{\tau} = \{ \tau_1, ...\tau_{m-1}\}$ that minimises the likelihood of a functional dependence between components. In each iteration, a prospective lag $\tau_i$ is tested for functional dependence with the existing lagged components corresponding to $\vec{\tau}$. If there is no significant functional dependence, then $\tau_i$ may be added to the collection of lags. To test the equality of Equation \ref{eq:functional_dependence}, the authors assume that $F$ is smooth and use the property of continuity to quantify functional dependence. 

\begin{figure}
    \centering
    \includegraphics[width = 0.4\textwidth]{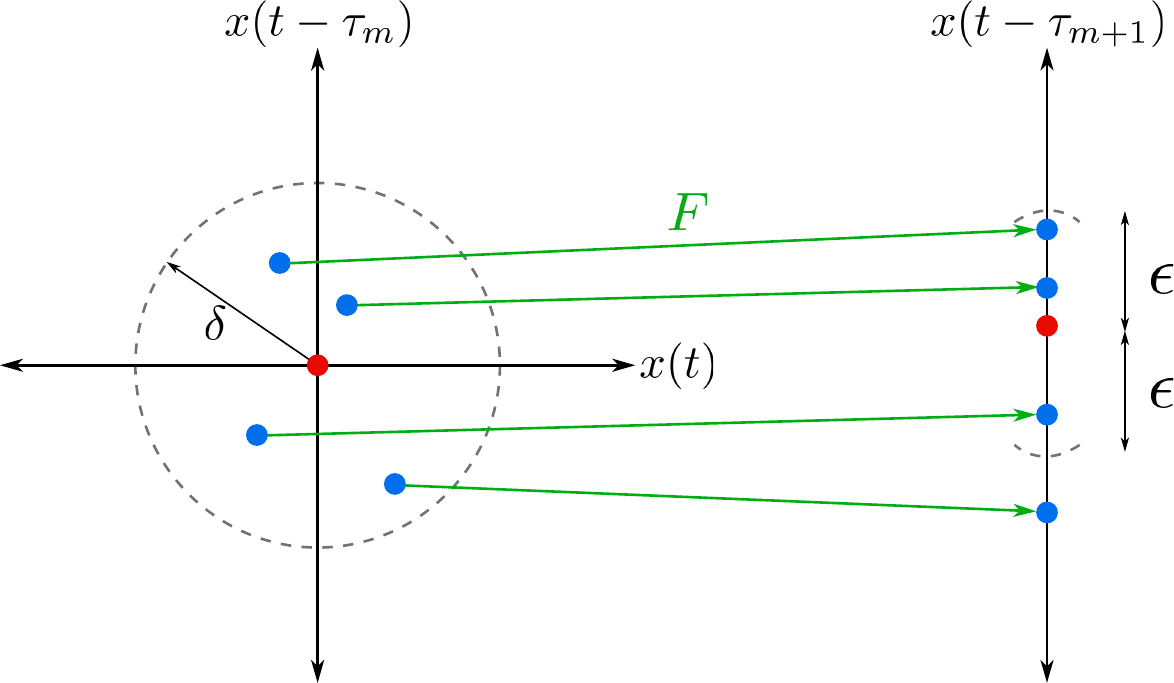}
    \caption{Illustration of the proposed functional dependence $F$ that is tested with the continuity statistic. The reference point (red) has 4 neighbours (blue) within a ball of radius $\delta$. Under the mapping $F$, only 3 of these neighbours lie within a distance $\epsilon$ from the image of the reference point. The value $p$ is the fraction of neighbours in embedded space (left) that are also neighbours in the potential new lagged component (right). The value of $\epsilon$ is adjusted until $p$ is insufficient to reject the null hypothesis that the probability of neighbours in embedded space (left) are also neighbours in the new lagged component (right) is binomially distributed with $p^*=0.5$.}
    \label{fig:continuity}
\end{figure}

Consider an existing $m$-dimensional embedding $\vec{x}_m(t) \in \mathbb{R}^m$ constructed from lag $\tau=\{ \tau_1, ... \tau_{m}\}$ and a potential new embedding lag to be tested $\tau_{m+1}$. To test the functional dependence of a new lag, select a reference reference point $\vec{x}_m(t_0)$ in embedded space. If a smooth functional dependence exists, then the continuity condition states that points $\vec{x}_{m,i}$ nearby the reference point ($|| \vec{x}_{m,i} - \vec{x}_m(t_0) ||<\delta$) in reconstructed space $\mathbb{R}^m$ should have lagged $m+1^{th}$ components that are also close by to each other ($|x_i(t-\tau_{m+1}) - x(t_0-\tau_{m+1})|<\epsilon$) (see Figure \ref{fig:continuity}).

The proportion $p$ of points $\vec{x}_{m,i}$ whose lagged components lie within $\epsilon$ of the reference point's lag component can be calculated. This proportion is then compared against a null hypothesis; that correspondence between these sets is purely by chance. Large values of $p$ suggest a strong relationship between the $m$-dimensional reconstruction and the new $\tau_{m+1}$ lagged component. Pecora et al. suggest the usage of a binomial distribution with a critical value of  $p^*=0.5$ in order to decide if a functional dependence exists with respect to some chosen $\epsilon$ due to its simplicity and robustness to noise \cite{pecora2007unified}. 

For a given $\tau_{m+1}$ to be tested and a sample of points, the continuity test is applied with decreasing values of $\epsilon$ until the null hypothesis fails to be rejected. The smallest possible value for rejecting the null hypothesis is given as $\epsilon^*$. This value is averaged over a collection of reference points sampled from the data to calculate the continuity statistic $\langle \epsilon^*\rangle(\tau_{m+1})$.

During each iteration of choosing a candidate lag $\tau_{m+1}$ for an existing collection of lags $\vec{\tau}=\{ \tau_1,..., \tau_{m}\}$, the continuity statistic profile $\langle \epsilon^*\rangle(\tau_{m+1})$ is calculated. The new lag $\tau_{m+1}$ is taken as the lag corresponding to the relative maxima of the continuity statistic profile. This is repeated until the desired embedding dimension (as per Whitney's theorem) is reached. Pecora et al. also propose an undersampling statistic that can be used as a termination criterion for iterative selection of time delays. Further details can be found in the original paper. 

This method was applied to a 2-torus, yielding embedding lags and embedding dimensions that were matching with theoretical expectations \cite{pecora2007unified}. The resulting reconstructed attractor was also found to be visually optimal. Similar results were gained when applied to the Lorenz chaotic time series. However, the resulting reconstruction appeared to be visually overfolded. 

\subsection{PECUZAL}
A criticism of the continuity statistic method is the ambiguity in selecting the optimal lag $\tau$ at each embedding iteration \cite{kramer2021unified}. In the original paper of Pecora et al., the definition of `relative maxima' is unclear and there is no objective criterion for selecting the best lag between multiple prospective local maxima \cite{kramer2021unified}. Additionally, the method also does not consider the effects of selecting different distances $\delta$ used to define nearby points in the reconstruction. Finally, the undersampling statistic originally proposed as a breaking condition for the embedding algorithm is computationally intensive, and does not inform on which of the prospective lags should be selected. A more detailed critique is provided by Kraemer et al \cite{kramer2021unified}.

Kraemer et al. suggested that the continuity statistics approach could be combined with Uzal's L-statistic \cite{uzal2011optimal} in order to provide a fully automated method of constructing non-uniform embedding delays. In their paper, they provide a workflow that uses the continuity statistic to perform a coarse search of multiple lag times and identify a small set of potential lags. These usually correspond to the various local maxima of the continuity statistic profile $\langle \epsilon^*\rangle(\tau_{m+1})$.

The L-statistic is then used as an assessment criteria to select which of the prospective lags should be selected in each embedding cycle. This addresses the problem of ambiguity of selecting lags that is present when using continuity statistics. The prospective lag whose new extended delay embedding resulted in the largest decrease of the L-statistic is selected in each embedding cycle. The L-statistic also provides a breaking condition for the embedding algorithm. The embedding cycles end when there is no achievable decrease in the L-statistic from the collection of prospective lags, i.e. $\Delta L >0$ between successive embedding cycles.

\subsection{Maximum Derivatives on Projection (MDOP)}
The maximum derivatives on projection (MDOP) method was first proposed by Nichkawde as a geometrical alternative to the statistics and information theoretic approaches of mutual information and continuity statistics \cite{nichkawde2013optimal}. MDOP optimises an embedding based on the criteria that a reconstructed attractor should be maximally unfolded and minimise redundancy in the delay components. Similar to the majority of non-uniform embedding methods, MDOP recursively constructs the delay vector through embedding cycles. Each cycle identifies a new time lag that maximises the directional derivative $\phi'_d(\tau)$ of points in reconstructed state space.

Like Pecora's approach in continuity statistics, MDOP begins with the criterion of functional dependence between each new prospective lag and an existing time delay reconstruction (see Equation \ref{eq:independence}). However, unlike in continuity statistic, Nichkawde suggests using the directional derivative of the functional dependence $F$ (see Equation \ref{eq:functional_dependence}) to quantify the degree of redundancy in the embedding and unfolding of the reconstructed attractor. This directional derivative is given by,

\begin{equation}
    F'_{m}= \left| \frac{\Delta F_m}{\Delta x_m} \right|
    \label{eq:independence}
\end{equation}
where $\Delta x_m$ corresponds to the spatial distance between a pair of nearby neighbour points $\vec{x}(t_i), \vec{x}(t_j)$ in  reconstructed space with $m-1$ dimensions,

\begin{align}
    \Delta x_{m,ij}&= ||\vec{x}(t_i)-\vec{x}(t_j)||\\
    &=\sqrt{\sum_{k=1}^{m-1} [x(t_i-\tau_k)-x(t_j-\tau_k)]^2}\\
    \tau_0=0.
\end{align}

The sampled pair of points should also be chosen such that spatial closeness is not due to them being virtually close in time \cite{theiler1986spurious, theiler1990estimating}. This is easily achieved by allowing for a Theiler window $l_{\rm Theiler}$ where $|i-j|>l_{\rm Theiler}$. 

Testing the inclusion of a new time lag $\tau_m$ requires evaluating the spatial variation in the prospective new lagged component, $\Delta F_{m,ij}$, and is given by Equation 
\ref{eq:spatial_variation},

\begin{equation}
    \Delta F_{m,ij}= |x(t_i-\tau_m)-x(t_j-\tau_m)|.
    \label{eq:spatial_variation}
\end{equation}

This quantity is used to evaluate the directional derivative of a small region on the reconstructed attractor,

\begin{equation}
    F'_{m,ij}= \left| \frac{\Delta F_{m,ij}}{\Delta x_{m,ij}} \right|.
\end{equation}

The directional derivative is evaluated with respect to each prospective new time lag $\tau_m$ and is averaged across randomly sample close pairs of points across the entire reconstructed attractor,

\begin{equation}
    \beta_d(\tau_m)=\langle F'_{m_{i,j}} \rangle _{i,j},
\end{equation}
where $\langle ...\rangle_{i,j}$ corresponds to the geometric mean across all sampled pairs of points. The author proposes using a geometric mean due to its robustness in the presence of outliers . In each recursive embedding cycle, the lag $\tau_m$ that maximises the directional derivative $\beta_d(\tau_m)$ is selected to be used for the reconstruction in the next cycle. This process is repeated until the desired number of embedding dimensions $m$ is reached, where $m$ is chosen via a number of different embedding dimension estimation methods such as false nearest neighbours etc.

\subsection{Reduced Autoregressive Models}
The reduced autoregressive model for non-uniform embedding was proposed by Judd and Mees \cite{judd1998embedding,small1999detecting} as a proposed method of constructing ideal models with respect to some information criterion. This method involves the construction of a pseudo-linear autoregressive predictive model with all $k$ possible lagged components as inputs or basis functions,
\begin{equation}
    x(t+1) = a_0 + \sum_{i=1}^k a_i x(t-\tau_i) + \epsilon (t),
\end{equation}
where $\lambda = (a_0,...,a_k)^T$ are the coefficients of each input to be determined and $\epsilon _t$ are the model prediction errors. For a time series of length $N$, construct a matrix $V_\tau$ with each row containing a vector of lags at a given time $t$,
\begin{equation}
    V_\tau = \begin{pmatrix}
    x(\tau_k+1) & x(\tau_k) & \dots & x(1)\\
    x(\tau_k+2) & x(\tau_k+1) & \dots & x(2)\\
    \vdots & \vdots & \ddots & \vdots \\
    x(N) & x(N-1) & \dots & x(N-\tau_k)
    \end{pmatrix}.
\end{equation}
The matrix $V_\tau$ has dimensions $(N-k)\times k$ and is defined with respect to a set of lags $B={\tau_1,...,\tau_k}$ and $\tau_1 < \tau_2 <...< \tau_k$.

Estimates for the coefficients of $\hat{\lambda}$ can be calculated using least squares regression,
\begin{equation}
    \hat{\lambda} = (V_\tau^T V_\tau)^{-1} V_\tau^T \xi,
\end{equation}
where $\xi=(x(\tau_k+1),x(\tau_k+2),...,x(N))$. Therefore, the resulting model errors utilising the set of all possible lagged components $B$ can be calculated as $e_\tau=V_\tau \hat{\lambda}_\tau$.

In order to reduce the number lagged components to a smaller selection, Judd and Mees propose the method of minimum description length. The principle of minimum description length is an application of Occam's razor to the context of model selection. It defines the best model for a given time series prediction task is one that achieves the most concise description of the data. For model selection, this would require achieving a compromise between model accuracy and model complexity (i.e. model description length). Model description length $L$ may be approximated by,
\begin{align*}
    L &= L_\epsilon + L_a\\
    &= N\log(\bar{\epsilon}) + L_\lambda,
\end{align*}
where $L_\epsilon$ is the description length of the model errors which is a function of the length of the time series $N$ and the mean square prediction error $\bar{\epsilon}$, and $L_\lambda$ is the description length of the model parameters.

The algorithm for reducing the number of lagged components is as follows:
\begin{enumerate}
    \item Construct an empty set $B'$ of chosen lags, and $B$ of candidate lags
    \item Define the prediction error with respect to the chosen set of lags $B'$ as $e_{B'} = \xi - V_{B'}^T \lambda_{B'}$ where $ V_{B'}^T$ and $\lambda_{B'}$ are defined with respect to the smaller set of lags $B'$. If $B'$ is empty, then $e_{B'} = \xi$.
    \item Calculate the vector $\mu=V_\tau^T e_{B'}$ and identify the index $p$ of the largest magnitude element corresponding to the most significant lag component. Add this lag $\{B'\cup {\tau_p}\} \to B'$
    \item Recalculate $\mu$ with respect to the new set of chosen lags $B'$ and verify that least significant lag component was the lag $\tau_p$ that was most recently added. Otherwise, return to step 2.
    \item Evaluate the model description lengths $L_B$ and $L_{B'}$ sets of lags $B$ and $B'$. If $L_{B'}<L_B$, return to step 2. Otherwise, end the algorithm and return the set of chosen lags $B' \to B$.
\end{enumerate}

An implementation of the minimum description length criterion for optimal embedding lag and window was done by Small and Tse \cite{small2003optimal,small2004optimal}.An extension of this method was proposed by Hirata et al. \cite{hirata2006reconstructing} where a normalised maximum likelihood $\mathcal{L}$ is used in place of minimum description length $L$ for model selection \cite{rissanen2000mdl,nakamura2006comparative}. Hirata et al. also propose a variation of the above algorithm by Judd and Mees - cross-validation - that utilises the radial basis modesl instead of pseudo-linear models.

\subsection{Search Optimisation Algorithms}
\textcolor{edit_col}{Many of the common non-uniform embedding methods involve a single optimisation step in each embedding cycle. In contrast, search optimisation algorithms attempt to search across state space of possible lags to identify ideal combination embedding lags without necessarily selecting the first local optima encountered. Two examples of such approaches are the ant colony optimisation (ACO) method \cite{shen2013optimal} and Monte Carlo decision tree search (MTCDS) \cite{kraemer2022optimal}.}

Ant colony optimisation (ACO) is a swarm intelligence method first proposed by Dorigo \cite{dorigo1996ant, dorigo1997ant} that is inspired by the foraging behaviour of ant colonies. Similar to other swarm optimisation methods, ACO initialises a number of agents (`ants') that simulteneously perform an initial search of the solution space. The quality of each attempted solution is assessed according to an objective function and a `pheremone' is assigned to the corresponding search path. These pheremones are able to accumulate and fade over time. This biases the search direction of subsequent iterations of the algorithm and is reminiscent of the optimal path finding behaviours of real world ant colonies.

The ant colony optimisation method applied to non-uniform embedding (ACO-NE) builds upon this framework in a few ways. (i) By using an objective function based on various notions of embedding quality (mean neighbourhood distance, minimum false nearest neighbours, minimum description length \cite{small2003optimal,small2004optimal}) to optimise parameters, (ii) Incorporating heuristics into the algorithm to speed up convergence. Interested readers are encouraged to refer to the original paper for more details \cite{shen2013optimal}.

\textcolor{edit_col}{Another search optimisation algorithm, MTCDS, proposed by Kraemer et al. \cite{kraemer2022optimal} reframes the non-uniform lag selection problem into a decision tree search. Each embedding cycle is represented by a collection leaves or nodes stemming from a root (the original time series) where each leaf is the selection of a particular candidate embedding lag. A Monte Carlo approach is used to randomly sample the tree and identify various local optima for a given objective function and backpropagation is then used to decide on the best selection of lags in each step.}

\section{Persistent Strands and Characteristic Times}
\label{Persistent_Strands}
The focus of this paper is on the problem of selecting time delays for delay embedding, and in particular for non-uniform delay embedding. Many of the proposed methods for optimising embedding delay both in the uniform and non-uniform case generally fall into the broad categories of dynamics (e.g. mutual information, continuity etc.) or topology (fill factor, MDOP). With the exception of the PECUZAL automated embedding framework, non-uniform delay embedding strategies only focus on one of these two broad aspects in their definition of a good embedding. 

Another weakness in non-uniform embedding is that they do not always provide a full picture on the relative significance of each delay. When operating under the iterative construction of delay vectors, each prospective new time delay must be reevaluated with respect to the most recently updated embedding. Hence, the significance of each subsequent delay is conditional on the previous selected sequence of delays. This is a weakness particularly in fully automated algorithms such as MDOP and PECUZAL where ideal embedding lags are selected automatically, with little to no reference on their relative impact on embedding quality. Additionally, there is often a lack of consistency between the results of different methods. This can be attributed to the fact that optimisation is done with respect to different notions of embedding quality. Often, these methods do not provide a dynamical explanation for each time lag's significance and ability to improve a given embedding.

In view of this, we argue that a good non-uniform embedding method should have two main qualities. Firstly, the embedding criteria should utilise both dynamical (irrelevance, periodicity, independence of coordinates) and topological (attractor unfolding) features in their selection of embedding lag. Secondly, the significance of each selected embedding lag should be explainable.

As our contribution, we propose a new method, `significant times on persistent strands' or SToPS, of identifying non-uniform embedding time lags using techniques drawn from persistent homology and recurrence analysis. We introduce the idea of a characteristic time scale spectrum of a signal based on the periodicities of time series and show how this may be used to identify ideal time delays. The selection of multiple time delays are also treated independently, marking a contrast to the iterative approach of constructing delay vectors that is common in most non-uniform delay embedding methods. We demonstrate the performance of our method on a collection standard periodic and chaotic time series from the literature. Additionally, we explore its performance on experiment neuron data containing fast-slow dynamics and show that SToPS is sensitive in identifying explainable time delays.

\subsection{Introduction to Persistent Homology}
Persistent homology has seen a recent growth in popularity particularly in the fields of dynamical systems \cite{adams2020fractal, Jaquette2020, Khasawneh2016, Myers2019, tan2021grading}. We also note that recent work has also been done attempting to automate the selection of delay embedding parameters with persistent homology \cite{myers2022delay}. In its essence, persistent homology aims to quantify and track the evolution of the topological properties of an object (network, point cloud data etc.) under an increasing notion of distance \cite{Jaquette2020}. The process of gradually increasing distances is referred to as ‘thickening’. Topological features that persist for a large interval of distance under the thickening process are observed to be significant. For simplicity, we refer to the coordinates or increasing distance $\epsilon$ as analogous to increasing time. Conversely, short-lived topological changes are typically perceived to be noise.

We describe the thickening process is as follows. Consider a point cloud $\mathcal{A}=\{ \boldsymbol{x}_i \}$ arranged in the pattern of a circle. Place an open ball $B_\epsilon(\boldsymbol{x}_i)$ of radius $\epsilon$ centred at each point and let the union of all open balls be the set of interest, 
\begin{align}
 \mathcal{B}(\epsilon) &= \bigcup\limits_{i=1}^{N} B_\epsilon(\boldsymbol{x}_i),
\end{align}
where we are interested in calculating the homology of the set $\mathcal{B}(\epsilon)$ whose complex is given by $K(\mathcal{B}(\epsilon))$. As $\epsilon$ is increased the set of open balls also increase in size, forming a filtration shown below, 
% \begin{align}
% K: \, \mathcal{B}(\epsilon_0) \subseteq \mathcal{B}(\epsilon_1)  \subseteq ... \subseteq \mathcal{B}(\epsilon_n), \, \epsilon_0 \leq \epsilon_1 \leq ... \leq \epsilon_n
% \label{eq:filtration}
% \end{align}
\begin{align}
K_0 \subseteq K_1 \subseteq ... \subseteq K_n, 
\label{eq:filtration}
\end{align}
where $K_n = K(\mathcal{B}(\epsilon_n))$ and $\epsilon_0 \leq \epsilon_1 \leq ... \leq \epsilon_n$. More precisely, each set $K_n$ represents a collection of simplicial complexes, with subsequent $\epsilon_n$ yielding a nested sequence called a filtration. Computation of the thickening process is well documented with two main algorithms employed, Vietoris-Rips \cite{vietoris1927} and Čech \cite{hatcher2005}, of which we will employ the former.

For this case, we are concerned with changes in the homology of the set as $\epsilon$ increases. Simply, persistent homology aims to enumerate and track the number of $n$-dimensional ‘holes’ in the set \cite{Maletic2016}. Namely, $H_0$, $H_1$ and $H_2$ for cases of low dimensional homology. The 0, 1 and 2 dimensional holes correspond to disjoint components, cycles and voids respectively. Simplices (triangles, tetrahedrons etc.) are considered solid components and are therefore not holes.

\textcolor{edit_col}{In the context of low-dimensional strange attractors, persistent homology tracks the persistence or lifetime (death time minus birth time) homological features (usually $H_1$ holes) of the attractor using spatial data from an embedding or otherwise. All the information pertaining to the birth and persistence of features can be represented easily in a persistence diagram. In a persistence diagram, each homological feature (i.e. $H_0$ - disjoint components, $H_1$ - holes) are represented by a plotted point with coordinates $(\epsilon_b, \epsilon_d)$ corresponding to its birth and death times respectively.}

\begin{figure*}
    \centering
    \includegraphics[width = 0.9\textwidth]{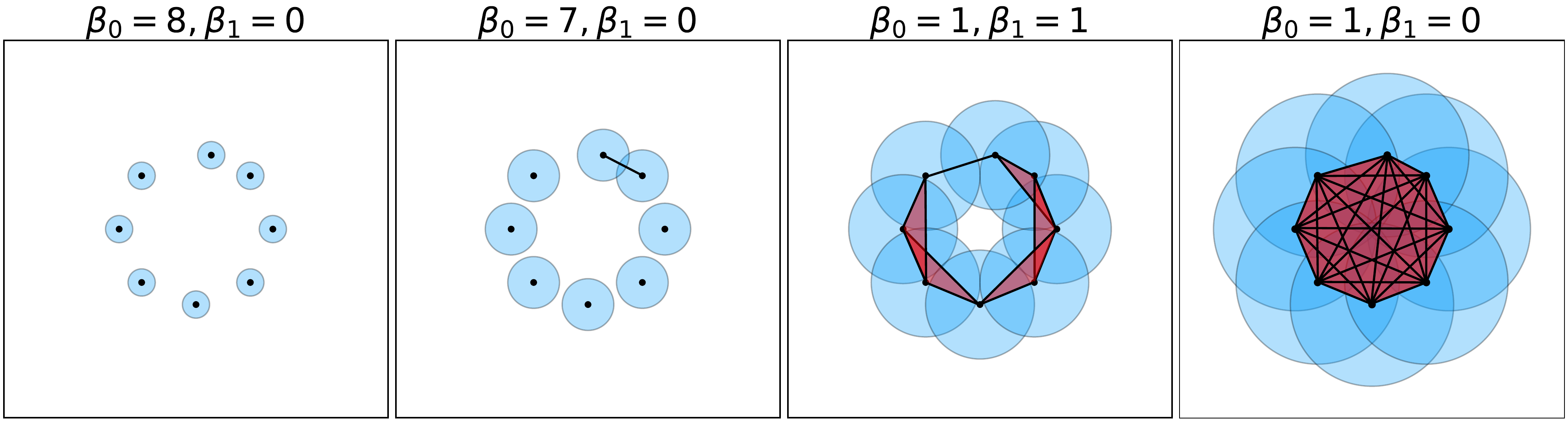}
    \caption{\textcolor{edit_col}{Illustration of the filtration process with Betti numbers $\beta_0, \beta_1$ illustrating the number of $H_0$ (disjoint components) and $H_1$ (holes) in the homology respectively. Filtration is with increasing open ball radius $\epsilon$. Initially for small $\epsilon$, the constructed balls form eight disjoint components with no holes ($\beta_0=8, \beta_1=0$). Further increases in $\epsilon$ eventually cause the intersection of two balls, resulting in a merging of two components. Subsequently, $\epsilon$ is increased until all components are merged into a single components with one unfilled center ($\beta_0=1, \beta_1=1$). Eventually, $\epsilon$ is increased until the hole in the center is filled in (shaded in red) resulting in a single solid component with no holes ($\beta_0=1, \beta_1=0$).}}
    \label{fig:homology}
\end{figure*}

The tracking of birth and death time of homological features presents two useful features. \textcolor{edit_col}{Firstly, it allows the tracking of the locations of holes within the data. When analysing phase space trajectories, these holes may correspond to short-term pseudo-periods or turning points in the time series (see Figure \ref{fig:circularity_efficiency}). By tracking $\epsilon$ and the representative cycles, one can also identify boundary points of holes. Secondly, the lifetimes of homological features allow an estimate of the relative size of the feature, which may or may not be related to its significance depending on the type of data.} All of these features may be represented in a persistent diagram (see Figure \ref{fig:homology}) where points further from the diagonal represent more persistent homological features.

\subsection{Characteristic Time Scales}
One approach to selecting lags for non-uniform embedding would be to select values related to the natural time scale of the system's dynamics. Picking lags that are much smaller or larger would logically correspond to the cases of high redundancy and irrelevance respectively. Consider the simplest case of a periodic signal. One can argue that an ideal embedding would require a delay that is related to the time scale of its main dynamics i.e. it's periodicity. For this we employ the quarter period heuristic in the definition of the characteristic time $\tau$,
$$
\tau = \frac{T}{4}.
$$

A natural progression of this concept into more complex dynamics would be to take the lag from the collection of natural frequencies,
$$
\tau_i=\frac{2\pi}{4\omega_i}.
$$
Because $\tau_i$ is evaluated at individual frequencies $\omega_i$, which may vary greatly in magnitude, it is possible to capture the dynamics of systems with multiple time scales, such as fast-slow dynamical systems. Relating embedding lag with natural frequencies and periodicity also introduces a degree of explainability to the selection of $\tau$ that also directly relates to the dynamics. Whilst this presents a potential advantage over typical non-uniform embedding methods, it requires one to be able to accurately measure the natural frequencies (or equivalent) from any given signal.

If the analysed signals were relatively smooth and easily decomposable into sinusoids, Fourier transforms would provide an excellent solution to this problem. However, this approach quickly fails when analysing discontinuous-like signals such as neuron voltages where the time series is characterised by alternating phases of bursting and resting dynamics. Fourier transforms also require the time signal to be stationary on its statistics. This property is not possible for chaotic time series where the phase and period of the signal varies over time. For example, the frequency power spectrum of a chaotic signal such as Lorenz and Chua produces a shape with an exponential decay. Similar arguments may be made when analysing experimental time series, where it is expected that drifting and oscillations in the phase and or period may occur.

Another alternative method to detecting periodicity would be to use the notion of recurrence distance employed in recurrence analysis and unstable periodic orbit detection \cite{marwan2007recurrence},

\begin{equation}
    d(\tau) = \frac{1}{N}\sum_{i=1}^N
    \left\| 
    \begin{pmatrix}
    x(t_i)\\
    x(t_i-\tau)
    \end{pmatrix}
    -
    \begin{pmatrix}
    x(t_i+4\tau)\\
    x(t_i+3\tau)
    \end{pmatrix}
    \right\|
    ,
% \sqrt{
% (x(t_i)-x(t_i+4\tau))^2+(x(t_i-\tau)-x(t_i+3\tau))^2
% },
\end{equation}
where $\tau$ is a characteristic time if the recurrence distance $d(\tau)$ is below some threshold $d^*$,
\begin{equation}
    d(\tau)<d^*.
\end{equation}

The recurrence distance tracks the displacement from a point in state space and its future trajectory. For a periodic orbit and correctly selected $\tau$, this will result in a local minima for $d(\tau)$ where the periodic orbit returns close to its initial position. This method of detecting characteristic times $\tau$ is unsuitable as it cannot distinguish between cases where $T=\tau/2$ and the resulting embedded trajectory clusters along the diagonal).

We propose a method of identifying and weighting the significance of characteristic times from a time series. The identification of characteristic times is done by \textcolor{edit_col}{sampling `strands' (short contiguous windows of 2D delay embedded time series) and calculating their persistent homology.} The representative cycles of the persistent strands' calculated homologies are used to assign a significance score to each identified characteristic time in order to construct a characteristic time spectrum. This spectrum can then be used to inform the selection of lags for non-uniform time delay embedding. This combined framework is named Significance Times on Persistent Strands (SToPS).

\section{Significance Times on Persistent Strands (ST\lowercase{o}PS)}
\subsection{Persistent Strands}
\label{SToPS}
The first challenge to tackle is the identification of all pseudo-periodic behaviour of period $T$. For a characteristic time $\tau$, the quarter period heuristic suggests that the corresponding 2D embedding with coordinates $(x(t), x(t-\tau))$ will result in a periodic orbit that is approximately maximally convex in reconstructed space. We note that the precise shape of this orbit is not guaranteed to be circular and varies depending on the signal.

To test for periodic behaviour with characteristic time $\tau$ for a signal $x(t)$, a 2D delay embedding 2D delay embedding of the time series with a single lagged component is constructed,

\begin{equation}
    \vec{x}(t) = (x(t), x(t-\tau)).
\end{equation}

A collection of $N$ strands of length $l=4\tau$ with random initial positions $t_i$ are then sampled from the time series, where each strand is given by,

\begin{equation}
    \vec{x}_i(t) = (x(t),x(t-\tau)), \quad t\in[t_i, t_i+4\tau].
\end{equation}

We argue that strands of this length should be approximately sufficient to detect loop structures in the 2D embedding based on the quarter period heuristic. The persistent homology of each strand can be calculated to detect the presence of orbits with period $l$. For each strand $\vec{x}_i(t)$, extract the maximum persistence from the resulting persistence diagram $PD$. A sample strand is said to contain an orbit of length $l$ if the maximum persistence of the corresponding diagram exceeds a critical value $\rho$. The value of $\rho$ is taken to be the average distance between two consecutive observations in phase space. Therefore, a na\"ive value $P_i$ that quantifies the significance of a characteristic time $\tau$ can be defined as,
\begin{equation}
    P_i = 
\begin{cases}
\text{maxpers}(PD_i), & \text{maxpers}(PD_i)>\rho\\
0 & \text{maxpers}(PD_i)<\rho\\
\end{cases}.
\end{equation}

This value can be used to calculate an overall maximum persistence spectrum $P(\tau)$ by averaging over all non-zero scores $P_i$ for each characteristic time $\tau$,

\begin{equation}
    P(\tau)=\langle P_i \rangle_{P_i \neq 0}.
\end{equation}

\begin{figure}

    \subfloat[Lags are too small and resulting strand length is too short to fully enclose resulting in the detection of spurious holes with very small number of boundary points.]{\includegraphics[width = 0.3\textwidth]{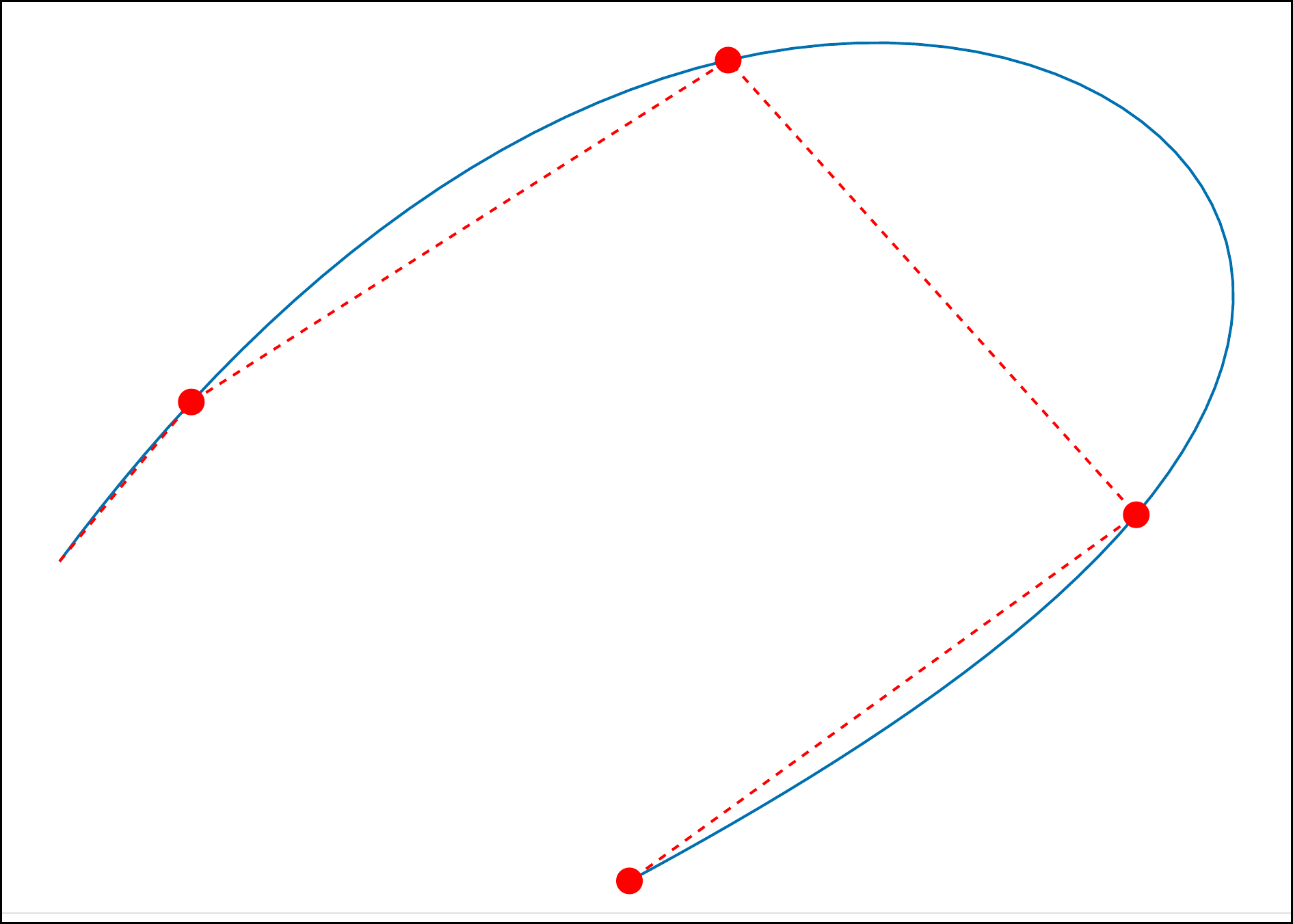}}
    
    \subfloat[Presence of a hole even at large lags due to overlapping trajectories on a similar orbit. Long strand does not efficiently use all the available points.]{\includegraphics[width = 0.3\textwidth]{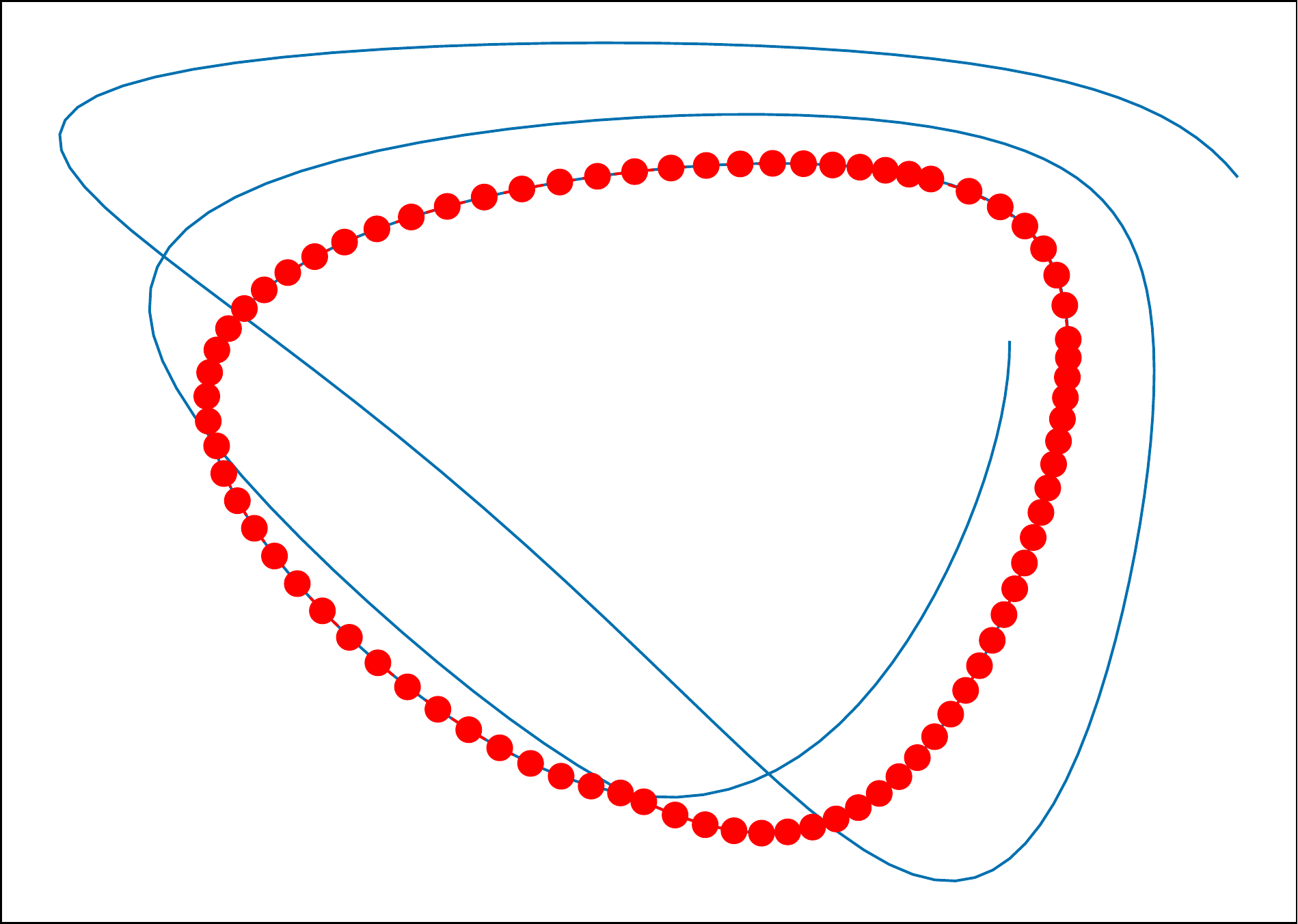}}
    
    \caption{\textcolor{edit_col}{Pathological problems of the SToPS method.}}
    \label{fig:spurious_lags}
\end{figure}

\textcolor{edit_col}{
We also impose an additional constraint that the number of points used to reconstruct the boundary of the hole ($H_1$ homology) in the orbit should be at least $N_{hole}$. This avoids the problem of including spurious holes where a small number of points suggest the existence of a hole even though the embedded strand is insufficiently long to close the orbit (Figure \ref{fig:spurious_lags}). We select a minimum value of $N_{hole} = 8$ for our analyses based on the argument that a minimum of 8 points should be sufficient in at least identifying a hole of a small lag $\tau=2$ without discounting higher lags. This value was found to work well in our analyses.}

% \begin{figure}
%     \centering
%     \includegraphics[width = 0.48\textwidth]{Pathological_Problems.pdf}
%     \caption{\textcolor{edit_col}{Pathological problems of the SToPS method: (1) Lags are too small and resulting strand length is too short to fully enclose resulting in the detection of spurious holes with very small number of boundary points. \textit{(left)} (2) Presence of a hole even at large lags due to overlapping trajectories on a similar orbit. Long strand does not efficiently use all the available points. \textit{(right)}}}
%     \label{fig:spurious_lags}
% \end{figure}

\begin{figure}
    \centering
    \includegraphics[width = 0.4\textwidth]{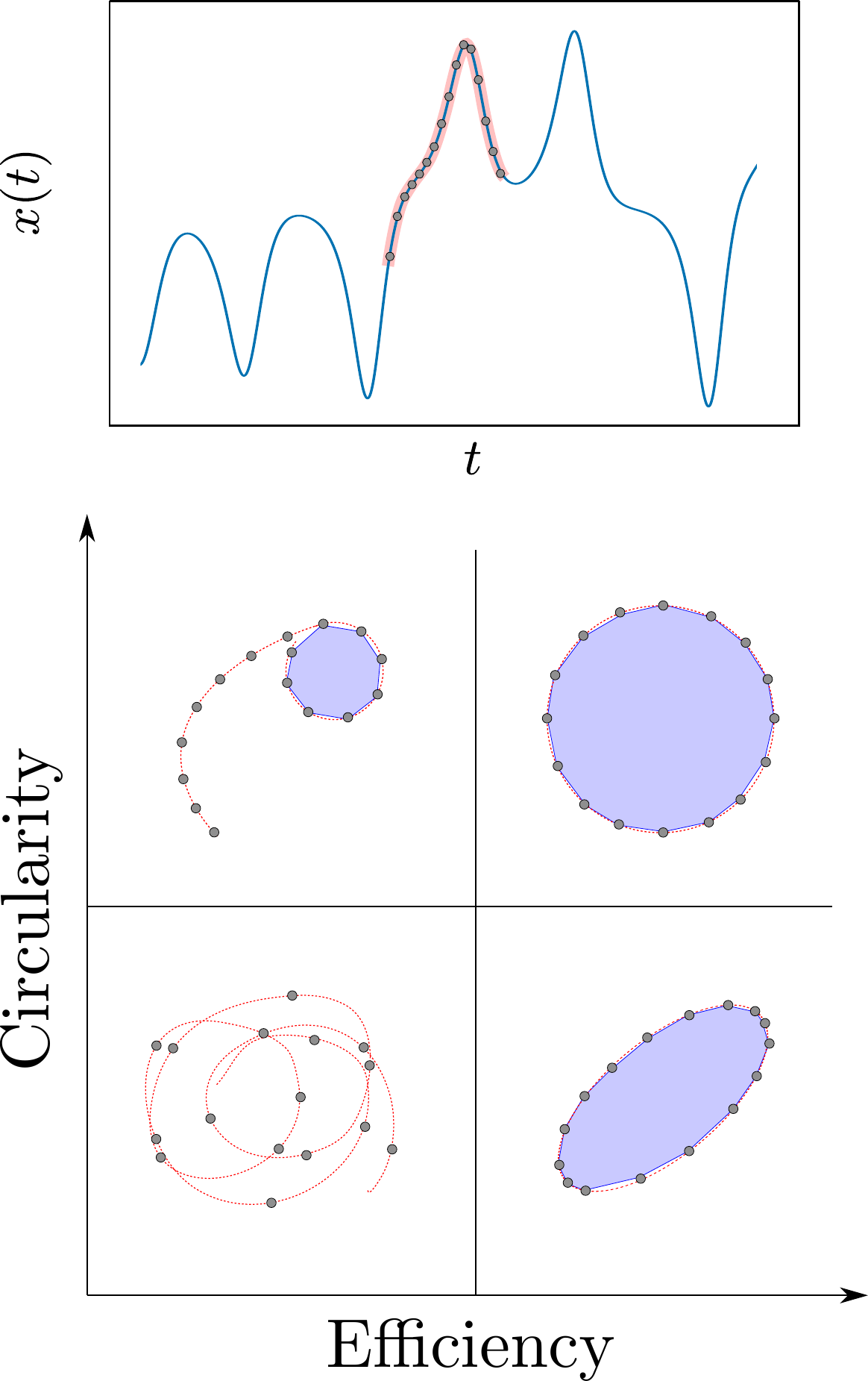}
    \caption{\textcolor{edit_col}{Illustration of a randomly sampled `strand' (shaded red) taken from a time series and the subsequent pathological cases of embeddings from different time scales $\tau$. Each case corresponds to non-optimal embeddings based on the notions of circularity and efficiency. Significant characteristic times should aim to maximise both metrics (top right quadrant). Blue shaded areas correspond to the polygon formed by the homology generators of the hole calculated from the embedded strand's representative homology.}}
    \label{fig:circularity_efficiency}
\end{figure}

The maximum persistence spectrum provides a simple way to quantify the degree to which $\tau$ corresponds to a dynamically significant timescale (i.e. periods). However, $P(\tau)$ is biased towards larger features due to its reliance on the lifetime of each homological feature. As a result, spatially small but dynamically significant features over small time scales are underrepresented. We also note that $P(\tau)$ is unable to differentiate between pathologically inefficient embeddings such as those whose loops are not maximally circular, or have trailing tails (see Figure \ref{fig:circularity_efficiency}). 

As we can see from Section \ref{embedding_quality}, a good embedding is one that has been maximally unfolded to best utilise the reconstructed state space, whilst being robust to the effects of noise amplification. Therefore, the significance of each characteristic time $\tau$ should be weighted according to how well the corresponding strand `unfolds' into a loop structure in the 2D delay embedding with lag $\tau$.

From this, we propose the significance score $S(\tau)$, which is a measure of the dynamical significance of each characteristic time $\tau$ that accounts for the quality of the unfolding of sampled strands. Using the topological notion of a good embedding, the significance score $S_i(\tau)$ for the $i^{th}$ sampled strand is defined as,
\begin{equation}
    S_i(\tau) = \alpha_i(\tau) \gamma_i(\tau),
\end{equation}
where $\alpha(\tau)$ and $\gamma(\tau)$ are two separate measures named the circularity and efficiency respectively. This score is also not biased towards larger homological features, allowing for the detection of both small and large pseudo-periodic dynamics.

The circularity $\alpha(\tau)$ tries to quantify the quality of the unfolding of a persistent strand in embedded space. Embeddings that yield circular loops imply a selected lag $\tau$ that maximally unfolds the dynamics (i.e. $\tau$ may correspond to a characteristic time). Therefore, loops that are more circular or regular have higher values compared to ellipses and other shapes with eccentricities. This allows circularity to also function as a measure of redundancy with lower $\alpha$ corresponding to high redundancy. 

To calculate $\alpha_i$ for a given sampled strand, the boundary points $\vec{b}_i$ corresponding to the birth of the most persistent homological feature are identified by examining representative cycles. \textcolor{edit_col}{Principal component analysis (PCA) is used to identify the major and minor axes of the embedded hole. Because each strand is used to evaluate a single lag in a 2D delay embedding, the first and second principal eigenvalues approximately correspond the relative sizes of the major and minor axes of the embedded points' bounding ellipse. Hence, we define circularity as the ratio between the first and second principal eigenvalues, averaged across all strands}

\begin{equation}
    \alpha(\tau) = \langle \alpha_i \rangle_{P_i \neq 0} = \biggl{<} \frac{\lambda_{2,i}}{\lambda_{1,i}} \biggr{>}_{P_i \neq 0},
\end{equation}
\textcolor{edit_col}{where $\alpha(\tau)\in(0,1]$. As $\alpha \to 1$, embedded holes are more uniformly circular.}

Efficiency $\gamma(\tau)$ is defined by the ratio of two areas,
\begin{equation}
    \gamma(\tau) = \frac{A_h}{A_{pc}},
\end{equation}
where $A_h$ is the area of the hole given by the ordered set of boundary points $\vec{b}_i$ and $A_{pc}$ is the area of the smallest convex polygon that includes all points in the strand $\vec{x}_i(t)$. \textcolor{edit_col}{The area of the hole $A_h$ can be simply calculated using the shoelace algorithm \cite{braden1986surveyor}. The area of the latter $A_{pc}$ can be similarly calculated by using a Graham scan algorithm to first identify the smallest bounding convex polygon \cite{graham1972efficient}. }Similar to $\alpha(\tau)$, the efficiency score is also bounded with $\gamma(\tau) \in [0,1]$. Efficiency is a measure of how well utilised the ambient space of an embedded strand is with respect to the hole. It is used to detect cases where the detected hole is circular, but does not utilise the full length of the strand (see Figure \ref{fig:circularity_efficiency}). 
% Efficiency also has the added benefit of filtering falsely toplogical noise from the persistent homology of the strand where the most persistent homological feature is an artificact of spatially close points.

\section{Testing Methodology}
\label{methodology}
SToPS was tested on three different types of time series covering periodic, chaotic and fast-slow dynamics. 

\begin{figure}

    \subfloat[Artificial time series constructed from sum of 3 sine terms with $\omega=\{1, 5, 30\}$, $\phi \{0,0.25,0.75\}$ and $A_i=\{1, 0.5, 0.2\}$.]{\includegraphics[width = 0.47\textwidth]{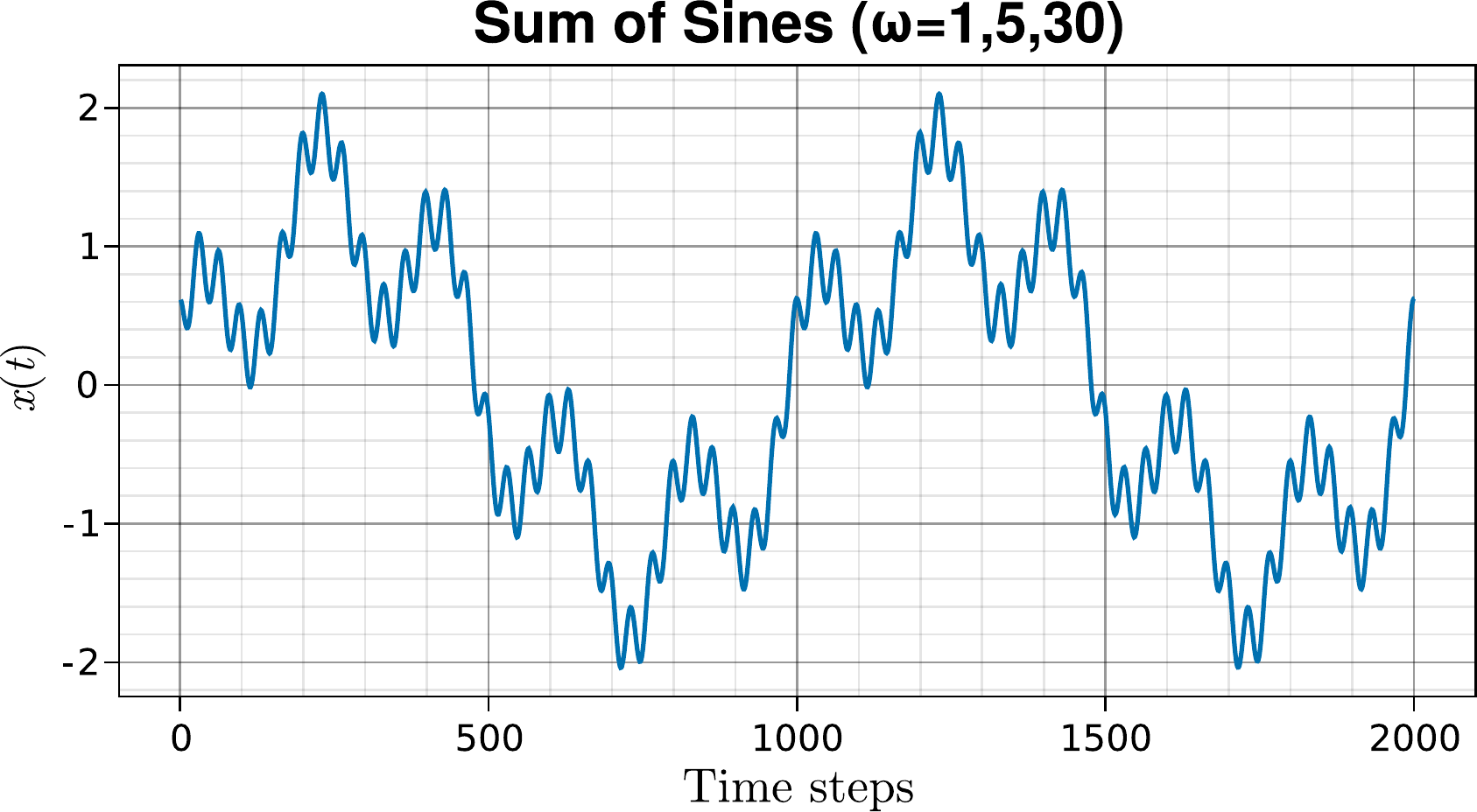}}
    
    \subfloat[First component ($x(t)$) of the chaotic Lorenz time series.]{\includegraphics[width = 0.47\textwidth]{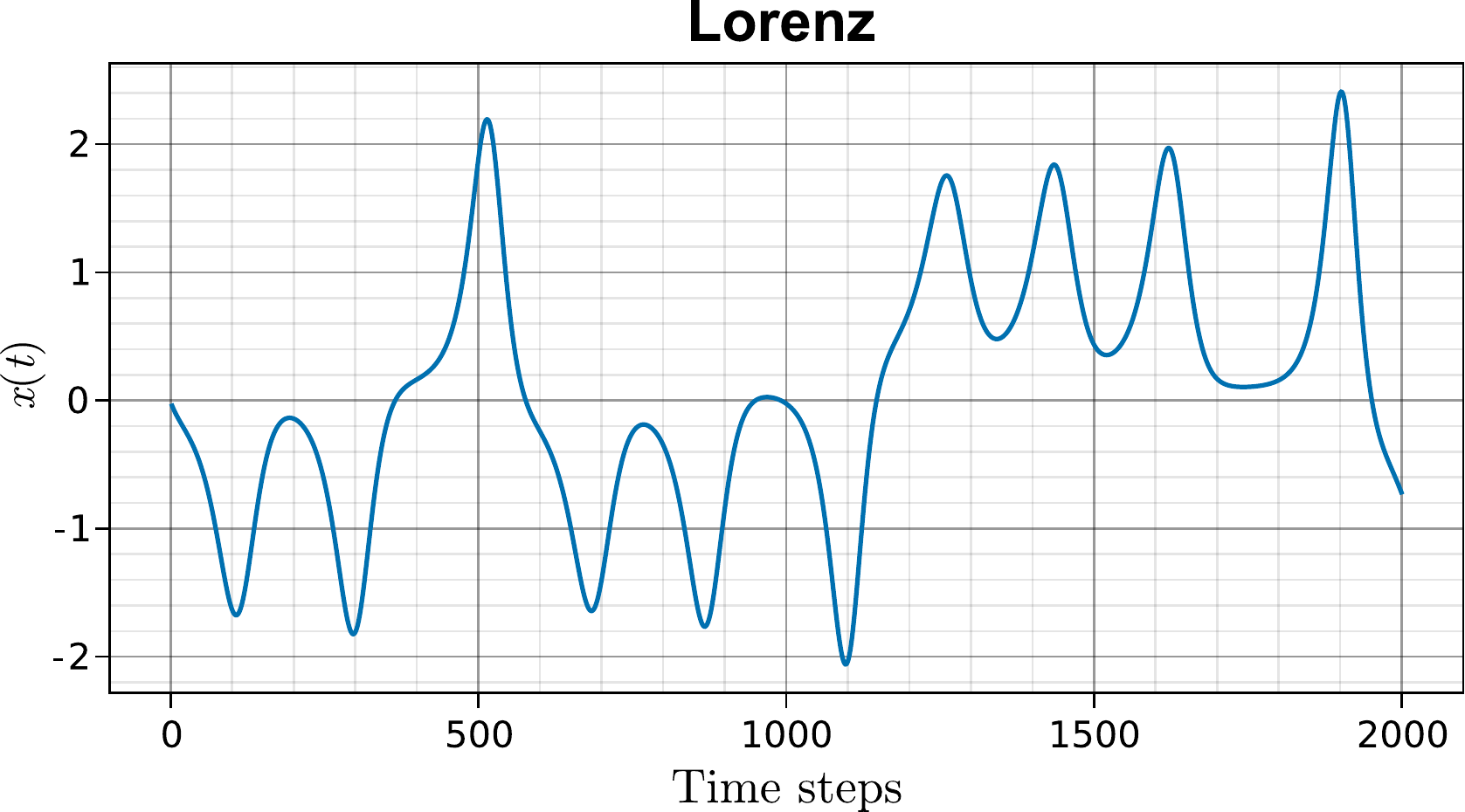}}
    
    \subfloat[Lobster LP neuron time series showing pseudo-periodic dynamics over two scales. Fast (spiking neurons) and slow (bursting neurons).]{\includegraphics[width = 0.47\textwidth]{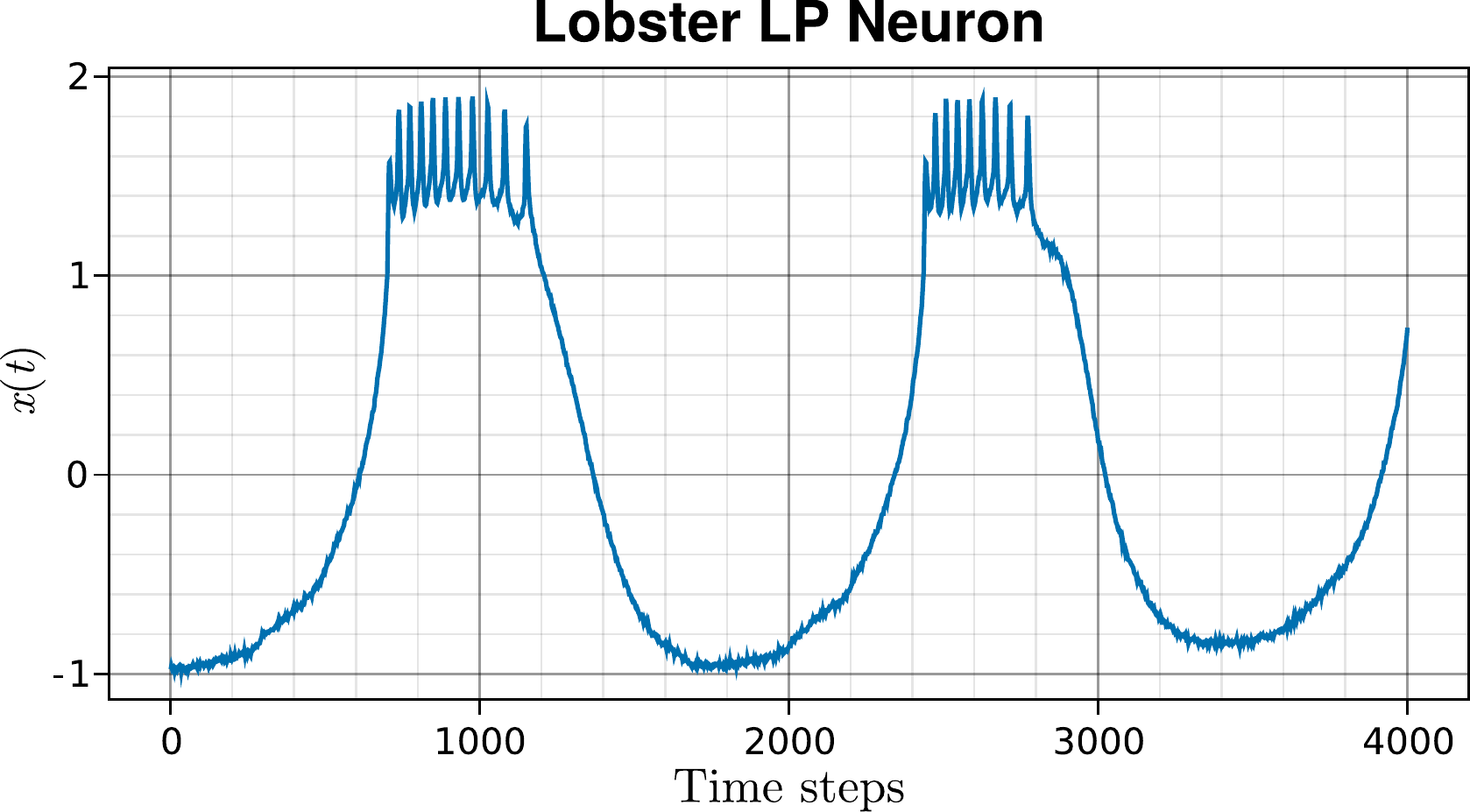}}
    
    \caption{Three classes of time series tested: sum of sines (periodic), Lorenz (chaotic), Lobster LP Neuron (multiple disparate time scales).}
    \label{fig:timeseries}
\end{figure}

Firstly, a sum of sines signal (see Figure \ref{fig:timeseries}) \textcolor{edit_col}{with step size $dt = 0.001$} was used to simulate the case of periodic time series,
\textcolor{edit_col}{
\begin{equation}
    x(t) = \sum_{i=1}^{N} A_i \sin(2\pi\omega_it+\pi\phi_i),
\end{equation}}
where $N=3$ is the number of superimposed of sine waves. Phases $\phi=\{0,0.25,0.75\}$ and magnitudes $A=\{1, 0.5, 0.2\}$ were selected to ensure that the spacial scales of the dynamics were distinct (see Figure \ref{fig:timeseries}).

% \begin{figure}
%     \centering
%     \includegraphics[width = 0.5\textwidth]{SINES_timeseries.pdf}
%     \caption{Artificial time series constructed from sum of 3 sine terms with $\omega=(1, 1/5, 1/30)$, $\phi=(0,0.25,0.75)$ and $A_i=(1, 0.5, 0.2)$}
%     \label{fig:SINES_timeseries}
% \end{figure}

The second time series analysed was the $x$-component of the canonical Lorenz chaotic time series. This is used to represent the case of chaotic time series. This time series was numerically integrated with a $4^{th}$ order Runge-Kutta scheme with an \textcolor{edit_col}{integration time step of 0.0004 and subsampled to an effective time step of $dt=0.004$ for 25000 steps.}

% \begin{figure}
%     \centering
%     \includegraphics[width = 0.5\textwidth]{LORENZ_timeseries.pdf}
%     \caption{First component ($x(t)$) of the chaotic Lorenz time series}
%     \label{fig:LORENZ_timeseries}
% \end{figure}

The third time series consisted of experimental data measured from a lobster somatogastric ganglion (STG) lateral pyloric (LP) neuron. This time series represents a third class of dynamics corresponding fast-slow dynamics with two different spatial and temporal scales. This time series was originally analysed by Abarbanel \cite{abarbanel1996synchronized} and includes two characteristic dynamics. These are small magnitude and time scale oscillations corresponding to neuron spiking dynamics, and a long time scale periodic behaviour for neuron bursting (see Figure \ref{fig:timeseries}). Additionally, the phase of the bursting dynamics also varies slightly over time. This results in a gradual shift of spatial position of the expected loop in a 2D embedding.  As a result, strands with lengths $T$ that are too long are penalised as these loops eventually get filled in.

% \begin{figure}
%     \centering
%     \includegraphics[width = 0.5\textwidth]{LOBSTER_timeseries.pdf}
%     \caption{Lobster LP neuron time series showing pseudo-periodic dynamics over two scales. Fast (spiking neurons) and slow (bursting neurons).}
%     \label{fig:LOBSTER_timeseries}
% \end{figure} 

\textcolor{edit_col}{All three input time series were normalised with zero mean and unit standard deviation before applying the embedding algorithms. A small amount of additive noise $\xi \sim N(0,0.001^2)$ was applied to the sum of sines and the digitised experimental Lobster LP neuron data to ensure there were not exact overlaps in values. In all cases, the input data was limited to 25000 points when calculating embedding lags to ensure the computational time was sufficiently short. This is due to the use of a $k$-means++ random sampler \cite{arthur2006k}, whose computation increases with the number of points, to select strand locations that uniformly explore the reconstructed state space.}

Different ranges of characteristic times were tested to calculate the efficiency, circularity and significance score profile depending on the type of time series. Each $\tau$ consists of sampling 250 strands of length $T=4\tau$. Due to the poor computational scaling of the Vietoris-Rips filtration, long strands where $T$ is large are subsampled by a factor of $k$,

\begin{equation}
    k = \left\lfloor \frac{4\tau}{N_p} \right\rfloor,
\end{equation}
where $N_p=250$ is the approximate scale of the maximum allowable strand length. \textcolor{edit_col}{This value was chosen in order to accommodate the slow computational efficiency of calculating the Vietoris-Rips filtration and associated persistent homology. }The threshold $\rho^*$ used to define the minimum lifespan needed to classify a homological feature as significant was chosen to be the average distance between temporally adjacent points,

\begin{equation}
    \rho^*=\langle ||\vec{x}(t_i)-\vec{x}(t_{i+1})|| \rangle_{i.k},
\end{equation}
where $\langle ... \rangle_{i,k}$ is the average across all points in the 2D delay embedded strand $\vec{x}(t)$ subsampled with a factor $k$. The circularity, efficiency and significance score profiles were calculated using SToPS. Both mean $S_\mu(\tau)$ and standard deviation $S_\sigma(\tau)$ profiles are calculated and compared against a selection of embedding delay optimisation measures. \textcolor{edit_col}{In our analyses, the peaks of each profile are selected by observation. However, automatic identification of peaks may be implemented by using a search algorithm to identify all local maxima in the significance score profiles.}

The calculated significance scores were compared against the mutual information used to select $\tau$ for uniform delay embedding. Two automated methods, MDOP \cite{nichkawde2013optimal} and PECUZAL \cite{kramer2021unified}, were also used as a comparison benchmark for non-uniform delay embedding. These calculations were done using implementations provided by the {\tt`DynamicalSystems.jl'} package \cite{Datseris2018}.

In addition to comparing the output lags from each method, it is of interest to find how different non-uniform embeddings impact performance in prediction tasks. For each of the non-uniform embedding methods, the first 2 dominant lags are taken to construct a 3D delay embedding. This is done by visually inspecting the $S_\mu(\tau), S_\sigma(\tau)$ profile in the case of SToPS. Time lags for PECUZAL and MDOP were taken as the first two detected timelags in their respective iterative procedures. 

The resulting delay embeddings used to train a simple 4 layer feed forward neural network consisting of 2 hidden layers, 1 input and 1 output layer. Each hidden layer has 128 nodes with a ReLU activation function with an overall network architecture of 3:128:128:3. 

The neural network is trained to output a one step prediction,

\begin{equation}
    \vec{x}(t+\delta t)=\text{NN}(\vec{x}(t)).
\end{equation}

The learning rate was set to 0.001 with a batch size of 512 and run for 30 epochs. In each instance, only the first half of the time series data is used for training. The second half is reserved for validation. Validation is done by calculating $n$-step freerun predictions. This is calculated by providing an initial condition and feeding back the neural network outputs $n$ times to get the final prediction (see Equation \ref{eq:freerun}), which is then used for calculating prediction error.

\begin{equation}
    \vec{x}(t+n \delta t)=\text{NN}^{(n)}(\vec{x}(t)).
    \label{eq:freerun}
\end{equation}

\section{Results and Discussion}
\subsection{Significant Times}
\subsubsection{Periodic Dynamics - Sum of Sines}
\begin{figure}
    \centering
    \includegraphics[width = 0.48\textwidth]{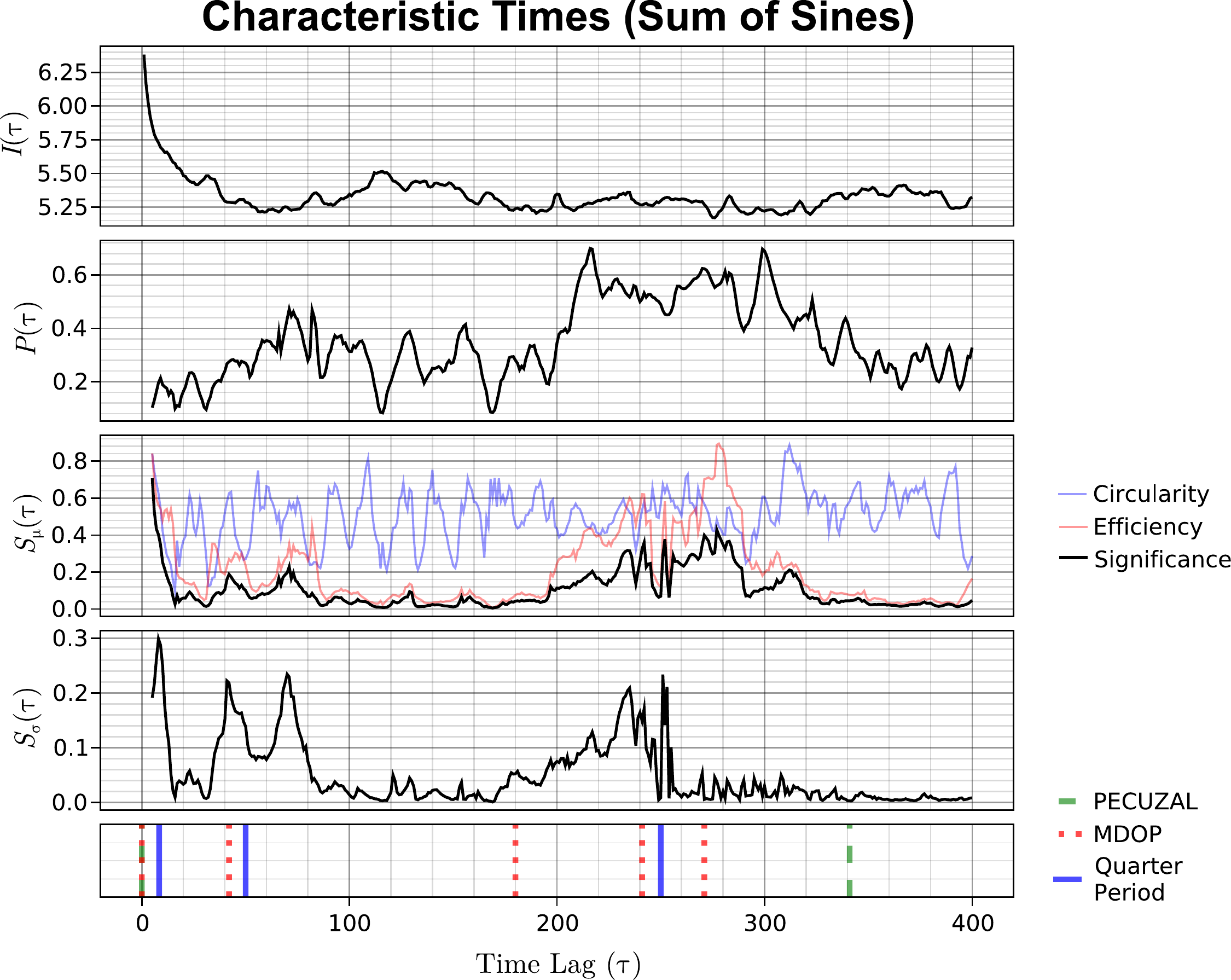}
    \caption{Comparison of various methods for estimating embedding delay for the sum of sines time series with component frequencies $\omega=(1, 5, 30)$. Top to bottom: mutual information, max persistence, significance score (black) with circularity (blue) and efficiency (red), standard deviation of significance score. \textcolor{edit_col}{Theoretical characteristic times based on one quarter of the period of component frequencies given by solid blue vertical lines.} Comparison non-uniform delays calculated from PECUZAL and MDOP given in green and red vertical lines. }
    \label{fig:SINES}
\end{figure}

The periodic sum of sines time series represents the case of periodic dynamics with multiple time scales. The component frequencies were selected as $\omega=\{1,5,30\}$, corresponding to 3 different characteristic times at lags $\tau=(250, 50, 8)$. The sum of sines times series was a major challenge for the baseline mutual information $I(\tau)$ measure (see Figure \ref{fig:SINES}), where the minima were only able to identify the highest frequency periodicity in the data. The max persistence $P(\tau)$ was also not useful in identifying any significant time lags from the time series. 

\textcolor{edit_col}{For automated non-uniform embedding methods, PECUZAL found a single time lag at $\tau=341$. MDOP returned four time lags at $\tau=(241,42,271,180)$, two of which are close to the expected lags $\tau=(250, 50, 8)$. The other remaining time lags detected by both methods did not bear any clear relation to the expected characteristic times.}

In contrast, the proposed significance score measures, $S_\mu(\tau)$ and $S_\sigma(\tau)$, both showed peaks around the expected time lags corresponding to approximately one quarter of the component periodicities (i.e. $T=1, 1/5, 1/30$) (see Figure \ref{fig:SINES}). However, the peak corresponding to the $\omega=5$ component (i.e. $\tau=45$) is not clear with two peaks occurring at nearby time lags instead. This anomaly may be because the spatial and time scales of the second frequency is not dynamically distinct enough from the large time scale. When calculating the 2D embedding, this can cause persistent strands to form spirals instead of circular holes at the characteristic times. We note that circularity $\alpha(\tau)$ and efficiency $\gamma(\tau)$ provide quite different profiles with the latter heavily influencing the shape of the resulting profile $S_\mu(\tau)$. \textcolor{edit_col}{A comparison of the phase space reconstructions between MDOP, PECUZAL and SToPS is provided in the Appendix.} 

\subsubsection{Chaotic Dynamics - Lorenz}
\textcolor{edit_col}{For Lorenz, the minima of mutual information yielded lags at $\tau\approx(40,150)$. This result matches closely with the maxima taken from the maximum persistence $P(\tau)$ profile. Of the two lags detected by PECUZAL  $\tau=(46,24)$, one was similar to a minima from the mutual information curve. Similarly, MDOP yielded 3 different lags $\tau=(148, 192, 37)$, one of which approximately matches the minima of the mutual information.}

\textcolor{edit_col}{The significance measures calculated using SToPS yielded 2 distinct maxima across $S_\mu(\tau), S_\sigma(\tau)$ with lags at $\tau=(41,130)$. The first lag is similar to calculated lags from PECUZAL and MDOP. However, the second lag $\tau = 130$ results in an overembedding of the time series and does not directly correspond to any lag output by either PECUZAL or MDOP. Closer inspection of the embedded time series at $\tau=120$ reveals the re-emergence of the lobes of the Lorenz attractor at $3/4$ of the period with boundaries created by multiple dense loops. Whilst this may produce well defined holes in the persistent homology near the the lobes, it does not efficiently utilise the points all the points in the sampled strand (i.e. $4\tau = 520$ points) (see Figure \ref{fig:spurious_lags}) and hence should not be classified as a characteristic time lag of the time series and is reflected in the much lower significance score. The usage of this time lag results in an overembedding of the time series. However, closer inspection of the sampled strands show that part of the increase in $S_\mu(\tau)$ is attributed to the hole formed from the spread of trajectories near the saddle point of the attractor.}

There is also an apparent correspondence between successive minima of mutual information and detected time lags in non-uniform embedding methods. This suggests that mutual information $I(\tau)$ may be useful for informing the selection of lags for non-uniform embedding. However, one difficulty is assessing if the lag of a minimum is within an acceptable embedding window $m\tau$ such that irrelevance is not a problem.

\begin{figure}
    \centering
    \includegraphics[width = 0.48\textwidth]{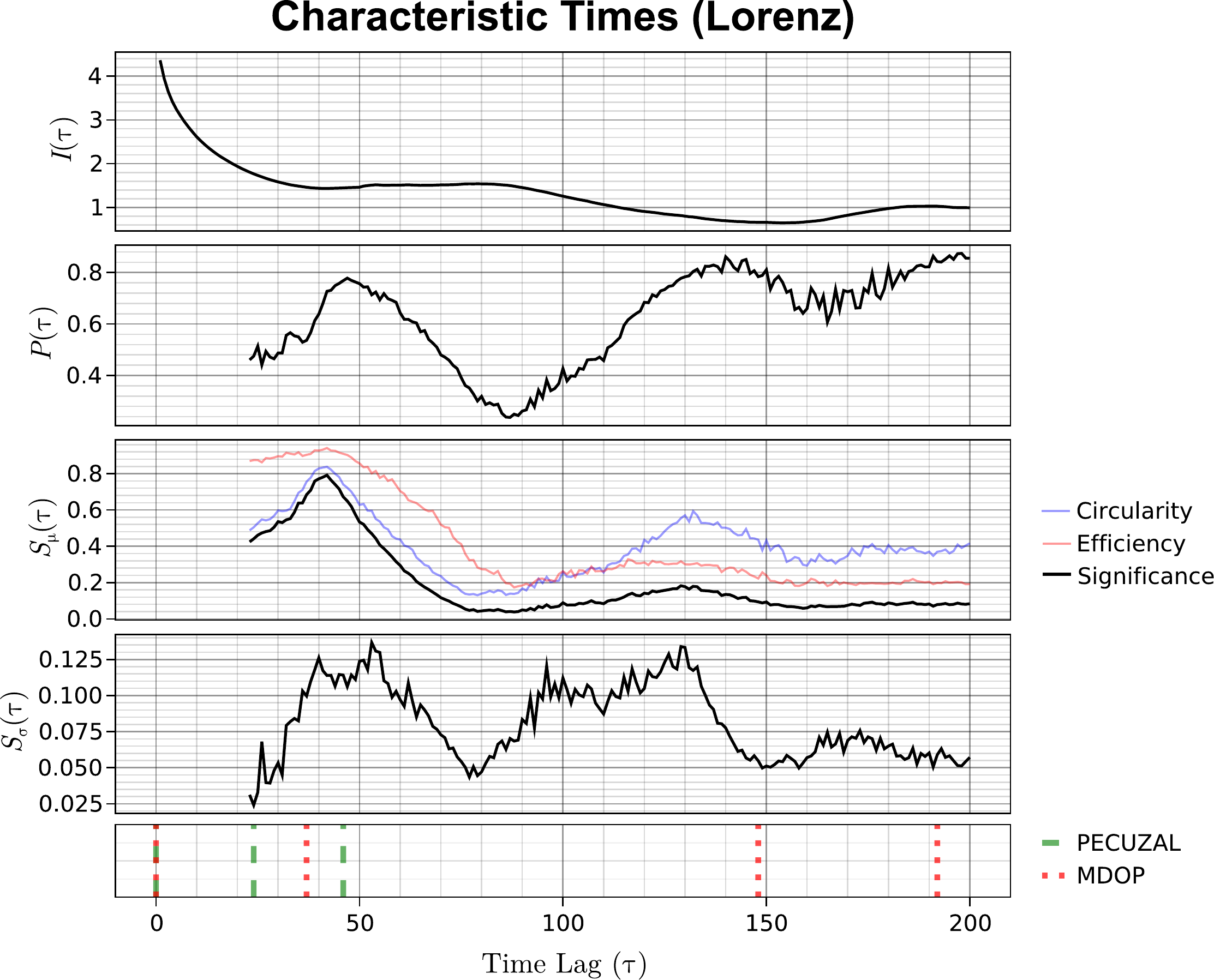}
    \caption{Comparison of various methods for estimating embedding delay for the chaotic Lorenz time series. Top to bottom: mutual information, max persistence, significance score (black) with circularity (blue) and efficiency (red), standard deviation of significance score. Comparison non-uniform delays calculated from PECUZAL and MDOP given in green and red vertical lines. }
    \label{fig:LORENZ}
\end{figure}

\subsubsection{Fast-Slow Dynamics - Lobster LP Neuron}
\begin{figure}
    \centering
    \includegraphics[width = 0.48\textwidth]{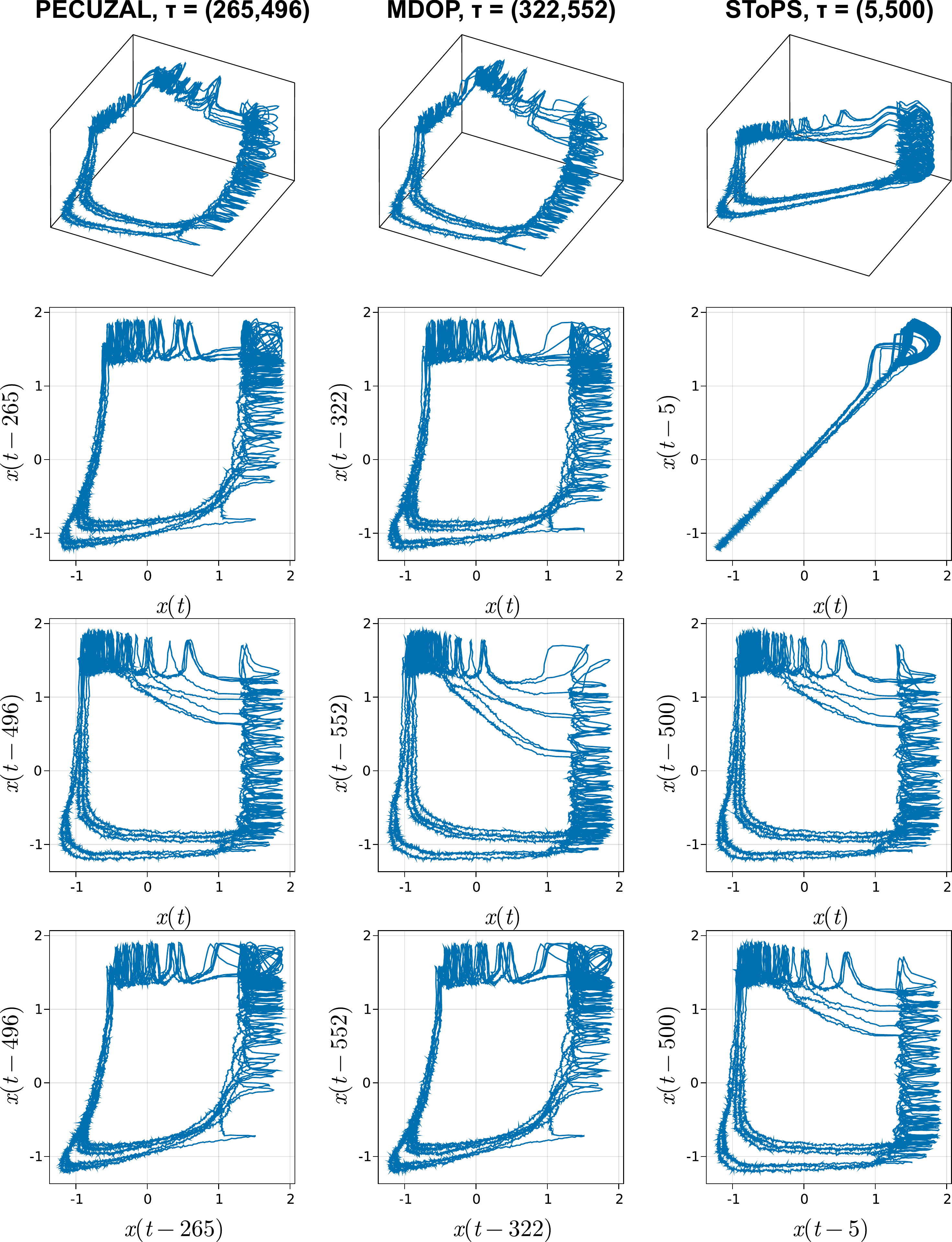}
    \caption{\textcolor{edit_col}{Comparison of 3D time delay embedding between PECUZAL, MDOP and SToPS with the the corresponding 2D projections. PECUZAL and MDOP time delays are selected from the first two detected delays. SToPS delays are visually selected from $S_\mu(\tau)$. The fast, small scale oscillations for spiking neurons are captured by SToPS but not by PECUZAL and MDOP.}}
    \label{fig:LOBSTER_Embedding}
\end{figure}

From observing the lobster neuron data (see Figure \ref{fig:timeseries}), it is possible to infer two dominant time scales corresponding to expected lags of approximately $\tau=(12,400)$. Temporal variations in these can be attributed to observational noise or the potentially chaotic dynamics. The uniform embedding measures of mutual information and max persistence were found to poorly in identifying lags. Although, max persistence $P(\tau)$ begins to quickly increase when approaching the expected characteristic lag $\tau=400$. This is unsurprising as a 2D embedding at those lags begins to unfold large orbits from the large time scale dynamics corresponding to the slow bursting phase of the neuron (see Figure \ref{fig:LOBSTER_Embedding}).

\begin{figure}
    \centering
    \includegraphics[width = 0.48\textwidth]{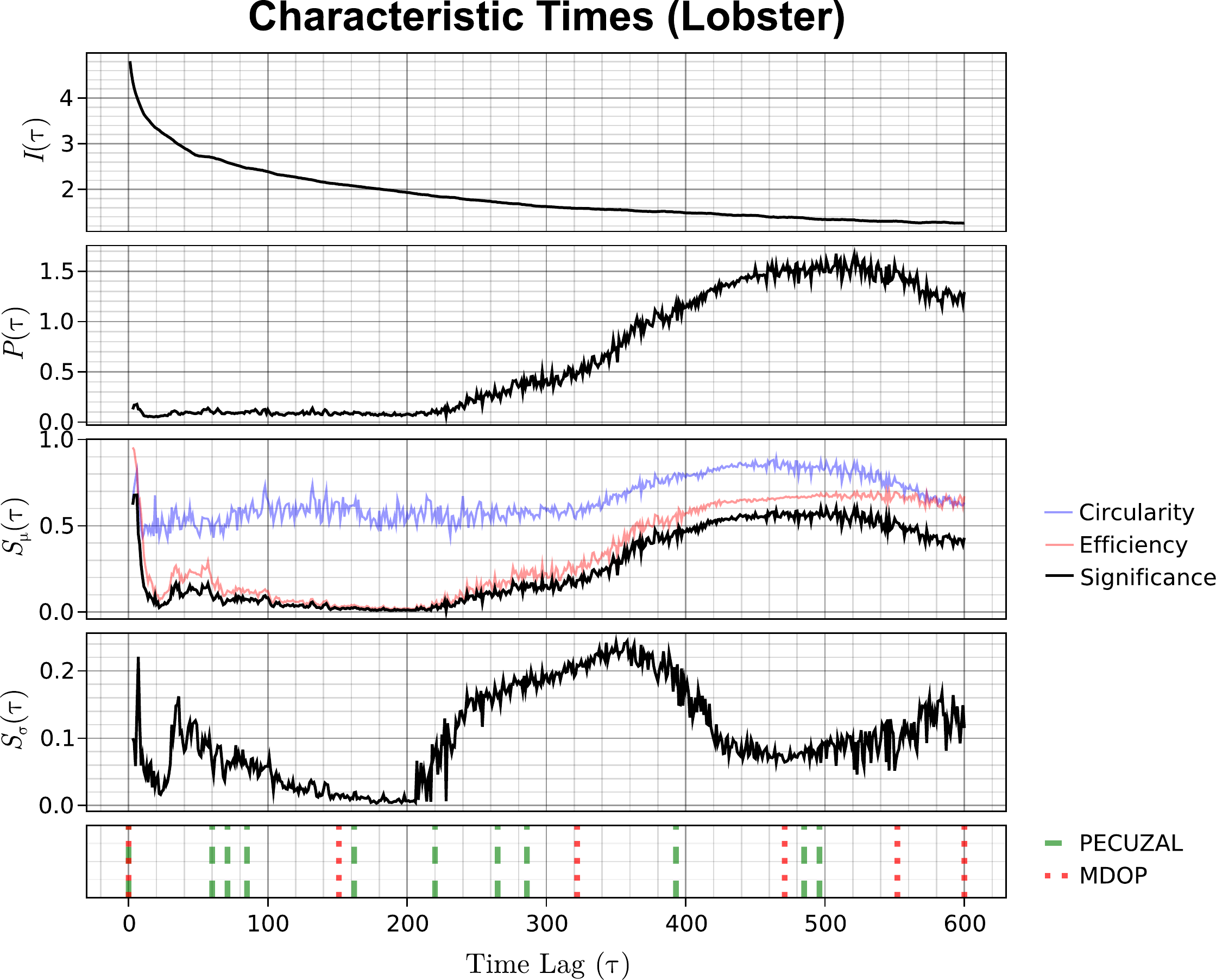}
    \caption{Comparison of various methods for estimating embedding delay for the chaotic Lorenz time series. Top to bottom: mutual information, max persistence, significance score (black) with circularity (blue) and efficiency (red), standard deviation of significance score.  Comparison non-uniform delays calculated from PECUZAL and MDOP given in green and red vertical lines. }
    \label{fig:LOBSTER}
\end{figure}

\textcolor{edit_col}{
For non-uniform embedding, both PECUZAL and MDOP yielded a large number of potential embedding lags. However, both methods failed to successfully recover the lags for the fast dynamics (see Figure \ref{fig:LOBSTER}). Additionally, apart from the lag at $\tau \approx 500$ corresponding to the slow dynamics, both PECUZAL and MDOP produce multiple additional time lags with no obvious explainable relation to the time scales of the data.}

The significance score using SToPS was able to retrieve the two main time lags present in the data at approximately $\tau=(5, 500)$. There appears to be a slight disagreement on the location of the larger time lag with $S_\mu(\tau)$ and $S_\sigma(\tau)$ producing slightly different lag times. \textcolor{edit_col}{The $S_\mu(\tau)$ also shows a small peak at approximately $\tau = (40,60)$. Further inspection into the representative homology of the sampled strands reveal that this is the result of multiple overlapping orbits from the fast spiking dynamics. However, this feature is captured by $\tau=5$ and the peaks at $\tau= (40,60)$ are an artefact of overembedded time series lying on similar orbits (see Figure \ref{fig:spurious_lags}).}

\textcolor{edit_col}{
Based on these results, only SToPS was able to detect both dominant time scales in the data. The lags from PECUZAL and MDOP are conditional on the selection of previously detected lags due to the iterative approach employed by the algorithm. In contrast, SToPS produces a single characteristic time spectrum $S(\tau)$ from which the significance of each potential lag can be assessed and selected independently. We visually compare the resulting 3D delay embeddings of these three methods in Figure \ref{fig:LOBSTER_Embedding}. \textcolor{edit_col}{Similar comparisons for the sum of sines and Lorenz data is provided in the Appendix.} For PECUZAL and MDOP, we select the first two non-zero delays in order of detection at $\tau=(265,496)$ and $\tau=(322,552)$ respectively. The lags for SToPS was selected visually from $S_\mu(\tau)$ and $S_\sigma(\tau)$ with lag times $\tau=(5, 500)$. From the projections of the reconstructed state space, we find that SToPS associates different projections with dynamics of different time scales. This results in an unfolding that is visually easier to interpret. In contrast, the PECUZAL requires a large number of dimensions before all dynamical components can be visually detected. Restricting the number of dimensions results in the fast dynamics being obscured at the expense of slow dynamics.}

\subsection{Freerun Predictions}
Freerun prediction errors with models constructed from different non-uniform embeddings were calculated for the non-periodic Lorenz and lobster LP neuron time series. In each case, a 3D delay embedding was constructed using the first two lags detected with PECUZAL and MDOP. \textcolor{edit_col}{For SToPS, the first two visually dominant maxima of the $S_\mu(\tau)$ and $S_\sigma(\tau)$ profiles were selected. In the case for Lorenz where only one relevant time lag exists (Figure \ref{fig:spurious_lags} and \ref{fig:LORENZ}) (i.e. $\tau = 41$), we use a uniform embedding scheme where the second lag is a multiple of the first.}

The prediction error for Lorenz was calculated using a \textcolor{edit_col}{25 step freerun prediction} with a trained feedforward neural network (see Section \ref{methodology}). For the lobster LP neuron, training was done on a subsampled data set that included every third point. This subsampling was done to reduce the training time. A 10 step freerun prediction was then evaluated with the subsampled data (i.e. equivalent to 30 step freerun prediction). The lower number of freerun prediction steps was used to allow a more accurate evaluation of the prediction performance within the fast neuron spiking regime of the time series. \textcolor{edit_col}{In both cases, the first half of the data was used for training and the second half used for testing and evaluation. Additional analyses for the freerun prediction with non-subsampled data is provided in the Appendix.}

The error of each prediction was calculated as the magnitude of the error between the predicted delay vector $n$-steps ahead,

\begin{equation}
    E(t)=||\text{NN}^{(n)}(\vec{x}(t))-\vec{x}(t+n)||.
\end{equation}

\begin{figure}
    \centering
    \includegraphics[width = 0.48\textwidth]{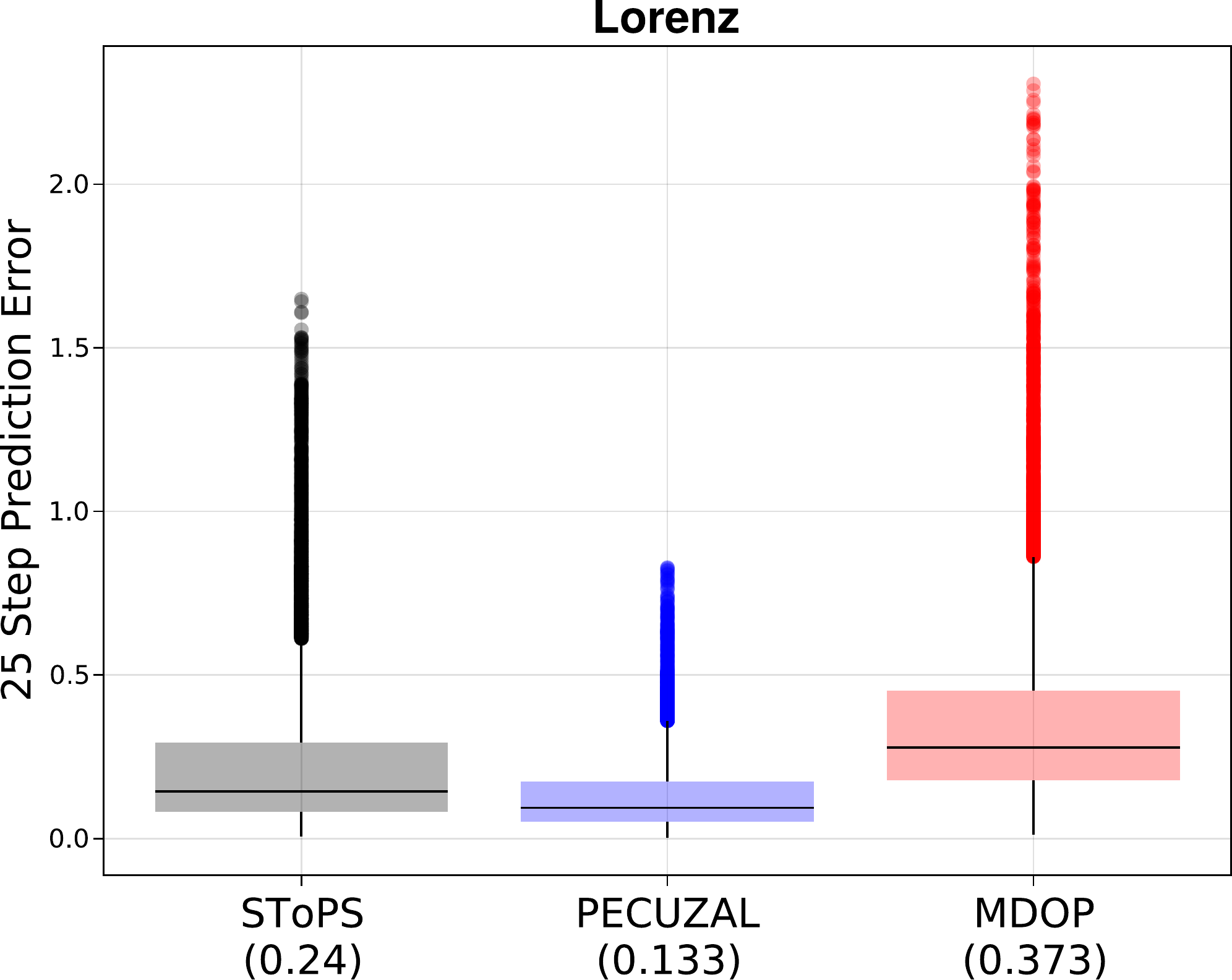}
    \caption{Distribution of the 25 step freerun prediction error for the Lorenz time series over 25000 steps using a neural network trained on the 3D delay embedding for one step prediction.}
    \label{fig:lorenz_pred}
\end{figure}

\begin{figure}
    \centering
    \includegraphics[width = 0.48\textwidth]{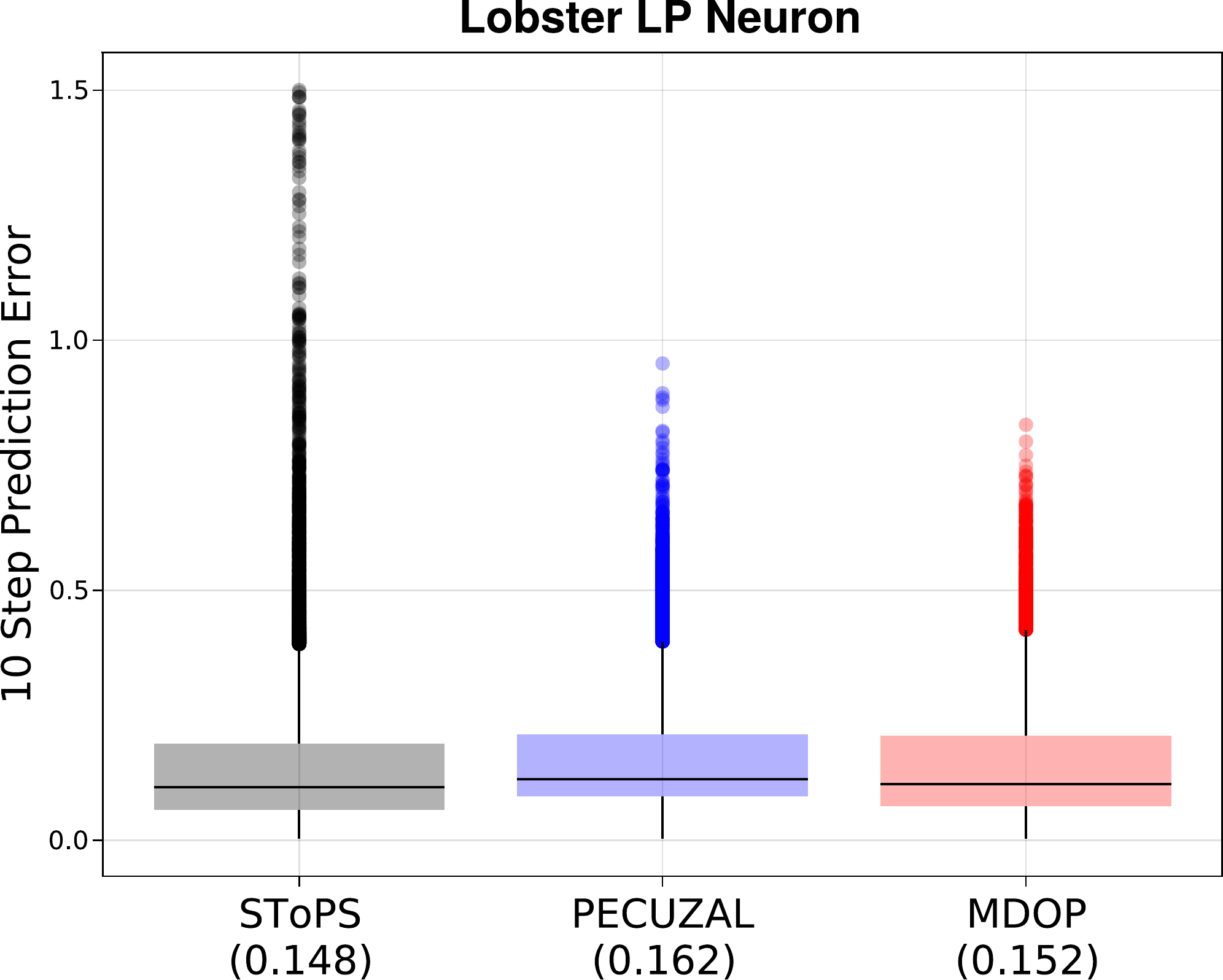}
    \caption{Distribution of the 10 step freerun prediction error for the lobster LP neuron time series over the entire subsampled data set (approximately 200000 steps) using a neural network trained on the 3D delay embedding for one step prediction.}
    \label{fig:lobster_pred}
\end{figure}

The resulting distributions of $E(t)$ for both time series with the three different embedding methods are given in Figures \ref{fig:lorenz_pred} and \ref{fig:lobster_pred}. \textcolor{edit_col}{For Lorenz, we see that SToPS (persistent strands) provides a mean error in between PECUZAL and MDOP. In the experimental data case (lobster LP neuron), SToPS outperforms both measures with a lower error. These findings are also reflected in the corresponding medians in both cases. The median prediction error in the Lorenz case was $0.106$ (SToPS), $0.122$ (PECUZAL) and $0.112$ (MDOP). The lobster LP neuron median prediction errors were $0.106$ (SToPS), $0.122$ (PECUZAL) and $0.112$ (MDOP). Additionally, SToPS shows an error distribution with a heavier tail for freerun predictions with the lobster LP neuron compared to PECUZAL and MDOP. Despite the potentially lower prediction error, we note that this improvement is not significant and should not be the targeted benefit of SToPS. Instead, we argue that the main advantage of SToPS is that it provides lags that are explainable in the context of the observed dynamics of the time series.}


\begin{figure}
    \centering
    \includegraphics[width = 0.45\textwidth]{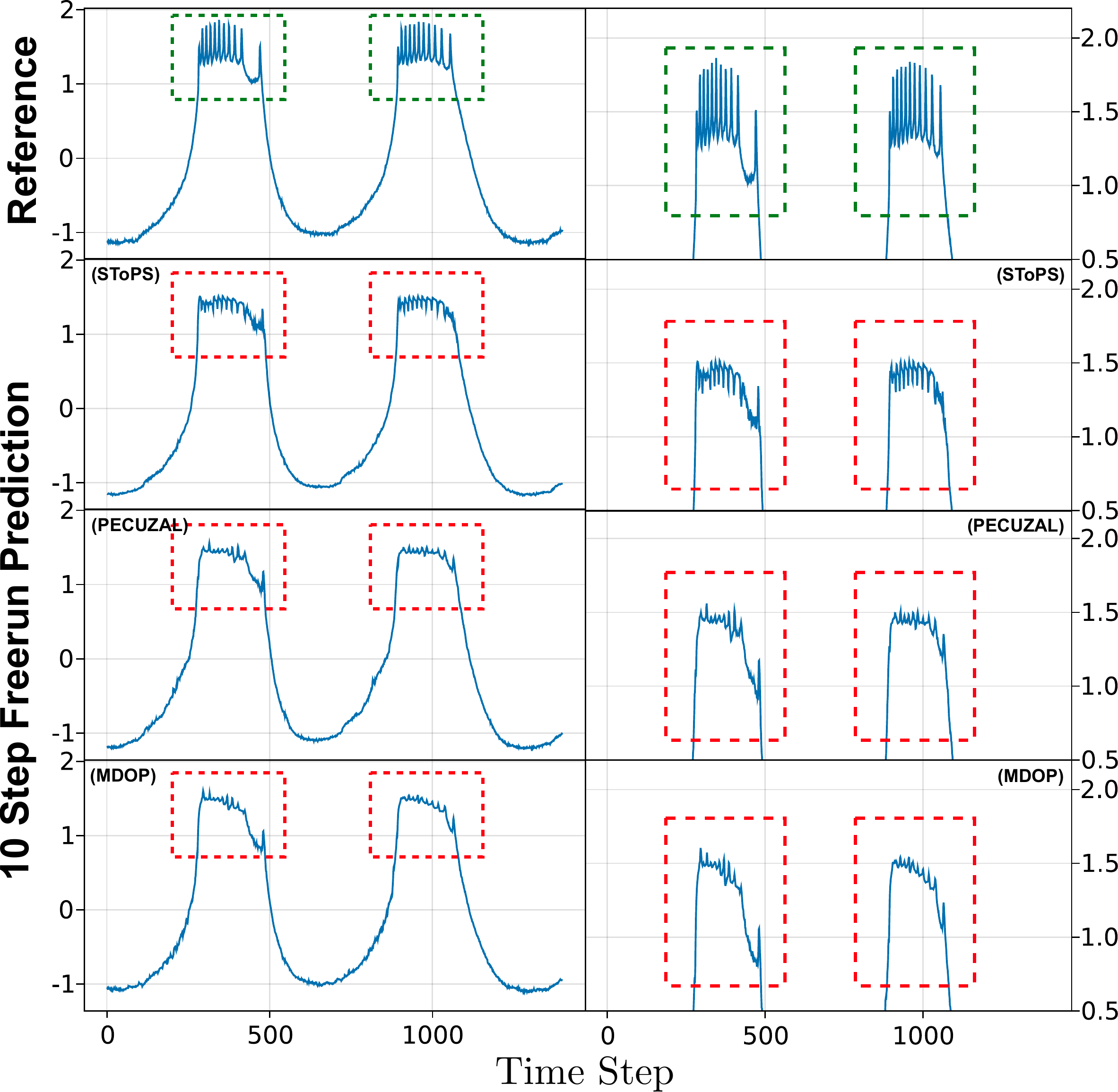}
    \caption{Comparison of 10 step freerun predictions of the lobster LP neuron time series across various embedding methods. The persistent strands method appears to be able to better detect and replicate the fast neuron spiking dynamics. Magnitude of predictions are normalised with zero mean and unit standard deviation.}
    \label{fig:lobster_pred_timeseries}
\end{figure}

One advantage of using SToPS is the deliberate inclusion of both fast and slow time scales within the embedding of the time series. Therefore, it is expected that models trained on this embedding should be able to better resolve the fast dynamics that would have otherwise be missed if larger delays were selected. This is verified in Figure \ref{fig:lobster_pred_timeseries} where SToPS produced a 10 step freerun prediction that is able to better replicate the small scale, fast spiking dynamics characteristic of the neuron time series. This is in contrast with models trained on the same number of lags from PECUZAL and MDOP where freerun predictions are not able to capture the same level of detail in the spiking dynamics. \textcolor{edit_col}{We note that this advantage is not as apparent when the full data set without subsampling is used to train a model for prediction. However, it was found that the prediction errors were slightly lower for SToPS in this case.}

\textcolor{edit_col}{
\section{Noise Effects}
One often cited benefit of persistent homology is its robustness to noise \cite{adams2017persistence,turkevs2021noise}. Because only geometric features (i.e. holes) are tracked, as long as the magnitude of the noise does not destroy the underlying structure of the embedded time series, the calculate homologies should be stable. To test this property in our method, we repeat analyses with the Lorenz data set for varying signal to noise ratios and observe the changes in the resulting significance score profile $S_{\mu}(\tau)$. Five noise levels of additive Gaussian noise with signal-noise ratios $SNR = (1000,500,100,50)$ were applied after normalising the input Lorenz time series with zero mean and unit standard deviation. Results are shown in Figure \ref{fig:noise}.
\begin{figure}
    \centering
    \includegraphics[width = 0.48\textwidth]{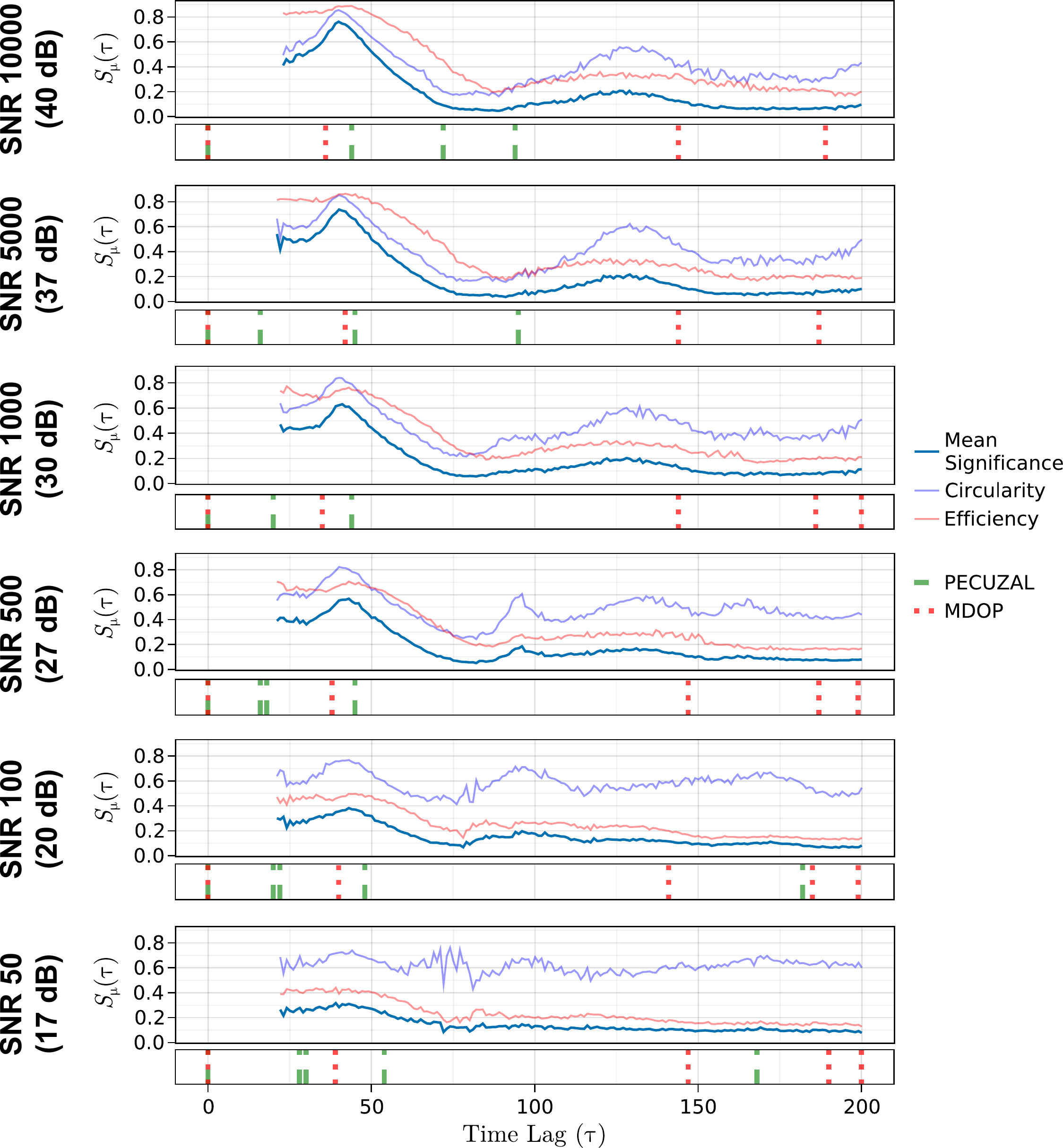}
    \caption{Comparison of significance score $S_\mu(\tau)$ against PECUZAL and MDOP lags at various additive noise levels. SToPS significance score is relatively stable with the peak locations being the same for all tested noise levels. PECUZAL and MDOP are stable for low levels of noise, but gradually drift for increasing noise levels.}
    \label{fig:noise}
\end{figure}
We found that each of the profiles was relatively stable and robust to increasing noise levels with little to no change the in the observed significant lags. For high levels of noise ($>20$ dB) the spurious lag artefact at $\tau=120$ disappeared. This is likely due to the effects of noise across multiple trajectories along a similar orbit destroying the underlying homology structure.}

\textcolor{edit_col}{
Similarly, PECUZAL and MDOP revealed similar lags for low noise levels ($<37$ dB). However, higher levels of noise resulted in gradual drifts in both methods. The lags from PECUZAL also differed significantly for small $\tau$ for noise levels above 30 dB. In contrast, MDOP was relatively unaffected by noise and drift effects were relatively small. Additionally, MDOP was also the only method able to produce any lag predictions for signal to noise ratios below 17 dB. Both SToPS and PECUZAL were unable to produce any results.}

\textcolor{edit_col}{
\section{Computation Complexity}
The computation of persistent homology is not parsimonious and still suffers from poor computational scaling. Whilst computation of persistent homologies for data sets with few points is relatively quick, the computation time grows rapidly for even moderately sized data sets. This presents a challenge of the SToPS method as it requires the computation of persistent homology of multiple sampled strands across a large collection of time series. From the proposed algorithm, the time complexity is approximately:}
\begin{equation*}
    % \mathcal{O}_{SToPS} \approx \mathcal{O}(N_\tau N_{s} (\mathcal{O}_{PH}(N_p)+N_p^2+4N_p + N_p\log(N_p)),
    \mathcal{O}_{SToPS} \approx \mathcal{O}(N_\tau N_{s} \phi_{PH}(N_p)),
\end{equation*}
\textcolor{edit_col}{
 where $N_\tau, N_s$ and $N_p$ are the number of lags tested, number of sampled strands and strand length respectively and $\phi_{PH}(N_p)$ is the time complexity for the persistent homology computation. This value varies depending on the implementation of the persistent homology algorithm and the type of filtration used (e.g. \v{C}ech, Rips, Delaunay etc.). The Vietoris-Rips filtration has a simplicial complex size $K$ that scales exponentially with $2^{\mathcal{O}(N_p)}$ in current formulations \cite{otter2017roadmap}. The computation of the Vietoris-Rips filtration can be split into two phases. However, calculating the computational complexity bounds is not straightforward \cite{zomorodian2010fast}. We note that there is ongoing work aimed at improving and optimising the persistent homology algorithm, which has resulted in significant gains in performance \cite{bauer2021ripser,otter2017roadmap}.}

\begin{table*}
\caption{\label{tab:table1} \textcolor{edit_col}{Computation times for PECUZAL, MDOP and SToPS on the Lorenz and lobster LP neuron time series with 25000 steps with varying maximum number of candidate lags (i.e. $\tau \in [1,\tau_{max}]$). Computation times are given in seconds. The computed lags for PECUZAL and MDOP are provided. SToPS does not have computed lags as results are given as a significance score profile}}
\begin{ruledtabular}
\begin{tabular}{ccccc}
% & & & Computation Time (s) and Lags & \\
% \hline
Time Series & $\tau_{max}$ & PECUZAL & MDOP & SToPS \\
\hline
Lorenz & 50 & 6.0 & 0.19 & 148\\
& & (45,23) & (42,23) & \\
& 100 & 10.2 & 0.23 & 982\\
& & (45,23) & (100,62) & \\
& 150 & 15.8 & 0.28 & 1689\\
& & (45,23) & (147,40) & \\
& 200 & 21.2 & 0.42 & 2396\\
& & (45,23) & (147,189,40) & \\
Lobster & 50 & 7.1 & 0.39 & 95.9\\
& & (50,24) & (50,31,11,44) & \\
& 100 & 46.9 & 0.48 & 881\\
& & (100,63,81,11,35) & (100,51,19,81) & \\
& 300 & 112 & 0.82 & 3376\\
& & (300,251,114) & (300,149,55,245) & \\
& 600 & 2705 & 1.62 & 7266\\
& & (515,304,132,505,32,20,253,394,81,232,55) & (544,274,599,114,423) & \\

% PECUZAL & 23.5\\
% MDOP & 1.8\\
% SToPS & 3178.1\\
\end{tabular}
\end{ruledtabular}
\end{table*}

\textcolor{edit_col}{
A comparison of the approximate run times for SToPS, PECUZAL and MDOP are provided in Table \ref{tab:table1}. Computation was done on a Ryzen 7 4800HS with 16GB of RAM using only one thread in order to allow results to be comparable. Despite having a much longer computational time compared to other methods, we note that SToPS allows multiple lags to be considered in parallel as the significance of each potential time lag is evaluated independently. This is in contrast with PECUZAL and MDOP where the iterative embedding cycles approach is used and each new lag is selected conditional upon previously selected lags. However, there is a new method proposed by Kr{\"a}mer et. al. based on Monte Carlo tree search that attempts to tackle this problem \cite{kraemer2022optimal}. The flexibility to assess embedding lags independently may provide large gains in computational speed where the computation for multiple lags may be distributed across multiple threads. Our current implementation of the algorithm does not yet provide support for this.}

\textcolor{edit_col}{The selection of new lags conditional on previously selected lags also means that the results of PECUZAL and MDOP are not robust to changes in the maximum allowable lag $\tau_{max}$. Changes to the range of potential lags can affect the order of the selection of future lags if a new candidate lag that better optimises the objective statistic is introduced. For complex data, this can result in widely varying results as the maximum lag changes. The presence of noise in the data also results in different lags as shown by the drifting lags found in Figure \ref{fig:noise}. Calculated lags for the lobster LP neuron time series were found to vary even across trials with different realisations of identically distributed noise.}

\textcolor{edit_col}{We also note that the computation time for these embedding algorithms can vary depending on the complexity of the input time series. For example, successive embedding cycles in PECUZAL typically increase in computation time. For complex or noisy time series, PECUZAL may produce multiple lags of varying significance. This is seen by the jump in estimated lags for the lobster data between max lags of 300 and 600. In contrast, the SToPS algorithm assesses the significance of each score and should grow linearly proportional to the number of lags. An exception is for small lags $\tau$ where strands are too short and have no holes' homology to track.}

\section{Conclusions}

One of the aims of this paper is to provide an overview of the embedding fundamentals and review existing methods for optimising embedding parameters. We first provide in Section \ref{embedding_case} a rough overview on the fundamental concepts of embedding theorems and its applications in the context of time series analysis. A simple case for the usage of embedding in time series prediction tasks is also given for new or uninitiated readers. In this paper, we focus on the problem of identifying good embedding parameters, specifically on the selection of embedding lags in non-uniform embedding. An overview on the considerations when selecting embedding parameters is provided in Section \ref{embedding_considerations} followed by a comprehensive review of various uniform delay embedding methods in Sections \ref{uniform_delay_embeddings}-\ref{optimising}. 

We argue that a non-uniform embedding approach provides more flexibility in reconstructing fast-slow dynamical systems. Following this, an overview of existing methods that attempt to automate the selection of non-uniform embedding lags is provided in Section \ref{non-uniform_embedding}. However, whilst many of these automated non-uniform embedding methods reliably return a collection of lags, they do not necessarily agree with each other or provide a satisfactory dynamical explanation for their selection. Furthermore, due to the iterative process used to select delays, the choice of each subsequent delay is conditional on previous selections.

We propose in Sections \ref{Persistent_Strands} and \ref{SToPS}, a new method of selecting non-uniform embedding lags, SToPS, that aims to produce lags that have more dynamical explainability and are independently selected. SToPS utilises persistent homology to detect loops formed by 2D delay embeddings of sampled windows of the time series, which we call `persistent strands'. This is done over multiple different lengths of sample windows to produce a characteristic time spectrum $S_\mu(\tau)$ where larger values of the significance score correspond to time scales that are dynamically significant (i.e. they relate to some notion of periodicity in the time series). The structure of each persistent strand loop is characterised by two quantities, circularity and efficiency, which are combined to give the significance score $S_i(\tau)$. Selection of time lags $\tau$ are done based on the mean and standard deviation profiles of the significance score $S_\mu(\tau), S_\sigma(\tau)$.

The SToPS embedding method is tested on three different classes of time series: periodic (sum of sines), chaotic (Lorenz) and fast-slow (lobster LP neuron). In all cases, SToPS was found to detect dynamically explainable time scales that were not reflected in other reference non-uniform embedding methods PECUZAL and MDOP. Additionally, SToPS was found to outperform PECUZAL and MDOP in identifying dominant time scales for the lobster LP neuron where fast-slow dynamics are present.

The impact of each different embedding method on the time series prediction performance was also tested. Embedded time series were used to train a one step neural network predictor. It was found that the resulting models performed similarly across all embedding methods for both the Lorenz and lobster LP neuron time series. However, freerun predictions of the lobster LP neuron time series with models trained using SToPS embedding lags were found to be able to replicate the fast neuron spiking dynamics better than reference embedding methods. \textcolor{edit_col}{We also provide a brief discussion and analysis on the computational efficiency of SToP, as well as its robustness to noisy input data.}

Overall, whilst the performance of SToPS is only marginally better than the existing non-uniform embedding methods PECUZAL and MDOP, we argue that SToPS provides lags that are more dynamically explainable compared to its counterparts.  \textcolor{edit_col}{The assessment individual time lags also allows the method to be applied to multivariate time series by considering each component as an independent scalar time series and identifying their respective time lags. This may then be used to construct a delay vector that utilises all components of the time series. The performance of SToPS in this context has not yet been tested and presents as an avenue of further research.} Additionally, the independent selection of lags via the characteristic time spectrum provides a clearer picture of the relative importance of each time lag when compared to existing iterative methods for automated non-uniform delay embedding where an explicit collection of lags is provided. 

However, the pursuit of more dynamically explainable delay lags introduces a level of subjectivity in the interpretation of the characteristic time spectrum. Nevertheless, we argue that the focusing on selecting dynamically relevant and explainable delay lags is potentially a more meaningful approach to constructing models that are more relatable to observed system dynamics. This advantage is especially evident in systems with multiple disparate time and spatial scales, as demonstrate by the dynamics of a the Lobster LP neuron. Therefore, we propose that dynamical relevance and explainability to be a key additional consideration in the future development of time delay embedding methods.

% The formulation of the significance score in SToPS also combines both dynamical and topological arguments in quantifying the significance of each time lag. 

\begin{acknowledgments}
D.C.C. is supported by the Australian Research Council through the Centre for Transforming Maintenance through Data Science (grant number IC180100030), funded by the Australian Government. D.W. and M.S. are supported by the Australian Research Council (grant number DP200102961). S.D.A is supported by the Forrest Foundation. E.T. is supported by the Robert \& Maude Gledden Foundation and the A.F. Pillow Applied Mathematics Trust. We would also like to thank H. Abarbanel and his team for providing us access to the experimental LP neuron data. 
\end{acknowledgments}

\section*{Data Availability Statement}

All the relevant data cited within this paper are publically available. For artificial data, generation methods have been provided in this paper. Other experimental data that supports the findings are available from the corresponding author at reasonable request. A Julia script containing the SToPS analysis method is also available at: \url{https://github.com/eugenetkj98/SToPS_Public} 

\bibliography{references}

\appendix
\section{Non-Uniform Embedding Profiles}
\begin{figure}

    \subfloat[Sum of 3 sine terms with $\omega=\{1, 5, 30\}$, $\phi=\{0,0.25,0.75\}$ and $A=\{1, 0.5, 0.2\}$.]{\includegraphics[width = 0.45\textwidth]{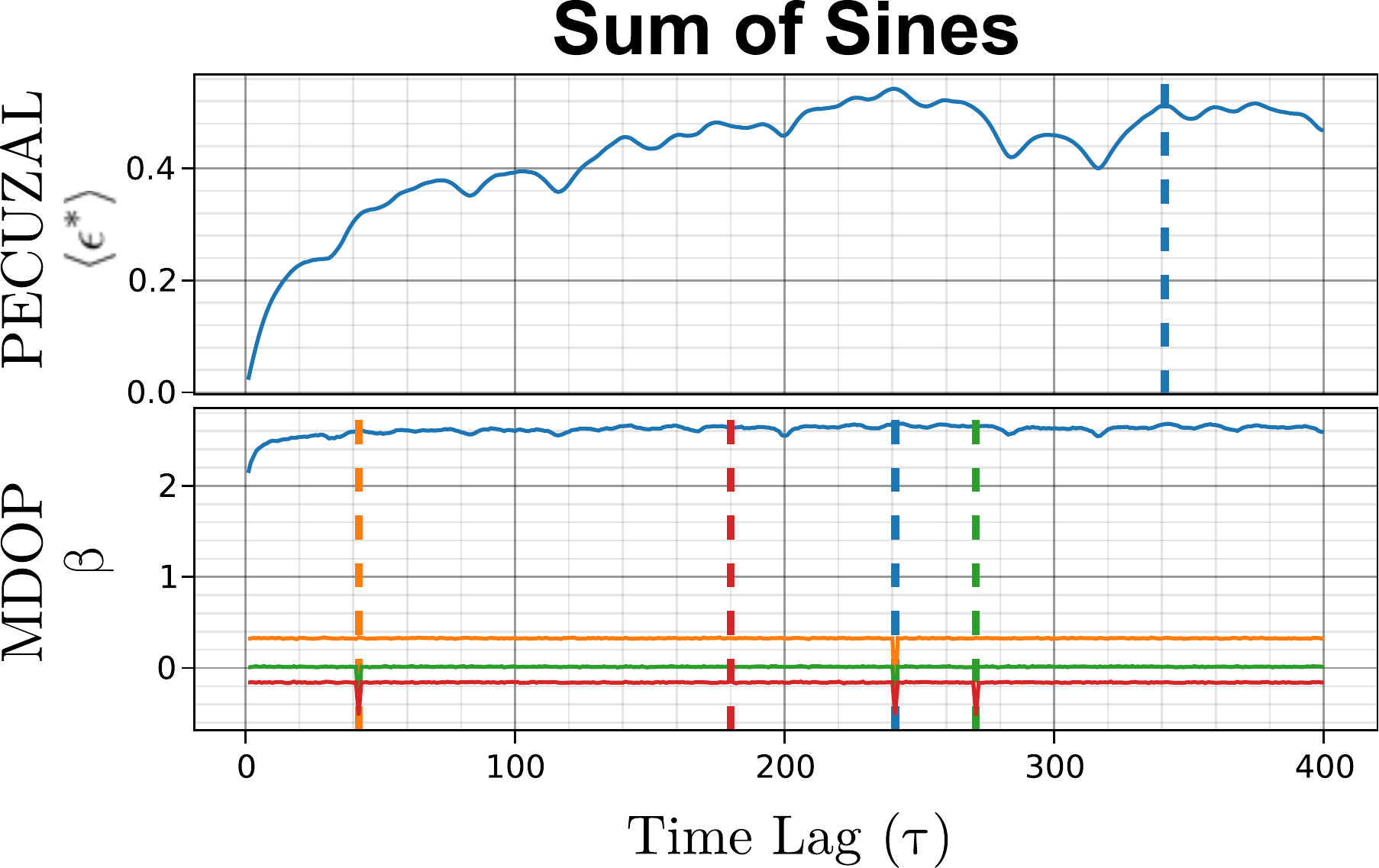}}\\
    
    \subfloat[First component of Lorenz.]{\includegraphics[width = 0.45\textwidth]{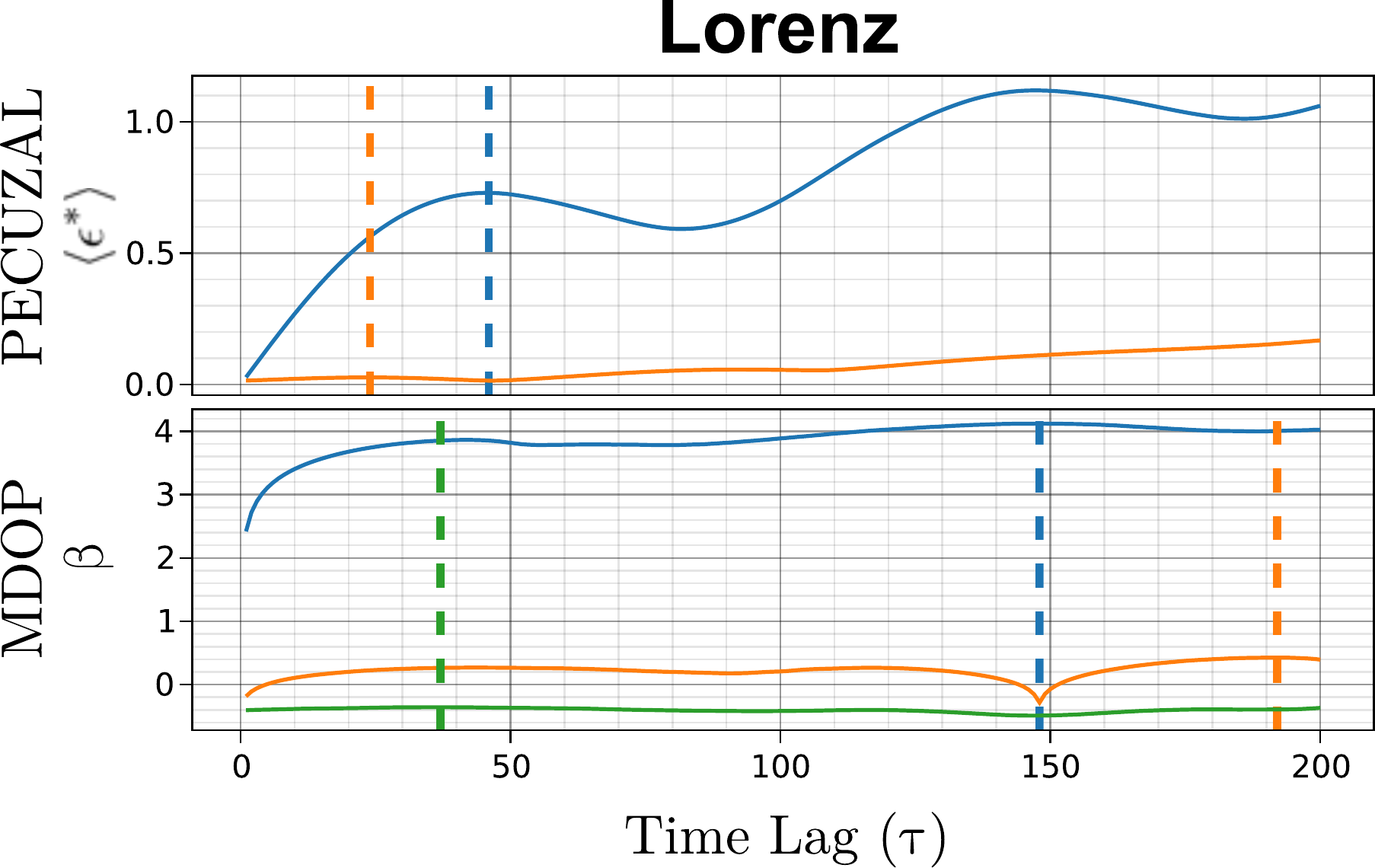}}
    
    \subfloat[Lobster LP neuron time series]{\includegraphics[width = 0.45\textwidth]{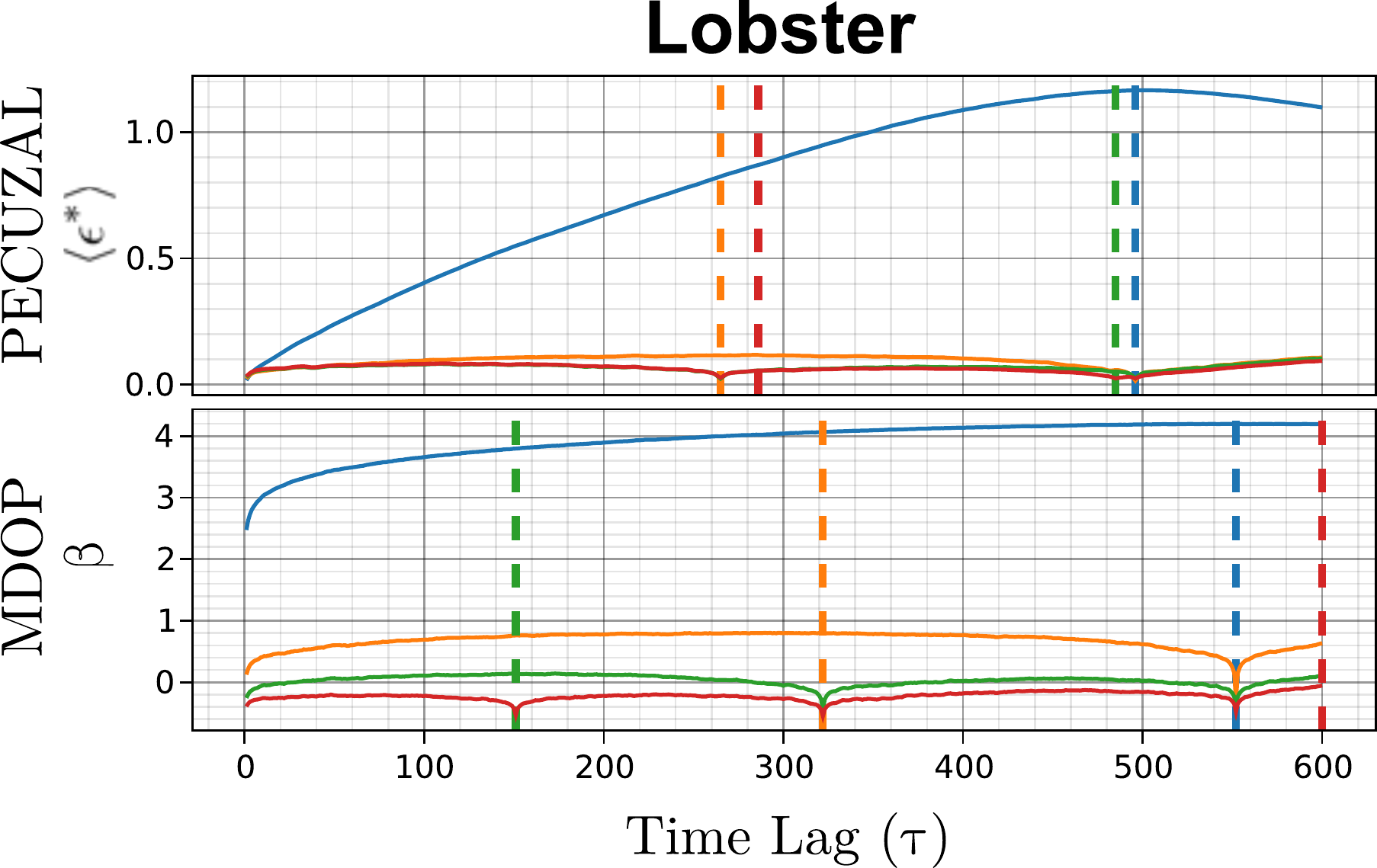}}
    
    \caption{Resulting embedding profiles for the calculated statistics in the comparison automated non-uniform embedding methods. }
    \label{fig:Embedding_Profiles}
\end{figure}

\textcolor{edit_col}{A collection of profiles of the calculated statistics used to select embedding lags the automated non-uniform embedding algorithms PECUZAL and MDOP are given in Figure \ref{fig:Embedding_Profiles}. Relevant statistics are the continuity statistic $\langle \epsilon^* \rangle$ for PECUZAL and the $\beta$ statistic for MDOP. For PECUZAL, the $\tau$ lags corresponding to local maxima in $\langle \epsilon^* \rangle$ are used candidate lags. The lag that results in the largest decrease in the $L$-statistic is chosen as the final embedding lag in each embedding cycle. A similar process is done for MDOP, but the global maxima of $\beta$ is chosen instead. A termination criterion based on the false nearest neighbour (FNN) statistic is used in conjunction.}

\section{Freerun Prediction of Non-Subsampled Data}
\begin{figure}
    \subfloat[Comparison of freerun prediction trajectories.]{\includegraphics[width = 0.48\textwidth]{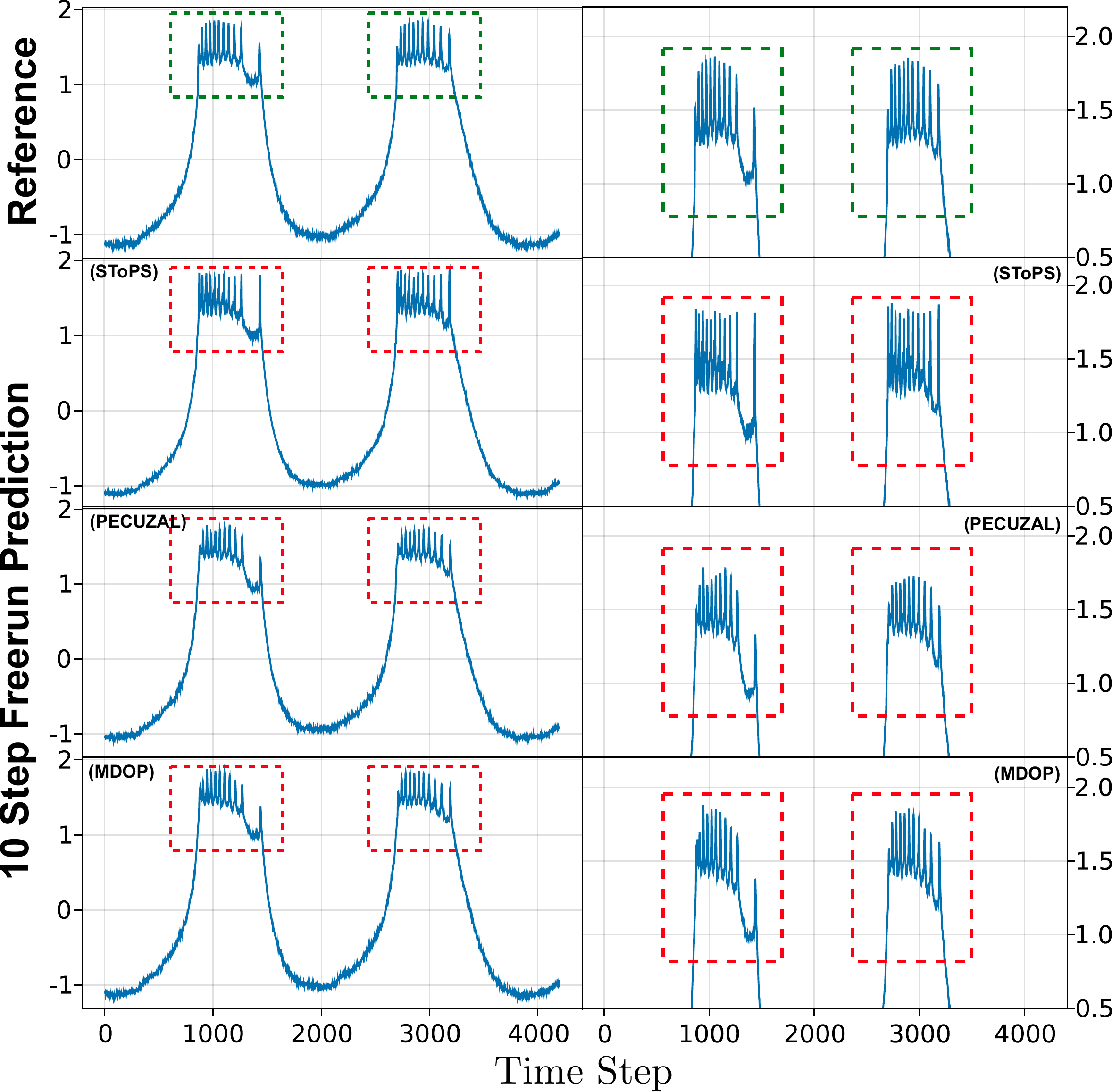}}\\
    
    \subfloat[Distribution of 10 step mean prediction error]{\includegraphics[width = 0.48\textwidth]{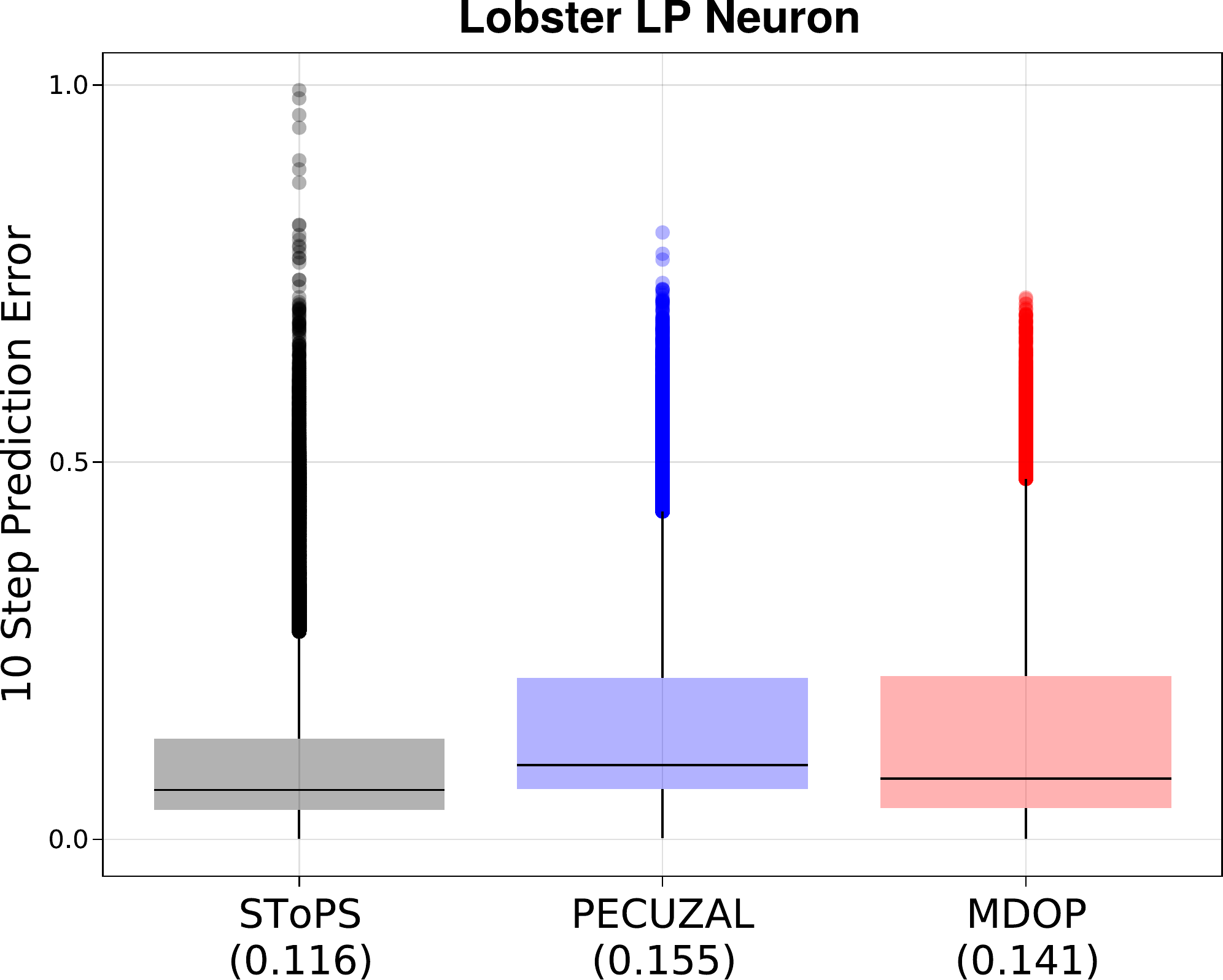}}

    \caption{10 step freerun prediction results for a neural network model trained on non-subsample data.}
    \label{fig:FreerunPred_10STEP_SS1}
\end{figure}

\begin{figure}

    \subfloat[Comparison of freerun prediction trajectories.]{\includegraphics[width = 0.48\textwidth]{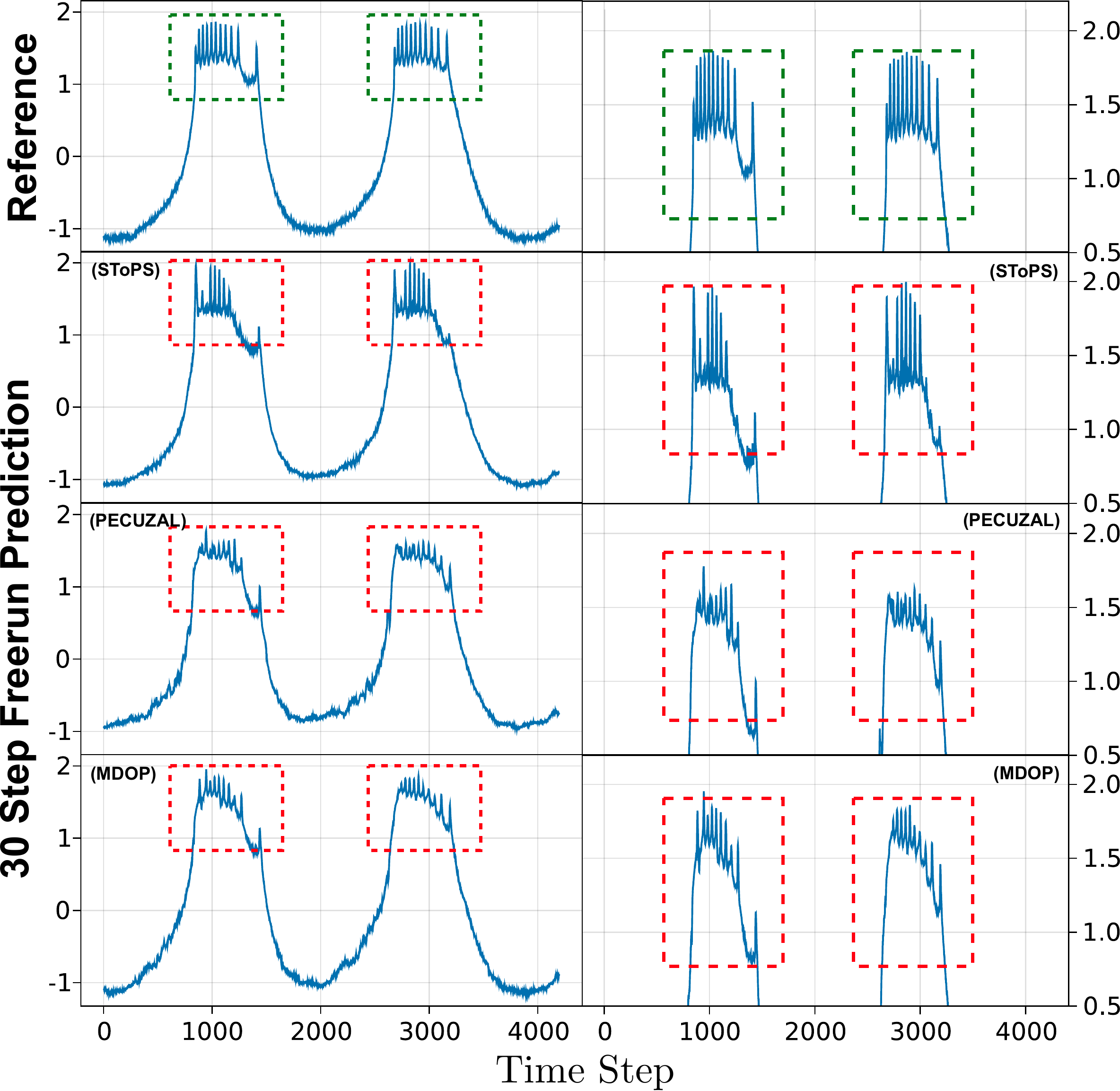}}\\
    
    \subfloat[Distribution of 30 step mean prediction error.]{\includegraphics[width = 0.48\textwidth]{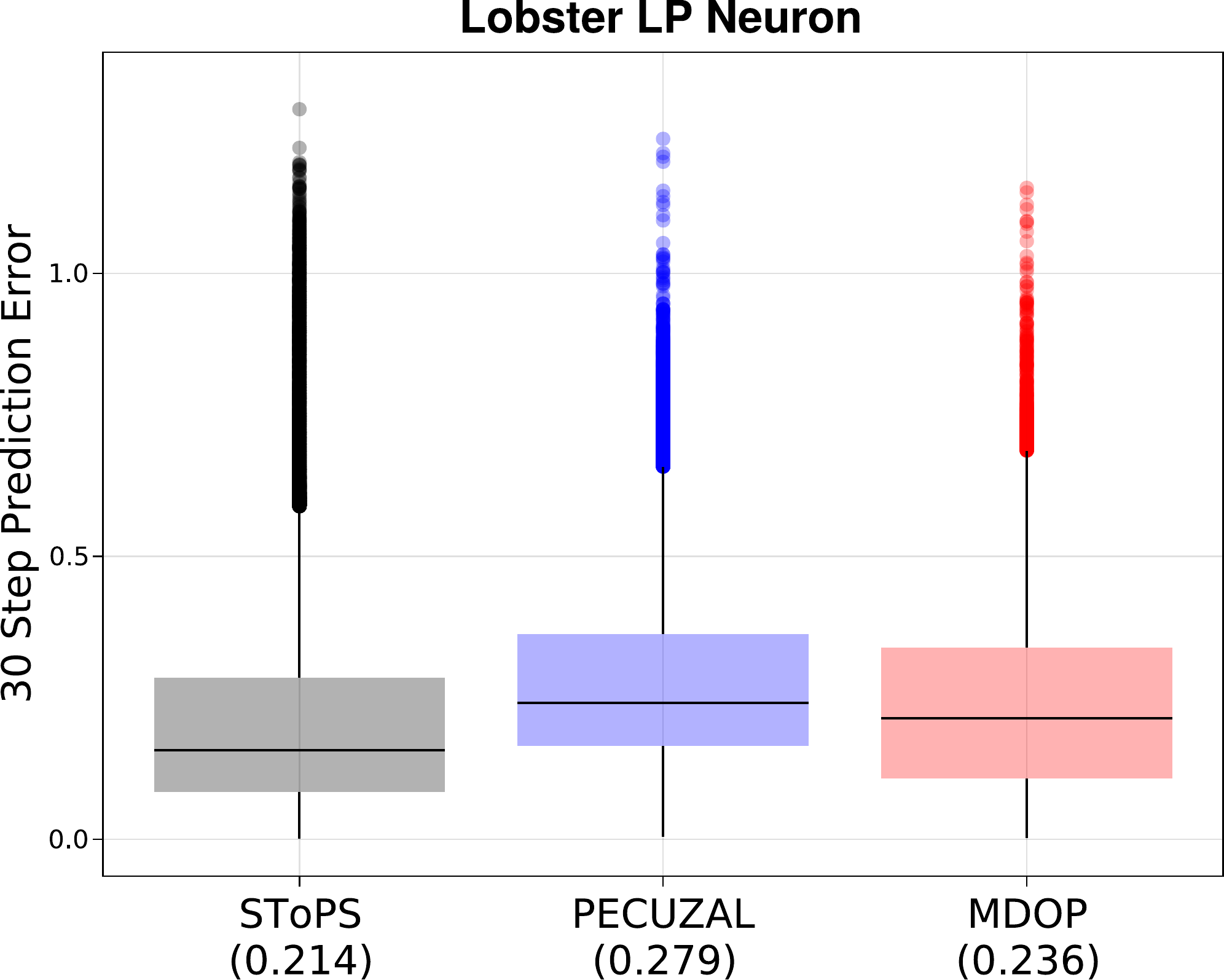}}

    \caption{30 step freerun prediction results for a neural network model trained on non-subsample data.}
    \label{fig:FreerunPred_30STEP_SS1}
\end{figure}

\textcolor{edit_col}{
This section contains results of additional freerun prediction tests for models trains on the full lobster LP neuron data set without any subsampling. This is theoretically an easier task as the model is provided with a larger amount of data and smaller with smaller magnitude predictions in each step. Predictions horizons of 10 steps (Figure \ref{fig:FreerunPred_10STEP_SS1}) and 30 steps (Figure \ref{fig:FreerunPred_30STEP_SS1}) were done. The latter's prediction horizon is equal to the 10 step prediction horizon models trained on the subsampled data. In the 10 step non-subsampled case, SToPS is found to yield a lower prediction error than PECUZAL and MDOP. The replication of the fast spiking dynamics is relatively similar between all methods. Similar results were found for the 30 step prediction case, However, the replication of the fast spiking dynamics is poorer than the 10 step prediction cases likely due to the accumulation of errors in successively predicted values.}

\section{Phase Space Reconstruction Comparisons}

\begin{figure*}
    \centering
    \includegraphics[width = 0.9\textwidth]{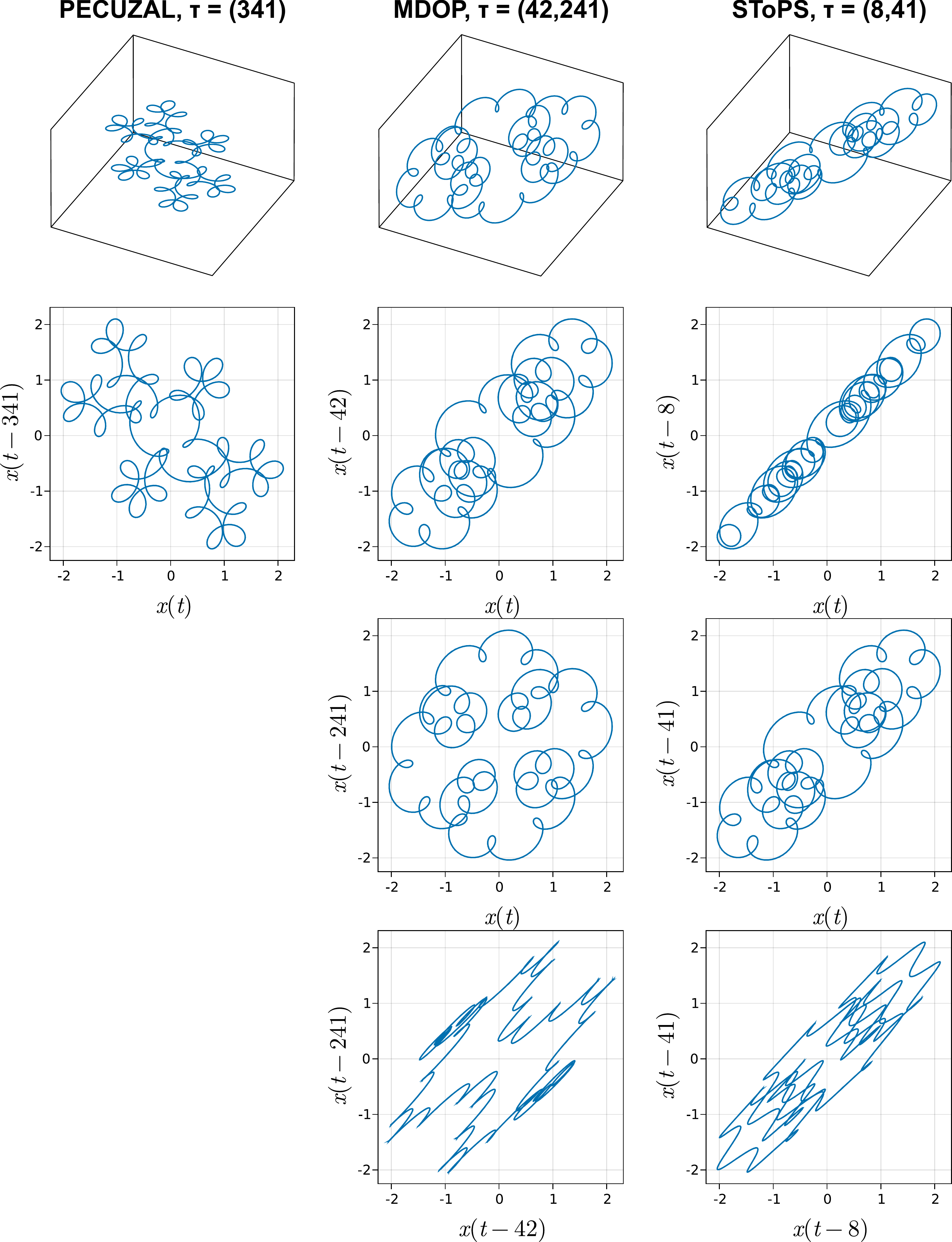}
    \caption{\textcolor{edit_col}{Sum of sines ($\omega = 1,5,30)$ phase space reconstructions with PECUZAL, MDOP and SToPS.}}
    \label{fig:SINES_Embedding}
\end{figure*}

\begin{figure*}
    \centering
    \includegraphics[width = 0.9\textwidth]{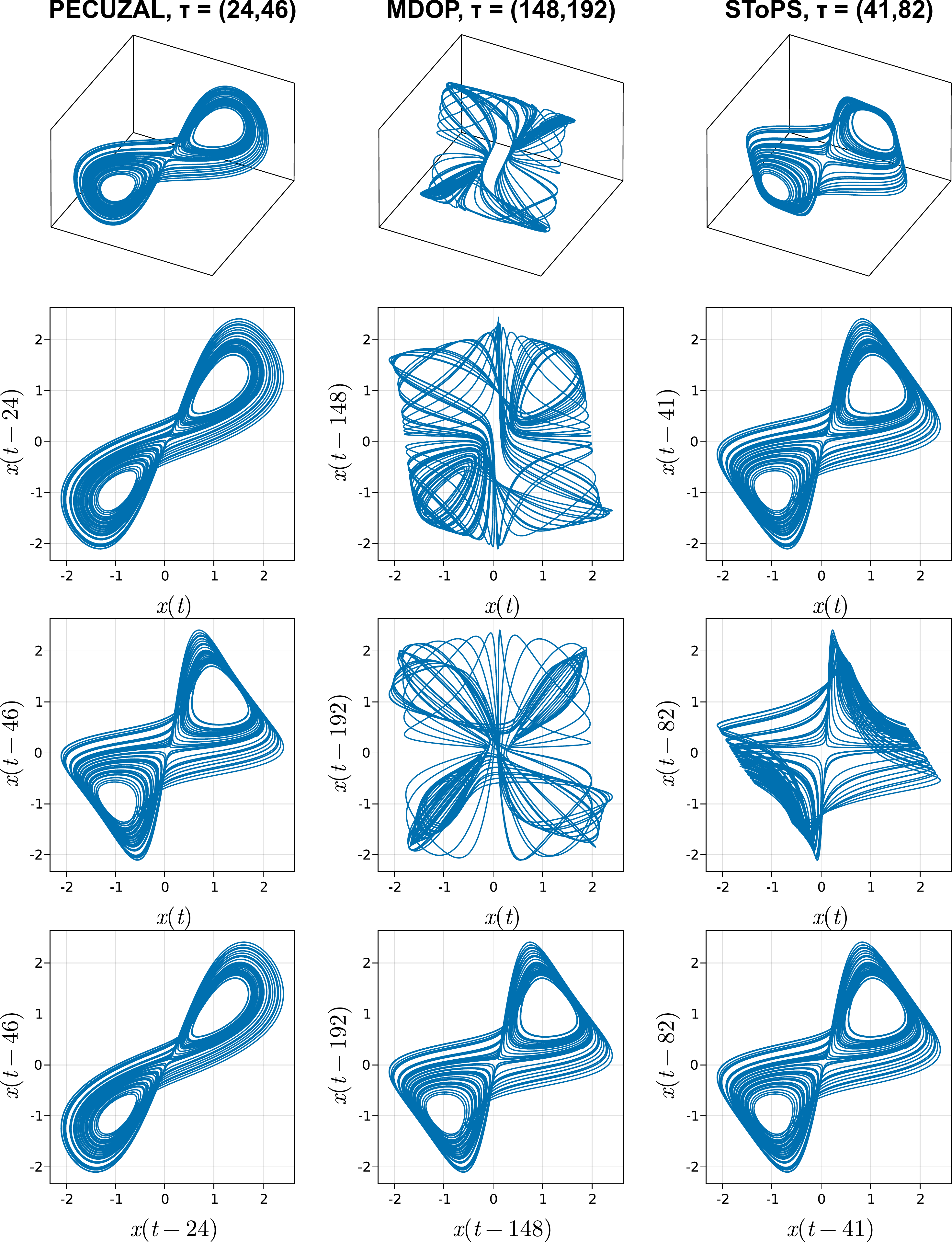}
    \caption{\textcolor{edit_col}{Lorenz phase space reconstructions with PECUZAL, MDOP and SToPS.}}
    \label{fig:LORENZ_Embedding}
\end{figure*}

\textcolor{edit_col}{A comparison of the various reconstructed attractors for the sum of sines and Lorenz time series is provided in Figures \ref{fig:SINES_Embedding} and \ref{fig:LORENZ_Embedding} with the first two detected lags are given. For the periodic sum of sines, PECUZAL yielded only a single lag. For the Lorenz time series, SToPs only yielded a single peak at $\tau = 41$ which was subsequently used for uniform embedding. Both PECUZAL and SToPS share a similar lag at $\tau \approx 40$. MDOP yields lags that cause overembedding.}

\begin{figure*}
    \centering
    \includegraphics[width = 0.9\textwidth]{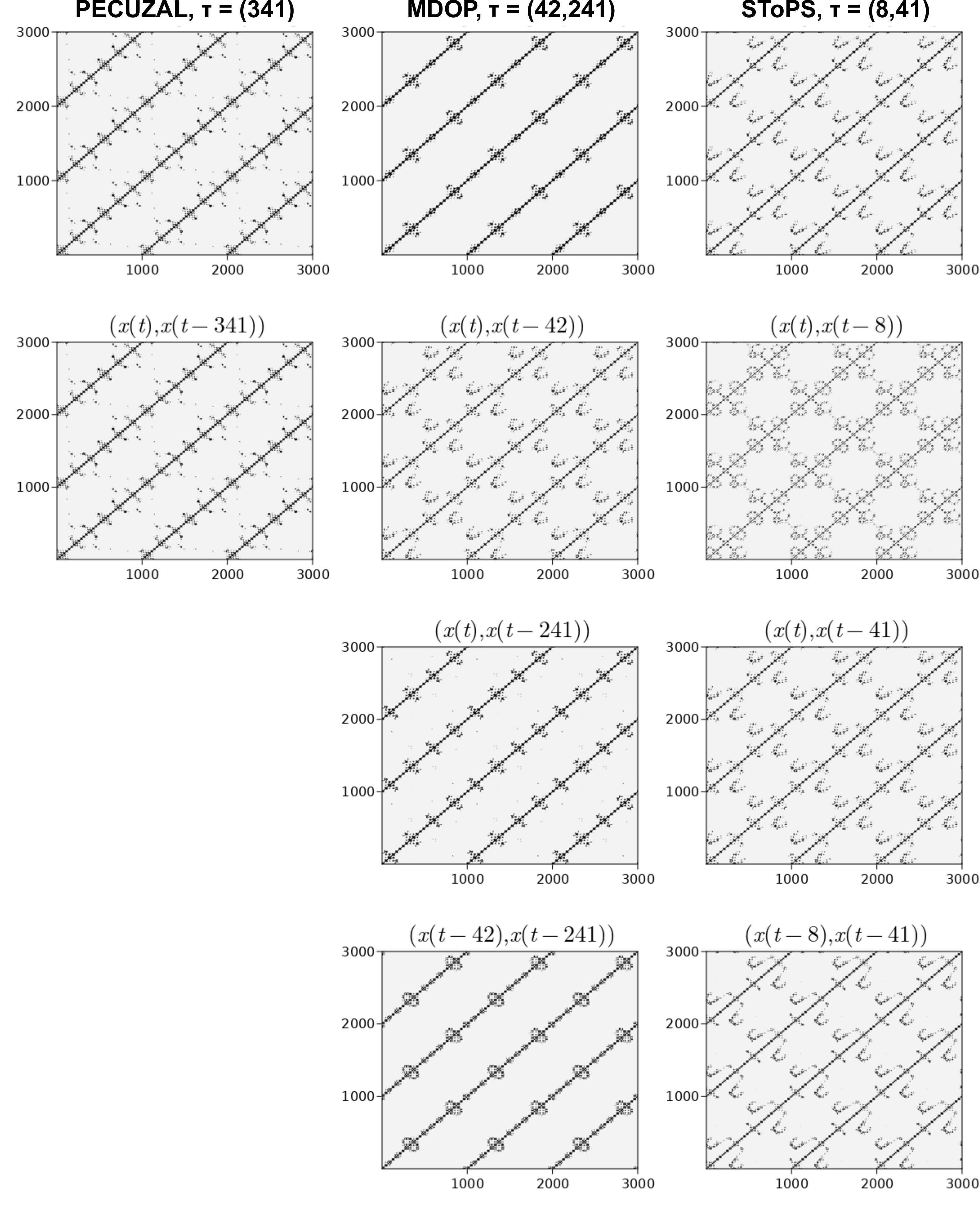}
    \caption{\textcolor{edit_col}{Sum of sines ($\omega = 1,5,30)$ recurrence plots of 2D and 3D embeddings.}}
    \label{fig:SINES_RP}
\end{figure*}

\begin{figure*}
    \centering
    \includegraphics[width = 0.9\textwidth]{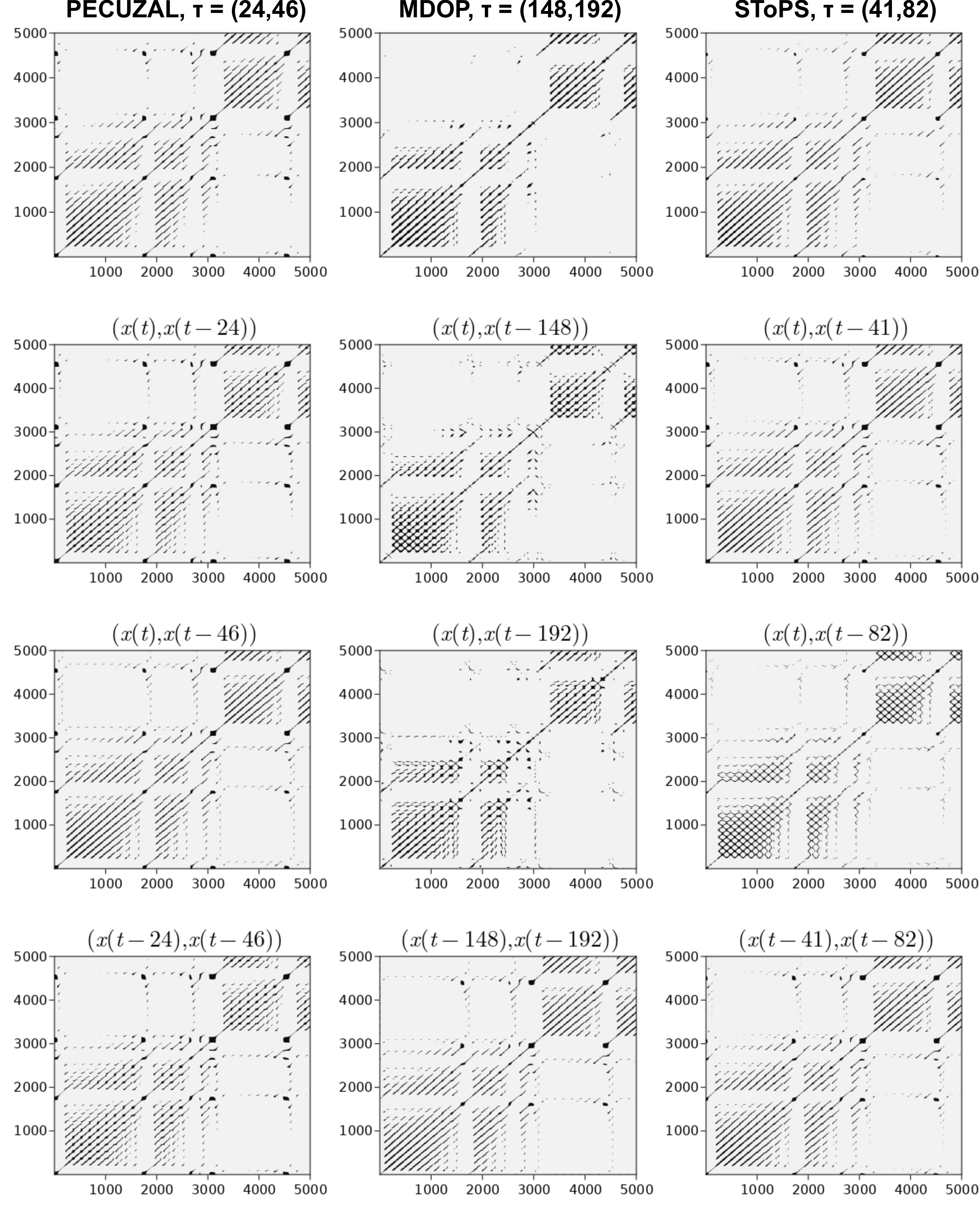}
    \caption{\textcolor{edit_col}{Lorenz recurrence plots of 2D and 3D embeddings.}}
    \label{fig:LORENZ_RP}
\end{figure*}

\begin{figure*}
    \centering
    \includegraphics[width = 0.9\textwidth]{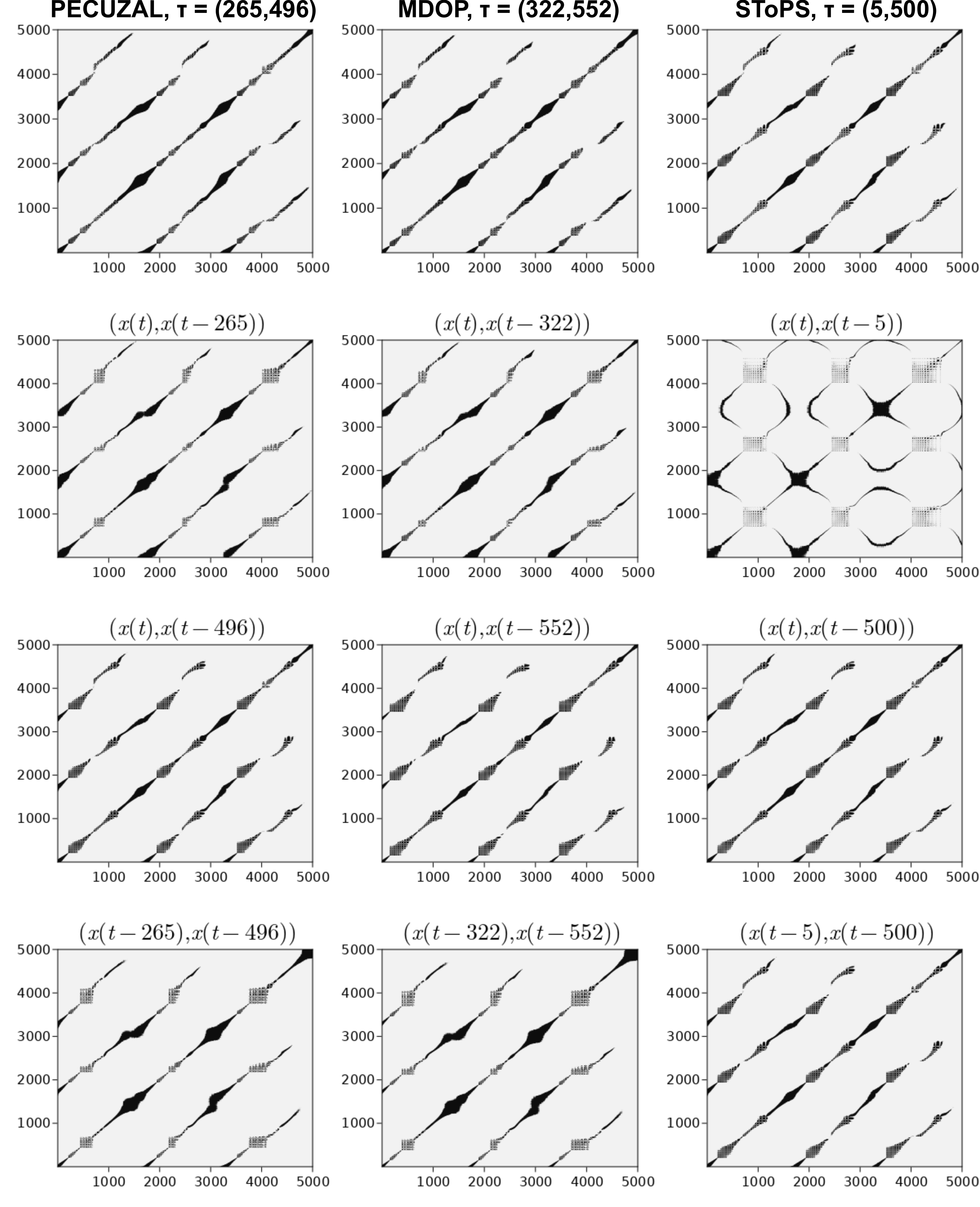}
    \caption{\textcolor{edit_col}{Lobster LP neuron recurrence plots of 2D and 3D embeddings.}}
    \label{fig:LOBSTER_RP}
\end{figure*}

\textcolor{edit_col}{Recurrence plots for the corresponding 2D and 3D embeddings for the sum of sines, Lorenz and lobster LP neuron time series are provided in Figures \ref{fig:SINES_RP}, \ref{fig:LORENZ_RP} and \ref{fig:LOBSTER_RP}. All embeddinig methods were able to preserve some of the periodic structure for the sum of sines data set. However, the recurrence plots for MDO and SToPS revealed more small scale structure than PECUZAL. For Lorenz, both PECUZAL and SToPS yield similar recurrence plots. MDOP recurrence plot lose some detail in comparison and is likely due to overfolding of the attractor caused by large embedding lags. In the lobster LP neuron time series, the selection of a small lag with SToPS reveals the expected periodic behaviour in several diagonal regions of the recurrence plot. This is not as clear in the PECUZAL and MDOP embeddings where the detected periodic behavour is dominated by the slow neuron dynamics.}

\end{document}